\documentclass[twocolumn, amsthm]{autart}    %

\usepackage{cite}
\usepackage{graphicx}
\usepackage{textcomp}
\usepackage{subcaption}
\usepackage{mathrsfs}
\usepackage{multirow}
\usepackage[flushleft]{threeparttable}
\usepackage{textgreek}
\usepackage[scr=boondox]{mathalpha}
\usepackage{array}

\usepackage{mystyle}
\usepackage{list_macro}
\usepackage{labelrel}

\newlist{lemitem}{enumerate}{10}
\setlist[lemitem]{label=(\alph*)}
\crefname{lemitemi}{}{lemitem}
\crefname{equation}{}{}
\crefname{thm}{theorem}{theorems}
\Crefname{thm}{Theorem}{Theorems}
\crefname{lem}{lemma}{lemmas}
\Crefname{lem}{Lemma}{Lemmas}
\crefname{def}{definition}{definitions}
\Crefname{def}{Definition}{Definitions}

\setboolean{showboth}{false}  %
\setboolean{showsecondboth}{false}

\begin{document}

\setlength{\aboverulesep}{-0.1pt}
\setlength{\belowrulesep}{0pt}
\setlength{\extrarowheight}{-.5ex}

\begin{frontmatter}

\title{\EMpp: A parameter learning framework for stochastic switching systems}

\thanks[footnoteinfo]{Work supported by
the Research Foundation Flanders (FWO) research projects \revnew{2}{postdoctoral grant 12Y7622N and}{}G081222N, G033822N, and G0A0920N; 
Research Council KU Leuven C1 project No. C14/24/103; 
European Union’s Horizon 2020 research and innovation programme under the Marie Skłodowska-Curie grant agreement No. 953348.
\revnew{2}{}{The work of Mathijs Schuurmans was supported by an FWO postdoctoral fellowship under Project 1248626N.}
}

\author[KUL]{Renzi Wang\corauthref{cor}}\revnew{2}{\ead{\href{mailto:renzi.wang@kuleuven.be}{renzi.wang@kuleuven.be}}}{},
\author[KUL]{Alexander Bodard}   \revnew{2}{\ead{\href{mailto:alexander.bodard@kuleuven.be}{alexander.bodard@kuleuven.be}}}{},  %
\author[KUL]{Mathijs Schuurmans} \revnew{2}{\ead{\href{mailto:mathijs.schuurmans@kuleuven.be}{mathijs.schuurmans@kuleuven.be}}}{},               %
\author[KUL]{Panagiotis Patrinos}\revnew{2}{\ead{\href{mailto:panos.patrinos@kuleuven.be}{panos.patrinos@kuleuven.be}}}{}

\corauth[cor]{Corresponding author. Tel: +32 16 32 03 64. Email: \href{mailto:renzi.wang@kuleuven.be}{renzi.wang@kuleuven.be}}

\address[KUL]{KU Leuven, Department of Electrical Engineering \textsc{esat-stadius} -- %
Kasteelpark Arenberg 10, B-3001 Leuven, Belgium}  %


\begin{abstract}                          %
This paper proposes a general switching dynamical system model, 
and a custom majorization-minimization-based algorithm \EMpp for identifying its parameters.
For certain families of distributions, such as Gaussian distributions, this algorithm reduces to the well-known expectation-maximization method.
We prove global convergence of the algorithm under suitable assumptions,
thus addressing an important open issue in the switching system identification literature.
The effectiveness of both the proposed model and algorithm is validated through extensive numerical experiments.

\end{abstract}

\end{frontmatter}
\endNoHyper %

\section{Introduction}

Obtaining a realistic model is of crucial importance in any model-based control application.
Yet, as the complexity of the underlying systems keeps on increasing, deriving such models from first principles becomes ever more difficult.
This evolution has motivated the development of data-driven modelling methods. 
Black-box models like neural networks can capture complex dynamics, but are challenging to analyze.
By contrast, \emph{switching systems} approximate complicated nonlinear systems as a combination of simple systems,
and as such offer a good balance between simplicity and expressiveness
\cite{paoletti2007identification,moradvandi2023models}.

Switching dynamical systems are dynamical systems which consist of a set of subsystems, 
and some switching mechanism that governs which of these subsystems is active at a given time.
The switching mechanism may be either \emph{subsystem-state-dependent}, or \emph{-independent} \cite{moradvandi2023models}.
A state-dependent switching mechanism often utilizes polyhedral partitioning,
in which case the activation of a subsystem depends on whether a regressor vector lies within a given region, either deterministically \cite{breschi2016piecewise,du2021online,bemporad2022piecewise}
or probabilistically \cite{jordan1994hierarchical,lecun2006tutorial}.
On the other hand, Markov jump systems are popular examples of state-independent switching systems.
The switching probability of these models is determined by a fixed transition matrix \cite{costa2005discrete,jin2012identification,fan2017robust}.
\revnew{2}{Moreover,}{}%
\rev{\cite{bemporad2018fitting,leoni2024explainableArxiv}}{\cite{bemporad2018fitting,leoni2025explainable}} 
propose general frameworks covering both \revnew{2}{types of switching}{}mechanisms,
which iteratively identify the subsystem parameters and the most probable subsystem combination weights for each data point.
However,
\revnew{2}{during}{in} inference
these models repeatedly solve an optimization problem to determine the weights for each predicted data point,
potentially limiting their integration with optimization-based controllers.

Although it is convenient to model subsystems as linear systems with additive Gaussian noise \cite{jin2012identification,kroemer2014learning,bianchi2020alternating},
this choice may not always be suitable for the data at hand.
For instance, hidden Markov models \cite{rabiner1989tutorial} typically involve observations from categorical distributions.
The exponential family \cite{wainwright2008graphical,efron2022exponential} is a general distribution class that includes both Gaussian and categorical distributions,
making it more versatile for modelling diverse types of data.
However, a limitation of the exponential family is its inability to 
capture many commonly occurring heavy-tailed distributions \cite{kodamana2018approaches,aravkin2017generalized}.
Such distributions offer an alternative option for modelling subsystems, 
and provide more robustness against outliers \cite{liu2020variational,fan2017robust}.

Algorithms for identifying switching systems are typically optimization-based \cite{moradvandi2023models}, and consist of two primary variants:
coordinate-descent methods 
\rev{\cite{bemporad2022piecewise,leoni2024explainableArxiv,bemporad2018fitting,bianchi2020alternating}}{\cite{bemporad2022piecewise,leoni2025explainable,bemporad2018fitting,bianchi2020alternating}}
and
expectation-maximization methods \cite{jin2012identification,fan2017robust,kroemer2014learning, sammaknejad2019review}.
Despite achieving good experimental performance, most methods lack convergence guarantees. 
We remark that some optimality guarantees for coordinate-descent-based methods have been provided \rev{\cite{bemporad2022piecewise,leoni2024explainableArxiv}}{\cite{bemporad2022piecewise,leoni2025explainable}}, but 
these do not ensure convergence to stationary points.

In this work, we present an identification method for generalized switching systems based on the majorization-minimization (MM) principle \cite{hunter2004tutorial,lange2016mm}.
This principle includes the EM algorithm as a special case.
The subsequential convergence of MM schemes has been proven under generic conditions
\cite[\S 7.3]{lange2016mm} \cite{razaviyayn2013unified}.
However, verifying whether these hold is usually 
only straightforward when the involved functions are smooth.
In such smooth case,
\cite{kunstner2021homeomorphic} connects the EM algorithm with the mirror descent algorithm \cite{beck2003mirror,bolte2018first},
which is similar to a classical gradient descent method, 
but replaces the Euclidean distance term by a so-called \emph{Bregman distance} \cite{bregman1967relaxation}.
In a similar way, the (subsequential) convergence of the EM algorithm has been analyzed in
\cite{tseng2004analysis,chretien2008algorithms,kunstner2021homeomorphic} by interpreting EM as a non-Euclidean descent method.
\revnew[tensor-based]{2}{}{%
The pursuit of such guarantees is a central theme in learning latent variable models,
where tensor-based decomposition method \cite{anandkumar2014tensor} offers an alternative path to global convergence.}
\revnew{2}{%
However, these works focus on identifying a static model with i.i.d. data.}{%
However, these works primarily operate in static settings by identifying the model from i.i.d. data.}
\revnew{2}{}{Furthermore, while tensor methods provide global results, they often rely on initialization conditions that are difficult to verify in practice,
and do not naturally account for the time-varying switching dynamics of the systems considered here.}
For a more general problem setup,
\cite{wu1983convergence,sun2016majorization} also prove sequential convergence of related methods, 
but under a discrete set assumption that may be difficult to verify.

Our contributions can be summarized as follows.
\vspace{-0.3\baselineskip}
\begin{enumerate}[wide]
    \item We present a \emph{general switching system model} that generalizes the one proposed in \cite{kroemer2014learning} to more \revnew{2}{generic}{}subsystem dynamics. The model includes:
    \begin{inlinelist}
        \item A parametric formulation for the switching mechanism that can be either subsystem-state-dependent or -independent. 
            In contrast to \rev{\cite{bemporad2018fitting,leoni2024explainableArxiv}}{\cite{bemporad2018fitting,leoni2025explainable}}, 
            this parametric formulation only requires the evaluation of the identified model during inference.
        \item A generic subsystem dynamics formulation covering various popular distributions, including the exponential family, 
        as well as a number of heavy-tailed distributions. 
    \end{inlinelist}
    \item\label{item: contribution_convergence} We introduce an \emph{MM-based identification algorithm} \EMpp for the proposed model and establish its global subsequential convergence to stationary points,
    even when the probability densities include nonsmooth functions.
    Our method generalizes the EM method to identify a significantly wider class of switching systems.
    Additionally, we extend the convergence analysis of the EM algorithm beyond the canonical exponential family and i.i.d. data presented in \cite{kunstner2021homeomorphic}. 
    We prove full sequence convergence to a stationary point under a mild Kurdyka-\L ojasiewicz (KL) condition \cite{lojasiewicz1963propriete}, which improves upon existing results even in the EM setting.
    \item We confirm the expressiveness of the proposed model and the effectiveness
        of the identification method through a series of numerical experiments.
\end{enumerate}

\textbf{Overview}.
\Cref{sec:problem-statement} introduces a general switching system and the corresponding identification problem.
\Cref{sec: mm} presents the proposed method for solving this problem, and
\cref{sec: conv_analysis} analyzes its convergence.
\Cref{sec: e_step} \rev{details how to effectively evaluate the key functions in the problem}{discusses the efficient implementation of our method}.
\Cref{sec: numerical_experiment} experimentally shows the efficiency of both the model and algorithm.

\textbf{Notation}.
Denote the cone of $m \times m$ positive definite matrices by $\mathbb{S}_{++}^{m}$.
Let $\norm{A}^2_B \!\dfn\! \tr(A^\top B A)$ for $B \in \mathbb{S}_{++}^m$.
Let $\N_{[a, b]} \!=\! \N \cap [a, b]$ and $\mathbf{1}_d \!=\! \bsmat{1 & \dots & 1} \in \re^d$.
Define $\lse: \R^{d} \to \R$ as
\(
    \lse(x) \! =\! \ln \big(
        \textstyle{\sum_{j=1}^{d}} \exp(x_j) 
    \big),
\)
and the softmax function $x \mapsto \sigma(x)$ 
with  
\(
    \sigma_i(x) \! = \!\exp\big(x_i\big) / \sum_{j=1}^{d}\exp\big(x_j\big)
\)
the $i$th entry of the output $\sigma(x)$.
Denote the class of $n$-times continuously differentiable functions by $\mathcal{C}^n$.
Given a function $f$: $\mathcal{O} \subseteq \R^n \to \R$, we use the convention $\forall x \not \in \mathcal{O}$: $f(x) = +\infty$, and follow the definition of strict continuity of \cite[Definition 9.1]{rockafellar2009variational}, which is equivalent to \emph{local} Lipschitz-continuity.
The directional derivative of $f$ at $x \in \mathcal{O}$ along \rev[typoNotationDirection]{$d$}{$v$} $\in \R^n$ is defined as 
\(
    f'(x\midsc v) = \lim_{\lambda \searrow 0}\frac{f(x + \lambda v) - f(x)}{\lambda}.    
\)
A function $f \in \mathcal{C}^1$ is Lipschitz smooth if $f$ has Lipschitz continuous gradients.

\section{Problem statement}\label{sec:problem-statement}
\newlength{\rowspace}
\setlength{\rowspace}{0.3pt}
\newlength{\columnspace}
\setlength{\columnspace}{1.8pt}

\begin{table*}[tb]
    \centering 
    \caption{List of pdf and pmf $p(y \mid z, \xi = j \midsc \dBeta)$ satisfying \eqref{eq: subsys_model}. The first block are members of canonical exponential family \cite{wainwright2008graphical,efron2022exponential}.}
    \label{tab: example_distr}
\revnew{2}{%
\begin{resetcolor}
  
\rev[reviewer3point1Table]{
    \begin{threeparttable}
    \addtolength{\tabcolsep}{-\columnspace}  
    \begin{tabular}{l l ccccc}
        \toprule \\ [-10pt]
                                                                  & Notation                                                         & $\beta_j$        & C                       & $f(x)$                & \rev{$\ell(y, z, \beta_j)$}{}                                     & \rev{$g(y, z, \beta_j)$}{} \\
        \midrule                                                                                                                                                                                     \\[-8pt]
        Exp. family\tnote{a}                                      &  $\quad\quad/$                                                   & $\beta$        & $\eta(y, z)$              & \rev{$x$}{}                   & $-\inner{\beta}{\sstat(y, z)}$                          & $\logpar(z, \beta)$    \\
        Categorical                                               & $\mathrm{\mathsf{Cat}}(\sigma(\Theta^\top z))$                   & $\Theta$       & 1                         & \rev{$x$}{}                   & $-\Theta_{y}^\top z$                                    & $\lse(\Theta^\top z)$\\ [\rowspace]
        Gaussian                                                  & $\gauss(\mu \!=\! Lz, \Sigma)$                                   & \tiny{$\begin{aligned} B =& \Lambda L\\[-3pt] \Lambda =& \Sigma^{-1}   \end{aligned}$}   & $(2\pi)^{-n_y / 2}$     & \rev{$x$}{}                   & $\frac{1}{2}\!\inner{\Lambda}{yy^\top\!}{-}\inner{B}{zy^\top\!}$  & $\frac{1}{2} \!\big(\!\norm{Bz}^2_{\Lambda^{{-}1}}\!{-}\ln\!\det(\!\Lambda\!)\!\big)$\\ 
        [\rowspace] \midrule \\ [-10pt]
        Student's t\rev{\tnote{b}}{\tnote{a}} \cite{roth2012multivariate}          & $\mathrm{\mathsf{St}}_{\nu}(\mu = Lz, \Sigma)$                   & \tiny{$\begin{aligned} B =& \Lambda L\\[-3pt] \Lambda =& \Sigma^{-1}\end{aligned}$}   & \scriptsize$\tfrac{\Gamma((\nu+n_y) / 2)}{\sqrt{(\pi\nu)^{n_y}}\Gamma(\nu / 2)}$ & \rev{$  \frac{\nu {+} n_y}{2} \! \ln (\!1 \!+\! \frac{2x}{\nu}\!)$}{} & $\frac{1}{2}\norm{\Lambda y - B z}^{2}_{\Lambda^{-1}}$ & $-\frac{1}{2} \ln\det(\Lambda)$ \\ [\rowspace]
        $\ell_1$-Laplace \cite{aravkin2011}                       & $\mathrm{\mathsf{Laplace}} (\mu \!=\! Lz,\! \Sigma)$             & \tiny{$\begin{aligned} M =& R L      \\[-3pt] R       =& \Sigma^{-\nicefrac{1}{2}}\end{aligned}$}        & $2^{-n_y / 2}$          & \rev{$\sqrt{2}x$}{}           & $\norm{R y - M z}_1$                                    & $-\sum_{i = 1}^{n_y}\ln (R_{ii})$\\ [\rowspace]
        Logistic \cite[\S 23]{johnson1995continuous}              & $\mathrm{\mathsf{Logistic}}(\mu \!=\! az,\! \varsigma)$          & \tiny{$\begin{aligned} b =& \lambda a\\[-3pt] \lambda =& \varsigma^{-1}\end{aligned}$}   & $\frac{1}{4}$           & \rev{$2\ln (x)$}{}            & $\mathrm{cosh} (\frac{1}{2}(\lambda y - b z))$          & $-\ln \lambda$\\ [\rowspace]
        Gumbel \cite[\S 22]{johnson1995continuous}                & $\mathrm{\mathsf{Gumbel}}  (\mu \!=\! az,\! \varsigma)$          & \tiny{$\begin{aligned} b =& \lambda a\\[-3pt] \lambda =& \varsigma^{-1}\end{aligned}$}   & $1$                     & \rev{$x$}{}                   & $\exp(-\lambda y + b z)$                                & $-\ln \lambda + \lambda y - b z$ \\ [\rowspace]
        \bottomrule
    \end{tabular}
    \addtolength{\tabcolsep}{\columnspace}  
    \begin{tablenotes}
        \item[a] Functions $\eta: \re^{n_y} \times \re^{n_z} \rightarrow \re_{++}$,
                 $\sstat: \re^{n_y} \times \re^{n_z} \rightarrow \mathcal{B}$,
                 $\logpar: \re^{n_z} \times \mathcal{X}_\beta \rightarrow \R$.
        \item[b] With given degrees of freedom $\nu$.
    \end{tablenotes}
\end{threeparttable}
}{
        \begin{threeparttable}
        \addtolength{\tabcolsep}{-\columnspace} 
        \begin{tabular}{l l ccccc}
            \toprule \\ [-10pt]
                                                                      & Notation                                                         & $\beta_j$        & C                       & $f(x)$                & \rev{}{$g_1(y, z, \beta_j)$}                                     & \rev{}{$g_2(y, z, \beta_j)$} \\
            \midrule                                                                                                                                                                                     \\[-8pt]
            Exp. family\rev{\tnote{a}}{}                              &  $\quad\quad/$                                                   & $\beta$        & $\eta(y, z)$            & \rev{}{$x_1 + x_2$}                   & $-\inner{\beta}{\sstat(y, z)}$                          & $\logpar(z, \beta)$    \\
            Categorical                                               & $\mathrm{\mathsf{Cat}}(\sigma(\Theta^\top z))$                   & $\Theta$       & 1                       & \rev{}{$x_1 + x_2$}                   & $-\Theta_{y}^\top z$                                    & $\lse(\Theta^\top z)$\\ [\rowspace]
            Gaussian                                                  & $\gauss(\mu \!=\! Lz, \Sigma)$                                   & \tiny{$\begin{aligned} B =& \Lambda L\\[-3pt] \Lambda =& \Sigma^{-1}   \end{aligned}$}   & $(2\pi)^{-n_y / 2}$     & \rev{}{$x_1 + x_2$}                   & $\frac{1}{2}\!\inner{\Lambda}{yy^{\!\top\!}\!}{-}\inner{B\!}{zy^{\!\top\!}\!}$  & $\frac{1}{2} \!\big(\!\norm{Bz}^2_{\Lambda^{\!{-}\!1}}\!{-}\!\ln\!\det(\!\Lambda\!)\!\big)$\\ 
            [\rowspace] \midrule \\ [-10pt]
            Student's t\rev{\tnote{b}}{\tnote{a}} \cite{roth2012multivariate}          & $\mathrm{\mathsf{St}}_{\nu}(\mu = Lz, \Sigma)$                   & \tiny{$\begin{aligned} B =& \Lambda L\\[-3pt] \Lambda =& \Sigma^{-1}\end{aligned}$}   & \scriptsize$\tfrac{\Gamma((\nu+n_y) / 2)}{\sqrt{\!(\pi\nu)^{n_y}}\!\Gamma(\nu / 2)}$ & \rev{}{$\frac{\nu {+} n_y}{2} \! \ln (\!1 \!{+}\! \frac{2x_1}{\nu}\!)\!{+}\!x_2$} & $\frac{1}{2}\norm{\Lambda y - B z}^{2}_{\Lambda^{-1}}$ & $-\frac{1}{2} \ln\det(\Lambda)$ \\ [\rowspace]
            $\ell_1$-Laplace \cite{aravkin2011}                       & $\mathrm{\mathsf{Laplace}} (\mu \!=\! Lz,\! \Sigma)$             & \tiny{$\begin{aligned} M {=}& R L      \\[-3pt] R       =& \Sigma^{-\nicefrac{1}{2}}\end{aligned}$}      & $2^{-n_y / 2}$          & \rev{}{$\sqrt{2}x_1 + x_2 $}           & $\norm{R y - M z}_1$                                    & $-\sum_{i = 1}^{n_y}\ln (R_{ii})$\\ [\rowspace]
            Logistic \cite[\S 23]{johnson1995continuous}              & $\mathrm{\mathsf{Logistic}}(\mu \!=\! az,\! \varsigma)$          & \tiny{$\begin{aligned} b =& \lambda a\\[-3pt] \lambda =& \varsigma^{-1}\end{aligned}$}                 & $\frac{1}{4}$           & \rev{}{$2\ln (x_1) + x_2$}            & $\mathrm{cosh} (\frac{1}{2}(\lambda y - b z))$          & $-\ln \lambda$\\ [\rowspace]
            Gumbel \cite[\S 22]{johnson1995continuous}                & $\mathrm{\mathsf{Gumbel}}  (\mu \!=\! az,\! \varsigma)$          & \tiny{$\begin{aligned} b =& \lambda a\\[-3pt] \lambda =& \varsigma^{-1}\end{aligned}$}                 & $1$                     & \rev{}{$x_1+x_2$}                   & $\exp(-\lambda y + b z)$                                & $-\ln \lambda + \lambda y - b z$ \\ [\rowspace]
            \bottomrule
        \end{tabular}
        \addtolength{\tabcolsep}{\columnspace}  
        \begin{tablenotes}
            \rev{\item[a] Functions $\eta: \re^{n_y} \times \re^{n_z} \rightarrow \re_{++}$,
                     $\sstat: \re^{n_y} \times \re^{n_z} \rightarrow \mathcal{B}$,
                     $\logpar: \re^{n_z} \times \mathcal{X}_\beta \rightarrow \R$.}{}
            \rev{\item[b]}{}\item[a] With given degrees of freedom \rev[rangeDoF]{$\nu$}{$\nu > 0$}.
        \end{tablenotes}
    \end{threeparttable}
}
\end{resetcolor}
}{}

\begin{threeparttable}
    
    \setlength{\tabcolsep}{7.5pt}      
    \renewcommand{\arraystretch}{1.65}
    \setlength\lightrulewidth{0.01pt}
    
    \begin{tabular}{l l c l >{\tiny}c >{\scriptsize}l >{\scriptsize}l}
        \toprule[1pt] 
        Dist. / Notation\revnew{2}{}{\tnote{a}} & $\beta_j$ & $C$ & $f(x)$ & \revnew{2}{}{$i$} & \revnew{2}{}{$g_i(\cdot)$} & \revnew{2}{}{$T_i(y, z, \beta)$} \\
        \midrule[1pt]
        
        \multirow{2}{*}{\shortstack[l]{Exp. family\revnew{2}{}{\tnote{b}}}} 
        & \multirow{2}{*}{$\beta$} 
        & \multirow{2}{*}{$\eta(y, z)$}
        & \multirow{2}{*}{$\revnew{2}{}{\sum_i x_i}$} 
          & \revnew{2}{}{1} & \revnew{2}{}{$-u$} & \revnew{2}{}{$\inner{\beta}{\sstat(y,\! z)}$} \\
        & & & & \revnew{2}{}{2} & \revnew{2}{}{$\logintexp^{u(y)}\diff y$} & \revnew{2}{}{$\inner{\beta}{\sstat(y,\!z)} + \ln\eta(y,\! z)$} \\
        \arrayrulecolor{black!30}\midrule
        
        \multirow{2}{*}{\shortstack[l]{Categorical\revnew{2}{}{\tnote{c}} \\ $\mathrm{\mathsf{Cat}}(\sigma(\Theta^\top z))$}} 
        & \multirow{2}{*}{$\Theta$} 
        & \multirow{2}{*}{$1$}
        & \multirow{2}{*}{$\revnew{2}{}{\sum_i x_i}$} 
          & \revnew{2}{}{1} & \revnew{2}{}{$-\tr(u)$} & \revnew{2}{}{$\Theta^\top (z e_y^\top)$} \\
        & & & & \revnew{2}{}{2} & \revnew{2}{}{$\lse(u)$} & \revnew{2}{}{$\Theta^\top z$} \\
        \arrayrulecolor{black!30}\midrule

        \multirow{3}{*}{\shortstack[l]{Gaussian \\ $\gauss(\revnew{2}{\mu=}{}Lz, \Sigma)$}} 
        & \multirow{3}{*}{$\begin{aligned} B =& \Lambda L\\[-2pt] \Lambda =& \Sigma^{-1}   \end{aligned}$}
        & \multirow{3}{*}{$(2\pi)^{\frac{-n_y}{2}}$}
        & \multirow{3}{*}{$\revnew{2}{}{\sum_i x_i}$} 
          & \revnew{2}{}{1} & \revnew{2}{}{$\text{tr}(u)$} & \revnew{2}{}{$\frac{1}{2}\Lambda (yy^\top) - B(zy^\top)$} \\
        & & & & \revnew{2}{}{2} & \revnew{2}{}{$\frac{1}{2} \norm{u}_{V^{-1}}^2$} & \revnew{2}{}{$(Bz, \Lambda)$} \\
        & & & & \revnew{2}{}{3} & \revnew{2}{}{$-\frac{1}{2}\ln \det(u)$} & \revnew{2}{}{$\Lambda$} \\
        \arrayrulecolor{black}\midrule[1pt]
        
        \multirow{2}{*}{\shortstack[l]{Student's t\tnote{d} \cite{roth2012multivariate} \\ $\mathrm{\mathsf{St}}_{\nu}(\revnew{2}{\mu=}{}Lz, \Sigma)$}} 
        & \multirow{2}{*}{$\begin{aligned} B =& \Lambda L\\[-2pt] \Lambda =& \Sigma^{-1}\end{aligned}$} 
        & \multirow{2}{*}{\scriptsize $\frac{\Gamma(\frac{\nu+n_y}{2})}{\sqrt{(\pi\nu)^{n_y}}\Gamma(\frac{\nu}{2})}$}
        & \multirow{2}{*}{{\tiny $\begin{aligned}\tfrac{\nu {+} n_y}{2} \ln (1 {+} \tfrac{2x_1}{\nu}) \\ {+}x_2\end{aligned}$}} 
          & \revnew{2}{}{1} & \revnew{2}{}{$\frac{1}{2} \norm{u}_{V^{-1}}^2$} & \revnew{2}{}{$(\Lambda y - Bz, \Lambda)$} \\
        & & & & \revnew{2}{}{2} & \revnew{2}{}{$-\frac{1}{2} \ln\det(u)$} & \revnew{2}{}{$\Lambda$} \\
        \arrayrulecolor{black!30}\midrule
        
        \multirow{2}{*}{\shortstack[l]{$\ell_1$-Laplace \cite{aravkin2011} \\ $\mathrm{\mathsf{Laplace}}(\revnew{2}{\mu=}{}Lz, \Sigma)$}} 
        & \multirow{2}{*}{$\begin{aligned} M {=}& R L      \\[-2pt] R       =& \Sigma^{-\nicefrac{1}{2}}\end{aligned}$} 
        & \multirow{2}{*}{$2^{-\frac{n_y}{2}}$}
        & \multirow{2}{*}{$\sqrt{2}x_1 {+} x_2$} 
          & \revnew{2}{}{1} & \revnew{2}{}{$\norm{u}_1$} & \revnew{2}{}{$Ry - Mz$} \\
        & & & & \revnew{2}{}{2} & \revnew{2}{}{$-\sum_{k=1}^{n_y} \ln(u_{kk})$} & \revnew{2}{}{$R$} \\
        \arrayrulecolor{black!30}\midrule
        
        \multirow{2}{*}{\shortstack[l]{Logistic \cite[\S 23]{johnson1995continuous} \\ $\mathrm{\mathsf{Logistic}}(\revnew{2}{\mu=}{}az, \varsigma)$}} 
        & \multirow{2}{*}{$\begin{aligned} b =& \lambda a\\[-2pt] \lambda =& \varsigma^{-1}\end{aligned}$} 
        & \multirow{2}{*}{$\frac{1}{4}$}
        & \multirow{2}{*}{$2\ln x_1 {+} x_2$} 
          & \revnew{2}{}{1} & \revnew{2}{}{$\cosh(\tfrac{1}{2} u)$} & \revnew{2}{}{$\lambda y - bz$} \\
        & & & & \revnew{2}{}{2} & \revnew{2}{}{$-\ln(u)$} & \revnew{2}{}{$\lambda$} \\
        \arrayrulecolor{black!30}\midrule

        \multirow{2}{*}{\shortstack[l]{Gumbel \cite[\S 22]{johnson1995continuous}\\ $\mathrm{\mathsf{Gumbel}}(\revnew{2}{\mu=}{}az, \varsigma)$}} 
        & \multirow{2}{*}{$\begin{aligned} b =& \lambda a\\[-2pt] \lambda =& \varsigma^{-1}\end{aligned}$} 
        & \multirow{2}{*}{$1$}
        & \multirow{2}{*}{$\revnew{2}{}{\sum_i x_i}$} 
          & \revnew{2}{}{1} & \revnew{2}{}{$\exp(u)-u$} & \revnew{2}{}{$-\lambda y + bz$} \\
        & & & & \revnew{2}{}{2} & \revnew{2}{}{$-\ln(u)$} & \revnew{2}{}{$\lambda$} \\
        \arrayrulecolor{black}\bottomrule[1pt]
    \end{tabular}
    \begin{tablenotes}
        \item[a] \revnew{2}{}{%
        In this table, $g_i$ can be a single- or two-argument function. We denote $u$, and $V$ if applicable, as the input(s) to $g_i$.}
        \item[b] \revnew{2}{}{%
        In this row, $g_2$ follows the definition of log-partition function $\mathcal{A}(z, \beta) = \logintexp^{\inner{\beta}{\sstat(y, z)} + \ln\eta(y, z)} \diff y$}
        \item[c] \revnew{2}{}{%
        In this row, $e_y$ represents a one-hot vector with $1$ at the $y$-th entry and $0$ elsewhere.}
        \item[d] With given degrees of freedom $\nu > 0$.
    \end{tablenotes}
\end{threeparttable}
\end{table*}
We consider a stochastic dynamical system that switches among $d$ linear subsystems, and is characterized by
\vspace{-0.3\baselineskip}
\begin{subequations}\label{eq: discrete_latent}
    \vspace{-0.5\baselineskip}
    \begin{alignat}{2}
        \xi_{t+1} 
            &\sim 
                p(\xi_{t+1} \mid z_t, \xi_t = i \midsc \dTheta)  
             = 
            \sigma_{\xi_{t+1}}\big(\Theta_i^\top z_t\big),
            \label{eq: switch_model} \\
        y_{t+1} 
            &\sim \nonumber 
            p(y_{t+1} \mid z_t, \xi_{t+1} = j\midsc \dBeta) \\
            \revalign  \rev{
                = C\!\exp\!\big({-}f(\ell(y_{t+1}, z_t, \beta_j))\!-\! g(y_{t+1}, z_t, \beta_j)\big){,} 
            \nonumber}{}\revnewline 
            & \rev{}{
                = C\exp\big[{-}f\big(G(y_{t+1}, z_t, \beta_j)\big)\big]{,}
             \label{eq: subsys_model}
            }
    \end{alignat}
\end{subequations}
where at time $t \geq 0$ the random variables $\xi_{t} \in \Xi = \{1, \dots, d\}$ and $y_t \in \re^{n_y}$
are the active subsystem index and the observation, respectively.
Both $\xi_t$ and $y_t$ depend on the trajectory history
\vspace{-0.3\baselineskip}%
\begin{equation}\label{eq: def_history} 
    z_t = \Upsilon(y_{t}, \dots, y_{t - t_y + 1})
    \vspace{-0.5\baselineskip}
\end{equation}
where $\Upsilon: \re^{t_y n_y} \rightarrow \re^{n_z}$ is a known (non-)linear mapping with the window length $t_y > 0$.
The switching mechanism \eqref{eq: switch_model} is modelled with a softmax function composed with a mapping that is linear w.r.t the parameter $\dTheta \dfn \{\Theta_i\}_{i = 1}^d$ with $\Theta_i \in \re^{n_{z} \times d}$.
With $\dBeta \dfn \{\beta_i\}_{i = 1}^d \in \mathcal{B}^d$ where $\mathcal{B}$ is a Euclidean parameter space,
the subsystem dynamics \eqref{eq: subsys_model} is modelled by a distribution with a probability density function (pdf)
or probability mass function (pmf)
defined by a parameter-invariant term $C > 0$, 
\rev{and functions $f$, $\ell$, and $g$\,}{ a function $f$ and a mapping $G$}
under the following assumption.
\begin{assumption}[Subsystem model]\label{asp: subsys_model}
    Regarding the functions in \eqref{eq: subsys_model}, i.e.,
    $f: \mathcal{X}_f \rightarrow \re$, 
    \rev{
        $\ell: \re^{n_y} \times \re^{n_z} \times \mathcal{X}_{\beta} \rightarrow \mathcal{X}_f$, and 
        $g: \re^{n_y} \times \re^{n_z} \times \mathcal{X}_{\beta} \rightarrow \re$}{
        $G: \re^{n_y} \times \re^{n_z} \times \mathcal{X}_{\beta} \rightarrow \mathcal{X}_f$
    } 
    with $\mathcal{X}_{\beta} \subseteq \mathcal{B}$ and \rev{$\mathcal{X}_{f} \subseteq \R$}{$\mathcal{X}_{f} \subseteq \R^{\revnew{2}{M}{n_g}}$} convex open sets,
    we assume that
    \begin{lemitem}[wide]
        \item \label{asp: affine_structure}\revnew{2}{
        \begin{resetcolor}
            \rev{
                $\beta \mapsto g(y, z, \beta)$ and $\beta \mapsto \ell(y, z, \beta)$ are}{
                $\beta \mapsto g_i(y, z, \beta)$ is%
            }
        \end{resetcolor}
        }{$\beta \mapsto g_i\big(T_i(y, z, \beta)\big)$ is}
        the $i$-th entry of mapping $G$ 
        \revnew{2}{and is}{with $g_i$ being}
        strictly continuous and convex\revnew{2}{}{, and $T_i$ being an affine mapping of $\beta$}
        for all $y \in \R^{n_y}$, $z \in \re^{n_z}$;
        \revnew{2}{%
        \begin{resetcolor}
            \rev{}{, where $g_i$ is the $i$-th entry of mapping $G$};
        \end{resetcolor}
        }{}
        \item $f$ is continuously differentiable, concave, and strictly increasing \revnew{2}{}{with $\nabla f(x) >0$ holding component-wise}. 
        The \rev{derivative $f'(\cdot)$ is upper-bounded by a constant $\bar{u} \in (0, \infty)$}{gradient norm is bounded, i.e., $\norm{\nabla f(x)}_{\infty} {\leq} \bar{u}$ for all $x \in \mathcal{X}_f$, where $\bar{u} \in (0, \infty)$}.
    \end{lemitem}
\end{assumption}
\revnew{2}{
    \begin{resetcolor}
\rev[new_reviewer2point1]{}{
    \begin{remark}[Translation invariance]
        Softmax functions are translation invariant, 
        i.e., $\sigma(x) = \sigma(x - a\mathbf{1}_d)$ with $a \neq 0$.
        This propety raises a challenge since infinitely many parameter values yield identical distributions.
        To resolve this,
        we apply a change of variable $\tilde{x}_i = x_i - x_d$ for all $i \in \N_{[1, d]}$,
        which maintains the same output while reducing the dimensionality of the unknown variable,
        This transformation ensures that $\sigma(\tilde{x}) = \sigma(\bar{x})$ if and only if $\tilde{x} = \bar{x}$ with $\tilde{x}_d = \bar{x}_d = 0$.
        Consequently,
        the parameter $\Theta_i \in \re^{n_z \times (d-1)}$ 
        and \eqref{eq: switch_model} takes the form $\dTheta \mapsto \sigma\big(\bsmat{z_t^\top \Theta_i  & 0}^\top\big)$.
        For the notational simplicity,
        we present the result without change of variable,
        though the derivations remain equivalent.
    \end{remark}
}
    \end{resetcolor}
}{}
\vspace{-0.3\baselineskip}

\begin{remark}[Translation invariance]\label{rem: translation_invariance}
        Softmax functions are translation invariant, 
        i.e., $\sigma(x) = \sigma(x - a\mathbf{1}_d)$ with $a \neq 0$.
        This property raises a challenge since infinitely many parameter values yield identical distributions.
        To resolve this,
        \revnew{2}{}{we enforce $x_d = 0$,
        which is equivalent to applying}
        \revnew{2}{we apply} a change of variable $\tilde{x}_i = x_i - x_d$ for all $i \in \N_{[1, d]}$,
        which maintains the same output while reducing the dimensionality of the unknown variable.
        \revnew{2}{This transformation}{Incorporating this constraint} ensures that $\sigma(\tilde{x}) = \sigma(\bar{x})$ if and only if $\tilde{x} = \bar{x}$ with $\tilde{x}_d = \bar{x}_d = 0$.
        \revnew{2}{Consequently,
        the parameter $\Theta_i \in \re^{n_z \times (d-1)}$ 
        and \eqref{eq: switch_model} takes the form $\dTheta \mapsto \sigma\big(\bsmat{z_t^\top \Theta_i  & 0}^\top\big)$.
        For the notational simplicity,
        we present the result without change of variable,
        though the derivations remain equivalent.}{%
        Restricted to this subspace,
        $x \mapsto \lse(x)$ is strictly convex,
        and its gradient is the projection of the standard gradient onto this subspace.
        Correspondingly,
        we set the last column of the parameter $\Theta_i$ of \eqref{eq: switch_model} to zero throughout the paper,
        i.e.,$\Theta_{i, d} = \boldsymbol{0}$ for all $i\in\Xi$.
        }
\end{remark}

\subsection{Connection with different models}\label{sec: connection_models}
The switching mechanism \eqref{eq: switch_model} depends on both the active subsystem index $\xi_t$ (mode) and the state history $z_t$.
This general formulation encompasses three special cases,
each corresponding to established models in the literature where the switching mechanism depends on only a subset of these variables:
\begin{enumerate}
    \item Static switching:
    \begin{equation}\label{eq: static_switching}
        p(\xi_{t+1} \!\mid\! z_t, \xi_t\!\midsc \Theta) = p(\xi_{t+1} \!\midsc \Theta) = \sigma(\Theta)  
        \vspace{-0.5\baselineskip}
    \end{equation}
    with $\Theta_i = \Theta \in \re^{1 \times d}$. 
    Such model is commonly used in mixture models such as Gaussian mixture models \cite[\S 9.2]{bishop2006pattern},
    see \cref{ex:gmm}.
    \item Mode-dependent switching:
    \begin{equation}\label{eq: mode_dependent_switching}
        p(\xi_{t+1} \!\mid\! z_t, \xi_t\!\midsc \Theta) = p(\xi_{t+1} \!\mid \!\xi_t \!\midsc \Theta) = \sigma(\Theta_{\xi_t})  
        \vspace{-0.5\baselineskip}
    \end{equation}
    with $\Theta_i \in \re^{1 \times d}$ for all $i \in \Xi$.
    Such model is commonly used in Markov jump systems \cite{costa2005discrete}.
    \item State-dependent switching:
    \begin{equation}\label{eq: state_dependent_switching}
        p(\xi_{t+1} \!\mid\! z_t, \xi_t\!\midsc \Theta) = p(\xi_{t+1} \!\mid\! z_t\! \midsc \Theta) = \sigma(\Theta^\top\! z_t)
        \vspace{-0.5\baselineskip}
    \end{equation}
    with $\Theta_i = \Theta \in \re^{ n_z \times d}$.
    Such model is commonly used in mixture of experts models \cite{jordan1994hierarchical}.
\end{enumerate}

\rev[reviewer3point1]{%
    Model \eqref{eq: subsys_model} covers various distributions, 
    including all members in the canonical exponential family \cite{wainwright2008graphical,efron2022exponential}
    (e.g., categorical, normal, gamma, chi-squared distributions).%
}{%
    Model \eqref{eq: subsys_model} represents a generalization that encompasses the popular canonical exponential family \cite{wainwright2008graphical,efron2022exponential}, 
    which includes categorical, normal, gamma, and chi-squared distributions, and is defined as
    \[
        p(y \mid z, \xi=j\midsc \boldsymbol{\beta}) = \eta(y, z) \exp[\inner{\beta_j}{\sstat(y, z)} - \mathcal{A}(z, \beta_j)],
        \vspace{-0.5\baselineskip}
    \]
    where
    $\eta: \re^{n_y} \times \re^{n_z} \rightarrow \re_{++}$ is the base measure,
              $\sstat: \re^{n_y} \times \re^{n_z} \rightarrow \mathcal{B}$ the sufficient statistics,
              and
              $\logpar: \re^{n_z} \times \mathcal{X}_\beta \rightarrow \R$ the log-partition function.
    This generalization extends beyond the canonical exponential family to incorporate non-exponential family distributions as well.
}
See \cref{tab: example_distr} for more examples.
By modifying the choice of \eqref{eq: def_history}, 
the subsystem \eqref{eq: subsys_model} covers many commonly studied models.
Here we list a few examples.
\begin{enumerate}[wide]
\revnew{2}{%
\subsubsection{Static distributions}}{%
\item Static distributions:}
Although our framework is designed to deal with dynamical systems, 
it also covers the classical case where the measurements $y_t, y_{t'}$, for $t \neq t'$,
are mutually independent. 
In this case, we have 
\begin{equation} \label{eq:static-subsys}
    p(y_{t+1} \mid z_t, \xi_{t+1} = j; \dBeta) = p(y_{t+1} \mid \xi_{t+1}=j ; \dBeta).
\end{equation}
This case is trivially recovered when setting $\Upsilon$ in \eqref{eq: def_history} equal 
to a constant, i.e., $\Upsilon \equiv \bar{z}$ for some $\bar{z} \in \R^{n_z}$.

\begin{example}[Gaussian mixture model]\label{ex:gmm}
    The Gaussian mixture model is recovered as a special case of 
    \eqref{eq: discrete_latent}, where \eqref{eq: switch_model} is given by 
    \eqref{eq: static_switching}, and 
    \eqref{eq: subsys_model} satisfies \eqref{eq:static-subsys} with\vspace{-2mm}
    \[ 
        p(y_{t+1} \mid \xi_{t+1} = j; \dBeta) = \gauss(\mu_j, \Sigma_j).
    \]
    Taking $n_z = 1$, and $\Upsilon \equiv 1$, we recover the Gaussian case 
    in \cref{tab: example_distr}, with \revnew{2}{natural}{} parameters
    \rev{$A$}{$L$}$ = \mu_j \in \R^{n_y}$, and $\Sigma = \Sigma_j \succ 0$.
    After the change of variables
    $B_j = \Sigma^{-1}\mu_j,\;
    \Lambda_j = \Sigma_j^{-1}$,
    we can define  
    $\beta_j = (B_j, \Lambda_j)$, and \rev[newRef-GMM]{$\ell, f$, and $g$}{$f$ and $G$} are given 
    by \cref{tab: example_distr}.
    (Note that the subscript $j$ is dropped from the parameters in the table.)
    Since $\Lambda_j \succ 0$ by construction,
    the \revnew{2}{natural}{} parameters $(\mu_j, \Sigma_j)$ can be recovered easily by inverting
    the change of the variables.
\end{example}

\revnew{2}{%
\subsubsection{Dynamical systems}}{
\item Dynamical systems:}
Besides the static case, \eqref{eq: subsys_model} additionally 
covers many classical models for dynamical systems. 
We illustrate this with some particular choices of the noise 
distribution, but the derivations can be carried out analogously
for other choices.

\begin{example}[State-space model ($\ell_1$-Laplace distributioned noise)] \label{ex: state-space}
    Consider a dynamical system with state $x_t \in \R^{n_x}$, governed by the dynamics\vspace{-1mm}
    \[
        x_{t+1} = A_{\xi_{t+1}} x_t + w_t,\; w_t \sim \mathrm{\mathsf{Laplace}}(\mu, \Sigma),
        \vspace{-0.25\baselineskip}
    \]
    where $A_{i}$ for $i \in \Xi$ denote parameters of the dynamics,
    $\mu$ and $\Sigma$ are parameters of the noise distribution,
    whose probability density function is given by \cite[eq. 2]{aravkin2011}
    \begin{equation} \label{eq: laplace-dens}
    p_w(w_t) \dfn \det(2 \Sigma)^{-\frac{1}{2}} \exp( -\sqrt{2} \|\Sigma^{-\frac{1}{2}} (w_t - \mu)\|_{1} ).
    \end{equation}
    Since $p(x_{t+1} \mid x_t, \xi_{t+1}=j) = $\rev{$p_w(y_{t+1} - A_{j} y_t)$}{$p_w(x_{t+1} - A_{j} x_t)$},
    it follows directly 
    from \eqref{eq: laplace-dens} 
    that
    \[
        (x_{t+1} \mid x_t, \xi_{t+1}=j) \sim \mathrm{\mathsf{Laplace}}(\mu + A_j x_t, \Sigma).
    \]
    Defining 
    $y_t = x_t$ and 
    $z_t = \Upsilon(y_{t}, \dots, y_{t-t_y+1}) = (y_t, 1)$,
    and $L_j = [ \begin{matrix} A_j & \mu \end{matrix} ]$ we obtain that
    \[
        (y_{t+1} \mid z_t, \xi_{t+1}=j) \sim \mathrm{\mathsf{Laplace}}(L_j z_t, \Sigma),
    \]
    which coincides with \cref{tab: example_distr}.
    As in \cref{ex:gmm}, the original parameters $(\mu, A_j, \Sigma)$ can be retrieved from $\beta_j$, 
    by inverting the change of variables.
\end{example} 

\begin{example}[Autoregressive model (Student's t-distributed noise)]\label{ex: autoregressive_model}
    Consider the \revnew{2}{switching}{} system
    \begin{align*}
        y_{t+1} = &\; a_{\xi_{t+1}, 1} y_{t-t_y+1} + \dots + a_{\xi_{t+1}, t_y} y_t + w_{t+1} \\
        w_{t+1} \sim&\; \mathrm{\mathsf{St}}_{\nu}(\mu, \Sigma),
    \end{align*}
    where the additive noise is distributed by Student's t distribution.
    Let use define $A_{j} \dfn (a_{j, 1} , \dots, a_{j, t_y})$, for $j \in \Xi$,
    and $x_t \dfn (y_{t-t_y+1}, \dots, y_t)$, so that $y_{t+1} = A_{\xi_{t+1}} x_t + w_{t+1}$.
    Since Student's t distribution is closed under affine mappings \cite[eq. (4.1)]{roth2012multivariate},
    i.e., 
    \[ 
        (y_{t+1} \mid x_t, \xi_{t+1} = j) \sim \mathrm{\mathsf{St}}_{\nu}(\mu + A_j x_t, \Sigma),
    \]
    we can proceed analogously to \cref{ex: state-space},
    denoting the parameters $L_j = (A_j, \mu_j)$, 
    and setting $z_t = \Upsilon(y_t, \dots, y_{t+1-t_y}) = (x_t, 1)$, so 
    we recover
    $$
        (y_{t+1} \mid x_t, \xi_{t+1} = j) \sim \mathrm{\mathsf{St}}_{\nu}(L_j z_t, \Sigma).
        \vspace{-0.5\baselineskip}
    $$
    The corresponding functions \rev[newRef-autoregressor]{$f, \ell, g$}{$f$ and $G$} in \cref{eq: subsys_model},
    are now obtained directly from \cref{tab: example_distr}.
\end{example}
\end{enumerate}

We highlight the fact that the model can be generalized into non-autonomous systems 
when $z_t$ in \eqref{eq: def_history} depends on the input-output history.
For clarity of exposition, however, 
we focus on autonomous systems.
The extension to the non-autonomous system follows analogously.

\revnew{2}{}{%
\begin{remark}[On identifiability]\label{rem: identifiability}
    Identifiability results exist for restricted settings,
    such as static switching with exponential family subsystems \cite{jiang1999identifiability},
    or considering subsystem parameters alone without considering the switching mechanism parameters \cite{petreczkyMinimalityIdentifiabilityDiscretetime2020,vidalObservabilityIdentifiabilityJump2002a}.
    For the general case of jointly identifying switching mechanism and subsystem parameters, 
    a unified analysis remains open.
    Following established practice in the field \cite{bemporad2018fitting,leoni2025explainable,ferrari-trecateClusteringTechniqueIdentication,pillonettoNewKernelbasedApproach2016,pigaEstimationJumpBox2020,lauerGlobalOptimizationLowdimensional2018,bianchiModelStructureSelection2021a,balenzuela2022parameter},
    this work focuses on algorithmic development,
    and we view a rigorous identifiability analysis as an important direction for future work.
\end{remark}}

\subsection{Regularized maximum likelihood estimation (MLE)}

Given a trajectory $\ve{y} = \{y_t\}_{t = 0}^T$ generated by system \eqref{eq: discrete_latent} with initialization $z_0$,
and $\modeSeq = \{\mode{t}\}_{t=0}^{T}$ as the latent mode sequence,
we aim to estimate the 
parameter $\theta = (\dTheta, \dBeta) \in \mathbb{T} \dfn \R^{d \times n_z \times d} \times \mathcal{B}^d$ in model \eqref{eq: discrete_latent} by solving 
\vspace{-1mm}%
\begin{equation}\label{eq: reg_loss}
	\minimize_{\theta} \bar{\loss}(\theta) =  \loss(\theta) + \reg(\theta)
    \vspace{-3mm}
\end{equation}
with the negative log-likelihood (NLL) $\loss: \mathcal{X}_{\theta} \subseteq \mathbb{T} \to \R$
\begin{equation}\label{eq: loss}
	\begin{aligned}
		\loss(\theta)
		=\!{-}\!\ln
		p(\ve{y}, z_0 \!\midsc \theta)
		=\!{-}\!\ln\!
		\big(\!\textstyle{\tlsum_{\modeSeq\in \Xi^{T{+}1}}}\!
		p(\ve{y}, \ve{\xi}, z_0 \!\midsc \theta)
		\!\big),
	\end{aligned}
\end{equation}
and $\reg$ being a regularizer.
The joint density $p(\ve{y}, \ve{\xi}, z_0 \midsc \theta)$ can be factorized as 
\vspace{-\baselineskip}
\vspace{-\baselineskip}
{\small%
    \begin{equation*}%
        \begin{aligned}
            p(&\ve{y}, \modeSeq, z_0 \midsc \theta) 
        = p(y_0, \xi_0, z_0) p(\modeSeq_{1:T}, \ve{y}_{1:T} \mid y_0, \xi_0, z_0 \midsc \theta) \\
        &= p(y_0, \xi_0, z_0) \prod\limits_{t = 0}^{T-1} p(\xi_{t+1}, y_{t+1} \mid \modeSeq_{0:t}, \ve{y}_{0:t}, z_0 \midsc \theta) \\
        &= p(\xi_0, z_0) \prod\limits_{t = 0}^{T-1} 
            p(y_{t+1} \mid z_{t}, \xi_{t+1} \midsc \dBeta) \; %
          p(\xi_{t+1} \mid z_t, \xi_t\midsc \dTheta),
        \end{aligned}
    \end{equation*}
}%
where the last equation follows from the conditional independence defined in model \eqref{eq: discrete_latent} and the definition of $z_t$ in \eqref{eq: def_history}.
Since our goal is to identify the parameter $\theta$
for a given sequence $\ve{y}$ and initialization $z_0$, 
we introduce a shorthand notation for \rev{$\ell$ and $g$}{$G$}, 
with a subscript $t$ denoting the dependence on data $z_t$ and $y_{t+1}$,
i.e.,
\vspace{-0.5\baselineskip}
\rev{
    \begin{align*}%
        \ell_{t}(\beta)          =&\; \ell(y_{t+1}, z_t, \beta), \\
        g_{t}(\beta)             =&\; g(y_{t+1}, z_t, \beta),    
    \end{align*}%
}{%
    \begin{equation}\label{eq: shortcut}
        G_{t}(\beta)             = G(y_{t+1}, z_t, \beta),
    \end{equation}
}
which, combined with the model definition \eqref{eq: discrete_latent}, yields
\vspace{-0.5\baselineskip}
{%
    \begin{equation}\label{eq: expansion_nll}
        \begin{alignedat}{3}
            \loss(\theta)
            =&\;  %
            -\lse
            \Big(
                \ve{c} - \Psi(\dTheta) - \Phi(\dBeta) 
            \Big),
        \end{alignedat}
    \end{equation}
}%
where we defined
$\ve{c}, \Psi(\dTheta), \Phi(\dBeta)\!\in\!\re^{d^{T{+}1}}$ with elements%
\footnote{
    We implicitly define an index function $i : \Xi^{T+1} \to \N$, and,
    for ease of notation, simply write subscripts $\ve{\xi}$ instead of $i(\ve{\xi})$ to index a vector by $i(\ve{\xi})$.
    E.g., we write $\ve{c}_{\ve{\xi}} = \ve{c}_{i(\ve{\xi})}$.
}%
\begin{subequations}\label{eq:defs-L}
\begin{alignat}{2}
    \ve{c}_{\modeSeq} &\dfn 
    \ln p(\xi_0, z_0) + T\ln C,\\
\Psi_{\modeSeq}(\ve{\Theta}) &\dfn 
        \tsum_{t=0}^{T-1} \lse(\Theta_{\xi_{t}}^\top z_t) - \Theta_{\xi_t, \xi_{t+1}}^\top z_t,
\label{eq:Lt2-def}\\
    \ifshowboth
    \rev{\Phi_{\modeSeq}(\ve{\beta})}{}
        & \rev{\dfn
        \tsum_{t=0}^{T-1} f(\ell_{t}(\beta_{\xi_{t+1}})) + g_{{t}}(\beta_{\xi_{t+1}}),\nonumber}{} \\
	\else
		{}
    \fi
    \rev{}{\Phi_{\modeSeq}(\ve{\beta})} & \rev{}{
            \dfn\tsum_{t=0}^{T-1} f(G_{t}(\beta_{\xi_{t+1}})),\label{eq:Lt1-def}
        }
\end{alignat}
\end{subequations}
for all $\ve{\xi} \in \Xi^{T+1}$\revnew[notation]{2}{.}{,%
where $\Theta_{\xi_t, \xi_{t+1}}$ is the $\xi_{t+1}$-th column of matrix $\Theta_{\xi_t}$.}
It is clear from this formulation that 
$\loss$ has
domain $\mathcal{X}_{\theta} \dfn \R^{d \times n_z \times d} \times \mathcal{X}_\beta^d$,
which under \cref{asp: subsys_model} is an open set.
We study the \emph{regularized} MLE problem, 
since in general $\loss$ can be unbounded below. 
\vspace{-0.5\baselineskip}
\begin{example}[Unbounded NLL {\cite[\S 9.2.1]{bishop2006pattern}}] \label{ex: unbounded_nll}
Consider the identification of a 1D Gaussian mixture model with $\reg(\theta) = 0$.
In this case, the switching mechanism \eqref{eq: switch_model} is static, as described in \eqref{eq: static_switching},
and
the subsystem dynamics \eqref{eq: subsys_model} is a static model with Gaussian distribution, as described in \cref{ex:gmm}.
As both the switching mechanism and subsystem dynamics are static, 
the NLL can be expressed without temporal dependencies, in the form
$\loss(\theta) = \sum_{t=0}^T-\ln (\sum_{\xi_t \in \Xi} p(y_t\mid \xi_t \midsc \dBeta)\, p(\xi_t \midsc \dTheta))$
where $p(y_t\mid \xi_t \midsc \dBeta)$ is the pdf for Gaussian distribution and $p(\xi_t \midsc \dTheta) \neq 0$ for all $\xi_t \in \Xi$.
Suppose that the mean of subsystem $i$ is equal to a specific data point $y_t$, i.e., $\mu_i = y_t$. 
Then, the corresponding term in the likelihood 
$p(y_{t} \mid \xi_{t}=i\midsc \dBeta) =(2 \pi)^{-\frac{1}{2}} \lambda_i^{\frac{1}{2}}$,
where $\lambda_i \dfn \sigma_i^{-2}$ represents the precision, as explained in \cref{tab: example_distr}.
When the precision $\lambda_i$ tends to $+\infty$,
this term increases towards $+\infty$,
driving the NLL $\loss$ to $-\infty$.
This behavior is illustrated in \cref{fig: gmm_nll}.
Importantly, with $d > 1$, this degenerate case cannot be mitigated by introducing more data,
since one can always achieve negative infinite cost by 
selecting the parameters of one subsystem to tightly fit around a single 
datapoint (i.e., with infinite precision/zero variance), while maintaining 
a finite likelihood for the remaining 
data using the remaining submodels.
\begin{figure}[!tbph]
    \centering
    \includegraphics[width=0.4\textwidth]{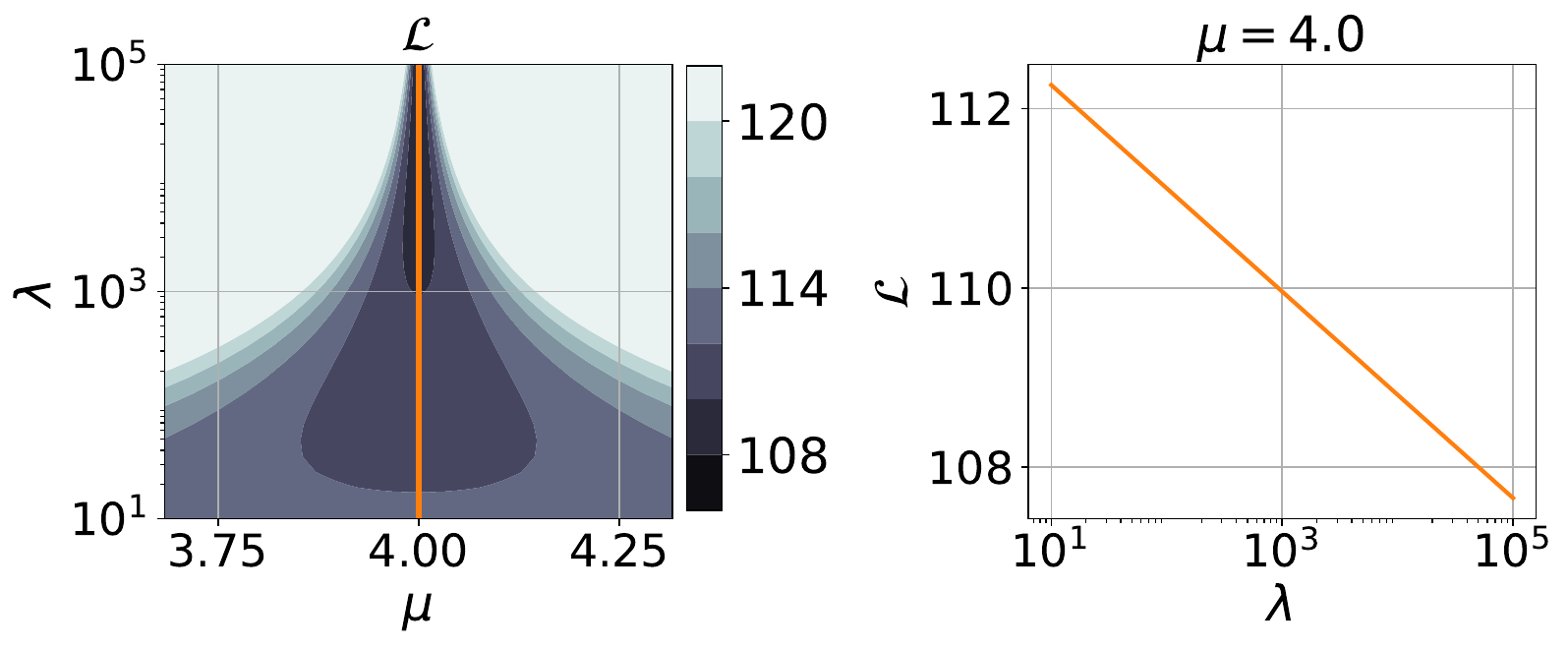}
    \caption{The NLL $\loss$ of a Gaussian mixture model 
    $p(y) = 0.5 \gauss(y\midsc 0, 1) + 0.5 \gauss(y \midsc \mu, \lambda^{-1})$
    evaluated on $T = 50$ data points with $y_T = 4$.}
    \label{fig: gmm_nll}
\end{figure}%
\end{example}%
\begin{assumption}[Regularizers]\label{asp: regularizers}
    The regularizer $\reg: \mathcal{X}_{\theta} \to \R$ is separable in
    $\dTheta, \dBeta$, i.e., 
    \rev{}{there exist functions $\reg_1$, $\reg_2$ such that}
    $\reg(\theta) = \reg_1(\dTheta) + \reg_2(\dBeta)$,
    and satisfies the following conditions
    \begin{inlinelist}
        \item $\reg_1, \reg_2$ are convex;
        \item there exist $\kappa_1, \kappa_2 > 0$ and $c_1, c_2 \in \R$ such that for any $\modeSeq \in \Xi^{T+1}$:
    \end{inlinelist} 
        \vspace{-2\baselineskip}
        {%
            \small %
            \begin{subequations}
                \begin{alignat}{4}
                    \Psi_{\modeSeq}(\dTheta) &+ \reg_1(\dTheta) &\geq&\; \kappa_1 \Vert \dTheta \Vert  &\,+\, c_1 \label{eq: reg_Theta_condition}, \\
                    \Phi_{\modeSeq}(\dBeta)  &+ \reg_2(\dBeta)  &\geq&\; \kappa_2 \Vert \dBeta \Vert   &\,+\, c_2. \label{eq: reg_beta_condition}
               \end{alignat}
            \end{subequations}
        }%
\end{assumption}%
\vspace{-\baselineskip}
In the upcoming analysis, 
we demonstrate that \cref{asp: regularizers} prevents the scenario described in \cref{ex: unbounded_nll}.
\rev[reviewer3point3]{
    A common regularizer satisfying \eqref{eq: reg_Theta_condition} is a quadratic function.
}{
    Common regularizers satisfying \eqref{eq: reg_Theta_condition} include 
    $\reg_1(\dTheta) = \tfrac{1}{2}\sum_{i = 1}^d \norm{\Theta_i}^2_{\mathrm{F}}$ (Ridge regression),
    and $\reg_1(\dTheta) = \sum_{i = 1}^d \norm{\Theta_i}_{\mathrm{*}}$ (LASSO regression) where $\norm{\cdot}_*$ is the nuclear norm. 
}
However, the selection of $\reg_2$ to satisfy \eqref{eq: reg_beta_condition}
depends on the specific form of \rev{$f, \ell_t$, and $g_t$}{$f$ and $G$}.
A potential candidate is the log-probability derived from the conjugate prior \cite[\S 2.4.2]{bishop2006pattern}, if it exists.
To illustrate this,
we present the following example.
\begin{example}[Regularization for Gaussian distributions]\label{ex: regularization_1d_gaussian}
    Consider \cref{eq: subsys_model} as a 1D Gaussian distribution. 
    Recall from \cref{tab: example_distr} that for $t \geq 0$
    \rev{
    \[
        \textstyle f(\ell_{t}(\beta_{i})) + g_{{t}}(\beta_{i}) = \frac{1}{2\lambda_{i}}(\lambda_{i} y_{t+1} - b_{i} z_t)^2 - \frac{1}{2}\ln \lambda_{i}
    \]}{
    \[
        \textstyle f(G_{t}(\beta_{i})) = \frac{1}{2\lambda_{i}}(\lambda_{i} y_{t+1} - b_{i} z_t)^2 - \frac{1}{2}\ln \lambda_{i}
    \]
    }
    for all $i \in \Xi$.
    Thus, by \eqref{eq:Lt1-def},
    $\Phi_{\modeSeq}(\dBeta) \geq \sum_{t = 0}^{T-1} -\frac{1}{2} \ln \lambda_{\xi_{t+1}}$.
    Let $\reg_2(\dBeta) = \frac{1}{2}\sum_{j = 1}^d  \lambda_j -\ln\lambda_j + b_j^2 / \lambda_j$,
    which is derived from the conjugate prior of a Gaussian distribution \cite[\S 2.3.6]{bishop2006pattern}.
    We first show that with this regularization, there exists some $\kappa > 0$, such that for all sequence $\modeSeq \in \Xi^{T+1}$,
    \begin{equation}\label{eq: example_proven}
        \textstyle
        \sum_{t= 0}^{T-1}\!-\frac{1}{2}\ln \! \lambda_{\xi_{t+1}}\! + \reg_2(\dBeta) \geq c + \kappa \sum_{j = 1}^d \abs{b_j} + \abs{\lambda_j}.
    \end{equation}
    This is equivalent showing for any $\alpha > 0$, there exists $\kappa > 0$ such that for all $a > 0$,
    \[ 
        u(\lambda, b) = -a \ln \lambda + \alpha \lambda + \alpha \tfrac{b^2}{\lambda} - \kappa \abs{b} - \kappa \lambda
    \]
    is lower bounded uniformly.
    By lower bounding the quadratic term of $b$, we obtain
    \(
        u(\lambda, b) \geq -a \ln \lambda + \big(\alpha - \frac{1}{4 \alpha}\kappa^2 - \kappa\big) \lambda.
    \)
    This lower bound has a unique minimum for all $\lambda > 0$ if $\alpha - \frac{1}{4 \alpha}\kappa^2 - \kappa > 0$.
    Given $a, \alpha > 0$,
    one can always choose $\kappa \in (0, 2(\sqrt{2} - 1)\alpha)$ to satisfy this condition.
    Thus,
    there exists $\kappa > 0$ such that 
    $u$ is lower-bounded. 
    Because 
    $\sum_{j = 1}^d \abs{b_j} + \abs{\lambda_j} \geq \norm{\dBeta} = \sqrt{\sum_{j = 1}^d \abs{b_j}^2 + \abs{\lambda_j}^2}$,
    \eqref{eq: example_proven} implies that \eqref{eq: reg_beta_condition} is satisfied.
\end{example}

It is evident from \eqref{eq: expansion_nll} and \eqref{eq:defs-L} that
the NLL $\loss$ is nonconvex, making
problem \eqref{eq: reg_loss} difficult to solve using off-the-shelf methods.
Moreover,
$\loss$ entails the marginalization over a discrete sequence $\modeSeq$,
of which the size grows exponentially with the horizon $T$.
To prevent an exponential growth in computational complexity, 
evaluating $\loss(\theta)$ typically involves a dynamic programming procedure,
with a complexity proportional to the sequence length.
However,
due to the sequential nature of dynamic programming,
each gradient evaluation $\nabla\loss(\theta)$ must be propagated through the entire sequence $\modeSeq$ in order,
leading to highly costly iterations.
The goal of this work is therefore to develop a solution method that reduces 
the number of inherently expensive iterations by maximally exploiting 
the known structure of the problem class.
By adopting the majorization-minimization (MM) principle \cite{hunter2004tutorial,lange2016mm}, with a
majorization model that provides a tighter upper bound of the cost landscape,
one can expect to make more progress per iteration than with general purpose (e.g., quadratic)
models underlying classical gradient-based methods.

\section{The identification method}\label{sec: mm}
This section introduces the proposed method \EMpp as a particular MM scheme for minimizing the regularized NLL.
An MM scheme iteratively performs two steps.
The \textit{majorization step} constructs
a surrogate function $\regSurrogate[k]$ at the current iterate $\thk \in \mathcal{X}_\theta$, satisfying
\begin{subequations}\label{eq: property_Q}
	\vspace{-0.5\baselineskip}
	\begin{align}
		\regSurrogate[k](\theta) \geq & \; \regLoss(\theta), \quad \forall \theta \in \mathcal{X}_\theta,  \label{eq: Q_pos} \\
		\regSurrogate[k](\thk) =      & \; \regLoss(\thk) \label{eq: Q_equal_func}.
		\vspace{-0.5\baselineskip}
	\end{align}
\end{subequations}
The \textit{minimization step},
consists of minimizing the surrogate function $\regSurrogate[k]$,
yielding the next iterate $\thkp$.
In particular, we construct a \emph{convex} surrogate function $\regSurrogate[k]$ satisfying \eqref{eq: property_Q}, such that a sequence of convex problems is solved instead of the original nonconvex problem \eqref{eq: reg_loss}.
The following lemma provides a first step towards the construction of such a surrogate function.
\begin{lemma}\label{lem: jensen_inequality}
    Let $\ve{y}$ be a trajectory
    of system \eqref{eq: discrete_latent} and $\thk \in \mathcal{X}_\theta$ an iterate.
    Then we can bound the NLL by
    \begin{equation*}
        \loss (\theta) \leq - \Pi(\thk)^\top \left( \ve{c} - \Psi(\dTheta) - \Phi(\dBeta)\right) + c_{\thk}   
    \end{equation*}
    \rev{with}{for all $\theta \in \mathcal{X}_\theta$, where} $c_{\thk} = \textstyle{\sum_{\modeSeq \in \Xi^{T+1}}} p(\modeSeq \mid \stateSeq, z_0 \midsc \thk) \ln p(\modeSeq \mid \stateSeq, z_0 \midsc \thk)$ 
	and $\Pi(\thk) \in \re^{d^{T+1}}$ with element 
	\(
		\Pi_{\modeSeq}(\thk) = p(\modeSeq \mid \stateSeq, z_0 \midsc \thk)
	\)
	for all $\modeSeq \in \Xi^{T+1}$.
    Moreover, this relation holds with equality whenever $\theta = \thk$.
\end{lemma}
\vspace{-0.5\baselineskip}
\begin{pf}
    Using \eqref{eq: loss} and Jensen's inequality we have
	\vspace{-\baselineskip}
	\vspace{-\baselineskip}
	{%
		\small%
		\begin{align*}
			\loss(\theta) &= -\ln \textstyle \sum_{\modeSeq \in \Xi^{T+1}} \tfrac{p(\modeSeq \mid \stateSeq, z_0 \midsc \thk)}{p(\modeSeq \mid \stateSeq, z_0 \midsc \thk)} p(\ve{y}, \modeSeq, z_0\midsc \theta)\\
			&\leq - \textstyle \sum_{\modeSeq \in \Xi^{T+1}} p(\modeSeq \mid \stateSeq, z_0 \midsc \thk) \ln \tfrac{p(\ve{y}, \modeSeq, z_0\midsc \theta)}{p(\modeSeq \mid \stateSeq, z_0 \midsc \thk)} \\
			&= - \textstyle \sum_{\modeSeq \in \Xi^{T+1}} p(\modeSeq \mid \stateSeq, z_0 \midsc \thk) \ln p(\ve{y}, \modeSeq, z_0\midsc \theta) + c_{\thk},
		\end{align*}
	}%
    and this holds with equality for $\theta = \theta^k$.
    The claims follows from the observation that
    \begin{equation*}
        \ln p(\ve{y}, \modeSeq, z_0\midsc \theta) = \ve{c}_{\modeSeq} - \Psi_{\modeSeq}(\dTheta) - \Phi_{\modeSeq}(\dBeta). \qed
		\vspace{-0.5\baselineskip}
	\end{equation*}
\end{pf}%
\vspace{-\baselineskip}
By exploiting concavity of $f$,
the next proposition constructs a surrogate $\regSurrogate[k]$ satisfying \eqref{eq: property_Q}.
\begin{proposition}\label{prop:separable}
	Let $\ve{y}$ be a trajectory
	of system \eqref{eq: discrete_latent} initialized at $z_0$.
	Then $\regSurrogate[k](\theta) := \surrogate(\theta) + \reg(\theta)$ with
	\begin{subequations}\label{eq: surrogate}
		\begin{equation}
			\surrogate(\theta) := \Qi[0][k]{} + \Qi[1][k]{\dTheta} + \Qi[2][k]{\dBeta}\rev{,}{}
			\vspace{-0.25\baselineskip}
		\end{equation}
		constitutes a convex surrogate of $\regLoss$ satisfying \eqref{eq: property_Q}, where%
		\vspace{-\baselineskip}
		\vspace{-\baselineskip}
		\vspace{-\baselineskip}
		{%
			\small%
			\begin{align}
				\label{eq: initial-surrogate}
				\Qi[0][k]{} \dfn & c_{\thk} - \Pi(\thk)^\top \ve{c}, \\
				\Qi[1][k]{\dTheta}
				\dfn             & \;
				\Pi(\thk)^\top \Psi(\dTheta), \label{eq:switch-surrogate}\\
				\Qi[2][k]{\dBeta}
				\dfn             & \;
				\Pi(\thk)^\top \hat{\Phi}^k(\dBeta),\label{eq:subsys-surrogate}
			\end{align}
		}%
	\end{subequations}
	and where for all $\modeSeq \in \Xi^{T+1}$ we defined
	\rev{
	\begin{align*}
		\hat{\Phi}^k_{\modeSeq}(\dBeta) &\dfn
            \tlsum_{t = 0}^{T-1} f[\ell_t(\beta^k_{\xi_{t+1}})](\ell_t(\beta_{\xi_{t+1}})) + g_{{t}}(\beta_{\xi_{t+1}}) \\
		f[\bar x](x) &\dfn f(\bar x) + f'(\bar x) (x - \bar x). \nonumber
	\end{align*}}{
		\vspace{-\baselineskip}
		\begin{align}
			\hat{\Phi}^k_{\modeSeq}(\dBeta) &\dfn
				\tlsum_{t = 0}^{T-1} f[G_t(\beta^k_{\xi_{t+1}})]\big(G_t(\beta_{\xi_{t+1}})\big)\label{eq: def_elem_surrogate_subsys}\\
			f[\bar x](x) &\dfn f(\bar x) + \inner{\nabla f(\bar x) }{x - \bar x}. \nonumber
		\end{align}
		\vspace{-\baselineskip}
	}
\end{proposition}
\vspace{-\baselineskip}
\begin{pf}
	By concavity of the function $f$, its linearization around a point $\bar x \in \dom f$
constitutes an upper bound, i.e., for all $x \in \dom f$: 
\rev{$f(x) \leq f(\bar x) + f'(\bar x) (x - \bar x)$,}{$f(x) \leq f(\bar x) + \inner{\nabla f(\bar x) }{x - \bar x}$,} 
and this holds with equality when $x = \bar x$.
Therefore, we have that
\(
    \Phi_{\modeSeq}(\dBeta) \leq \hat{\Phi}^k_{\modeSeq}(\dBeta)
\),
with equality holding for $\dBeta = \dBeta^k$.
In combination with the bound from \cref{lem: jensen_inequality}, this proves the claim.

\qed\end{pf}%
\vspace{-\baselineskip}
As the mappings $\Psi$, \rev{$g_{t}$, $\ell_{t}$}{$G_t$}, 
and $\reg$ are convex, so is
the surrogate problem
\(
\argmin_{\theta} \regSurrogate[k](\theta)
\).
We remark that the direct computation of $\surrogate$ based on 
\cref{eq: surrogate} is expensive,
as it involves evaluating the vector $\Pi(\theta^k)$ of dimension $d^{T+1}$.
\Cref{sec: e_step} describes a different representation of $\surrogate$ (cf.\,\cref{prop:separable-explicit}) that allows
for efficient evaluation, once certain quantities are computed (cf.\,\cref{alg:construct-Q}).
Exploiting separability of $\surrogate$ and $\reg$ w.r.t. $\dTheta, \dBeta$,
the resulting scheme is summarized in \cref{alg:em-conceptual}.
\begin{algorithm}
	\caption{\EMpp}
	\label{alg:em-conceptual}
	\begin{algorithmic}[1]
        \Require Data $\ve{y}$; Initialize $\theta^{0} \in \mathcal{X}_{\theta}$
		\For{$k = 0,1,\dots$ (until convergence)}
		\State Run \cref{alg:construct-Q} to construct $\Qi[1][k]{}, \Qi[2][k]{}$ 
		\State
		Solve (in parallel) %
		\vspace{-0.5\baselineskip}
		\begin{subequations}
			\vspace{-0.5\baselineskip}
			\begin{align}
				\dTheta^{k+1} \gets & \; \textstyle \argmin_{\dTheta} \Qi[1][k]{\dTheta} + \reg_1(\dTheta) \label{eq: prob_switch_iden} \\
				\dBeta^{k+1} \gets  & \; \textstyle \argmin_{\dBeta} \Qi[2][k]{\dBeta}  \,\, + \reg_2(\dBeta) \label{eq: prob_subsys_iden}
			\end{align}%
		\end{subequations}%
		\vspace{-2\baselineskip}
		\EndFor
	\end{algorithmic}
\end{algorithm}

We highlight that the loss evaluated at $\{\thk\}_{k\in\N}$ from \EMpp is nonincreasing, as formally stated below.
\begin{corollary}[Nonincreasing sequence]\label{cor: nonincreasing_sequence}
	The iterates $\{\thk\}_{k \in \N}$ generated by \EMpp (cf. \cref{alg:em-conceptual}) satisfy\vspace{-1mm}
	\begin{equation}\label{eq: nonincreasing}
		\forall k \in \N: \regLoss(\theta^k) \geq \regLoss(\theta^{k+1}).
		\vspace{-\baselineskip}
	\end{equation}
\end{corollary}%
\begin{pf}
	By \cref{prop:separable}, \eqref{eq: property_Q} holds, and 
    we have
    \begin{equation*}
        \forall k \in \N: \regLoss(\thkp) \leq \regSurrogate[k](\thkp) \leq \regSurrogate[k](\thk) = \regLoss(\thk).\qed
		\vspace{-0.5\baselineskip}
    \end{equation*}
\end{pf}

\section{Convergence analysis}\label{sec: conv_analysis}

\subsection{Subsequential convergence}
We start by presenting two lemmata that will prove useful in the upcoming convergence analysis.
The first states that the loss function $\regLoss$ is lower bounded, and that solutions to the surrogate problems 
$\argmin_{\theta} \regSurrogate[k](\theta)$ remain in a compact set.
The second relates to the directional differentiability of $\regLoss$ and $\regSurrogate[]{}$.

\begin{lemma}\label{lem: general_compactness}
    Under \cref{asp: subsys_model,asp: regularizers},
    \begin{lemitem}[wide]
        \item $\regLoss$ is lower bounded;
        \item there exists a compact set $\Omega \subseteq \mathcal{X}_\theta$ 
    that contains the iterates $\left\{ \theta^{k} \right\}_{k \in \N}$ generated by \EMpp (cf. \cref{alg:em-conceptual}).
    \end{lemitem}
\end{lemma}
\vspace{-\baselineskip}

\begin{pf}
    \rev[general-compactness]{%
    By \cref{asp: regularizers}, \cref{eq:Lt1-def,eq:Lt2-def} we have
\begin{alignat*}{6}
    &\min_{\modeSeq \in \Xi^{T+1}}  \Psi_{\modeSeq} (\dTheta) \,&+&\, \reg_1(\dTheta) 
    &\,\geq&\, \kappa_1 \norm{\dTheta}  &+\, c_1, \\
    &\min_{\modeSeq \in \Xi^{T+1}} \Phi_{\modeSeq}(\dBeta)    \,&+&\, \reg_2(\dBeta)
    &\,\geq&\,  \kappa_2 \norm{\dBeta}  &+\, c_2.
\end{alignat*}%
Following from $\lse(x) \leq \max_i x_i + \ln n_x$,
\vspace{-\baselineskip}
\vspace{-\baselineskip}
{%
\small%
\begin{align*}
    \regLoss(\theta) = &\;-\lse\left( \ve{c} - \Psi(\dTheta) - \Phi(\dBeta) - \reg(\theta) \right)\\
                     \geq&\;  \min_{\modeSeq \in \Xi^{T+1}} -\ve{c}_{\modeSeq} + \Psi_{\modeSeq}(\dTheta) + \Phi_{\modeSeq}(\dBeta) + \reg(\theta)  - \ln d^{T+1},
\end{align*}
}%
where $\reg(\theta) = \reg_1(\dTheta) + \reg_2(\dBeta)$.
Thus,
there exist $\kappa_3 > 0$, $c_3 \in \R$ for which
\begin{align*}
    \regLoss(\theta) 
    &\geq \kappa_1 \Vert \dTheta \Vert + \kappa_2 \Vert \dBeta \Vert + c_3 \geq \kappa_3 \Vert \theta \Vert + c_3,
\end{align*}
and the last step follows by equivalence of norms in Euclidean spaces,
demonstrating that $\regLoss$ is bounded below and that
\(
    \liminf_{\Vert \theta \Vert \to \infty} \frac{\regLoss(\theta)}{\Vert \theta \Vert} = \kappa_3 > 0.
\)
Thus, $\regLoss$ is \textit{level-coercive} \cite[Definition 3.25]{rockafellar2009variational},
implying it has bounded level-sets \cite[Corollary 3.27]{rockafellar2009variational}.
As shown in \cref{cor: nonincreasing_sequence}, the sequence $\left\{ \regLoss(\theta^k) \right\}_{k \in \N}$ is nonincreasing.
Hence, the iterates $\left\{ \theta^k \right\}_{k \in \N}$ 
remain in the bounded level-set $\Omega \dfn \left\{ \theta \mid \regLoss(\theta) \leq \regLoss(\theta^0) \right\} \in \mathcal{X}_{\theta}$.
By continuity of $\regLoss$, $\Omega$ is closed \cite[Theorem 1.6(c)]{rockafellar2009variational}, and thus compact.

    }{%
        See \cref{prf: general_compactness}.
    }
\qed\end{pf}

\begin{lemma}\label{lem: directional_derivative}
    If \cref{asp: subsys_model,asp: regularizers} hold, then $\regLoss$ and $\regSurrogate[]{}$ are directionally differentiable at $\theta^k$ along any $v = (v_{\dTheta}, v_{\dBeta}) \in \mathbb{T}$.
    Moreover,
    \vspace{-\baselineskip}
    \begin{equation*}
        \regLoss ' (\thk \midsc v) = 
        \regSurrogate[k'](\thk \midsc v).
    \end{equation*}%
\end{lemma}%
\vspace{-2.3\baselineskip}

\begin{pf}%
    \rev[directional-derivative]{
    Denote $L_{\modeSeq}(\theta) = \ve{c}_{\modeSeq} - \Psi_{\modeSeq}(\dTheta) - \Phi_{\modeSeq}(\dBeta)$ for all $\modeSeq \in \Xi^{T+1}$, such that $\loss(\theta) = -\lse(L(\theta))$.
Under \cref{asp: subsys_model}, \rev{the functions $\ell_t, g_t$, and $-f$}{the mapping $G_t$ and the function $-f$} are strictly continuous and convex ($f$ is concave).
Since convexity ensures directional differentiability\cite[Th.\,23.1]{rockafellar1997convex}, it follows that \rev{$f, \ell_t, g_t$}{$f$ and $G_t$} are B-differentiable on their domains \cite[Definition 3.1.2]{facchinei2003finite}.
Then, by the composition rule \cite[proposition 3.1.6]{facchinei2003finite} also $\loss$ is B-differentiable on its domain, and for all $v \in \mathbb{T}$ the directional derivative at $\theta$ along $v$ equals
\vspace{-\baselineskip}
\vspace{-\baselineskip}
{%
    \small%
        \begin{align}
            \loss'(\theta \midsc v) 
                &= -\lse'\Big(
                    L(\theta) \midsc 
                    L'(\theta\midsc v)
                \Big)
                = -\inner{
                    \nabla \lse\left(
                        L(\theta)
                    \right)
                }{
                    L'(\theta\midsc v)
                } \nonumber \\
                &= \Big\langle
                    \nabla \lse\left(
                        L(\theta)
                    \right), 
                \Psi'(\dTheta \midsc v_{\dTheta})
                + \Phi'(\dBeta \midsc v_{\dBeta}) 
                \Big\rangle. \label{eq: directional_dev_L}
        \end{align}
}%
Since by \eqref{eq: loss} and \eqref{eq: expansion_nll} we have
\(
    L_{\modeSeq}(\theta) = \ln \pjoint ,
\)
it follows that
\(
    \exp(L_{\modeSeq}(\theta)) = \pjoint.
\)
As the gradient of $\lse$ is the softmax function, i.e.,
\(
    \nabla \lse\left( L(\theta)\right)  = \sigma(L(\theta)),
\)
it follows that
\vspace{-\baselineskip}
\vspace{-\baselineskip}
{%
    \small%
    \begin{align*}
        \sigma_{\modeSeq}(L(\theta))
            =&\; \frac{
                \pjoint 
            }{
                \sum_{\modeSeq' \in \Xi^{T+1}} p(\stateSeq, \modeSeq', z_0 \midsc \theta)
            } \\
            =&\; \frac{
                \pjoint 
            }{
                p(\ve{y}, z_0 \midsc \theta)
            } 
            = p(\modeSeq \mid \ve{y}, z_0\midsc \theta)  = \Pi_{\modeSeq}(\theta)
    \end{align*}
}%
for all $\modeSeq \in \Xi^{T+1}$.
We thus obtain that
\begin{equation*}
    \loss'(\theta^k \midsc v) = \Pi(\theta^k)^\top 
    \left(\Psi'(\dTheta^k \midsc v_{\dTheta}) + \Phi'(\dBeta^k \midsc v_{\dBeta})\right).
\end{equation*}
On the other hand, we have by definition of $\surrogate$ that
\begin{equation*}
    \surrogate[k'](\theta^k \midsc v) 
    = \Pi(\theta^k)^\top 
    \left(\Psi'(\dTheta^k \midsc v_{\dTheta}) + \hat{\Phi}^{k'}(\dBeta^k \midsc v_{\dBeta})\right).
\end{equation*}
At any point $\bar x \in \dom f$, 
\rev{the derivative of $f$ equals the derivative of its linearization around $\bar x$, i.e.,
\(
    f'(\bar x) = [f(\bar x) + f'(\bar x) (\cdot - \bar x)]'(\bar x).
\)}{
    the gradient of $f$ equals the gradient of its linearization around $\bar x$, i.e.,
    \(
    \nabla f(\bar x) = \nabla [f(\bar x) + \inner{\nabla f(\bar x)}{\cdot - \bar x}](\bar x).
    \)
}
It follows that
\(
    \Phi_{\modeSeq}'(\dBeta^k \midsc v_{\dBeta}) = \hat{\Phi}^{k'}_{\modeSeq}(\dBeta^k \midsc v_{\dBeta}) 
\)
for all $\modeSeq \in \Xi^{T+1}$, and hence that
\(
    \loss' (\thk \midsc v) = 
    \surrogate[k'](\thk \midsc v)
    \) for all $v \in \mathbb{T}$.
Since $\reg$ is convex, it is also directionally differentiable \cite[Theorem 23.1]{rockafellar1997convex}.
Thus, $\regLoss = \loss + \reg$ and $\regSurrogate[k] = \surrogate + \reg$ are directionally differentiable on their domain and for all $v \in \mathbb{T}$ we obtain that 
\(
    \regLoss ' (\thk \midsc v) = 
    \regSurrogate[k'](\thk \midsc v).
\)

    }{See \cref{prf: directional_derivative}}
\qed\end{pf}
\vspace{-\baselineskip}
\begin{theorem}[Subsequential convergence]\label{thm: subsequential_convergence}
    Under \cref{asp: subsys_model,asp: regularizers}, 
    every limit point of the iterates $\left\{ \theta^{k} \right\}_{k \in \N}$ generated by \EMpp (cf. \cref{alg:em-conceptual}) is a stationary point of $\regLoss$.
\end{theorem}
\vspace{-\baselineskip}
\begin{pf}
    The proof follows \cite[Theorem 1]{razaviyayn2013unified}.
    By \cref{asp: regularizers} and \cref{lem: general_compactness},
    all iterates $\thk$ stay in a compact set $\Omega$.
    Thus, there exists a subsequence $\{\theta^{k_j}\}_{{k_j} \in \mathcal{K}\subseteq \N}$ that converges to some $\theta^\infty \in \Omega$.
    \rev{It satisfies}{We have}
    \begin{equation}\label{eq: subsequent_nonincreasing}
        \begin{aligned}
            &\; \regSurrogate[k_{j+1}](\theta^{k_{j+1}}) 
                = \regLoss(\theta^{k_{j+1}}) 
                \leq \regLoss(\theta^{k_j + 1}) \\ 
            &\; \quad\quad \leq  \regSurrogate[k_j](\theta^{k_j + 1}) \leq \regSurrogate[k_j](\theta), 
            \quad \forall \theta \in \mathbb{T},
        \end{aligned}        
    \end{equation}
    where we consecutively used \eqref{eq: Q_equal_func}, \eqref{eq: nonincreasing}, and \eqref{eq: Q_pos}.
    Taking the limit as $k \to \infty$ yields
    \begin{equation}\label{eq: limit_subsequent}
        - \infty < \regLoss(\theta^\infty) = \regSurrogate[\infty](\theta^{\infty}) \leq \regSurrogate[\infty](\theta), \quad\forall \theta \in \mathbb{T}. 
    \end{equation}
    The last inequality implies that $\theta^\infty$ is a (global) minimizer of $\regSurrogate[\infty]$.
    By the first-order necessary conditions of optimality
    we therefore have that
    \(
        (\regSurrogate[\infty])'(\theta^\infty \midsc v) \geq 0,
    \)
    for all $v \in \mathbb{T}$. 
    Thus, $\theta^\infty$ is a stationary point of $\regLoss$, since by \cref{lem: directional_derivative} we have
    \(
        \regLoss'(\theta^\infty \midsc v) \geq 0,
    \)
    for all $v \in \mathbb{T}$.
\qed\end{pf}

\subsection{Global sequential convergence}

Interpreting \EMpp (cf. \cref{alg:em-conceptual}) as a \emph{mirror descent} method, this section establishes global sequential convergence under a slightly more restrictive assumption.
\begin{assumption}\label{asp: smoothness}
    We assume that
    \begin{lemitem}
        \item \label{asp: local_Lipschitz_smoothness_g}\rev{$\ell_{t}$ and $g_{t}$ are}{$G_t$ is} Lipschitz smooth on any compact set $\Omega_\theta \subseteq \mathcal{X}_\theta$ for all $t \in \N_{[0, T-1]}$;
        \item \label{asp: local_Lipschitz_smoothness_f}$f$ is Lipschitz smooth on any compact set $\Omega_f \subseteq \mathcal{X}_f$;
        \item \label{asp: local_strong_convexity_of_kernels}\revnew{2}{}{In addition, one of the following conditions holds:}
        \begin{enumerate}[label=(\roman*), wide]
            \item \label{asp: reg_c2}the regularizer $\reg$ is of class $\mathcal{C}^2$ with $\nabla^2 \reg(\theta) \succ 0$ for all $\theta \in \mathcal{X}_\theta$.
            \item \label{asp: sufficient_data}\revnew{2}{}{%
            $\reg$ is Lipschitz smooth on any compact set $\Omega_\theta \subseteq \mathcal{X}_\theta$.
            There exists a nonempty set $\mathcal{I} \subseteq \N_{[1, n_g]}$,
            such that for all $i \in \mathcal{I}$,
            $g_i$ is of class $\mathcal{C}^2$ with
            $\nabla^2 g_i(u) \succ 0$ for all $u\in \mathcal{U}$,
            and
            \[
                \textstyle \sum_{t=0}^{T-1} z_t z_t^{\top} \succ 0, \quad \sum_{t=0}^{T-1}\sum_{i \in \mathcal{I}}\jac{\!}^\top_{i, t}\jac{\!}_{i, t} \succ 0,
            \]
            where 
            $\jac{\!}_{i, t}$ is the Jacobian of $T_i(y_{t+1}, z_t, \beta)$ w.r.t. $\beta$.}
        \end{enumerate}
    \end{lemitem}
\end{assumption}
\revnew{2}{}{%
Note that from \cref{asp: subsys_model},
$T_i$ is an affine mapping w.r.t. $\beta$ for all $i \in \N_{[1, n_g]}$.
Hence, the Jacobian $\jac{\!}_{i, t}$ only depends on data $(y_{t+1}, z_t)$.
The structure of \cref{asp: smoothness}~\labelcref{asp: local_strong_convexity_of_kernels}~\labelcref{asp: sufficient_data} resembles a persistent excitation condition.
This property is crucial for establishing the local strong convexity in \cref{lem: local_strong_convexity_of_kernels} below.
}%
We briefly introduce Bregman distances, which play a central role in the context of mirror descent methods.
\begin{definition}[Bregman distance {\cite{bregman1967relaxation}}] \label{def:bregman-distance}
    For a convex function 
    $h: \re^n \rightarrow \re\cup \{\infty\}$,
    which is continuously differentiable on $\intdom h \neq \emptyset$,
    the Bregman distance $D_h: 
        \re^n \times \re^n \rightarrow 
        \ere
        $ 
    is given by
    \vspace{-\baselineskip}
    \vspace{-\baselineskip}
    {\small%
    \begin{equation*}
        D_{h}(x, \tilde{x}) \dfn \begin{cases}
            h(x) - h(\tilde{x}) - \inner{\nabla h(\tilde{x})}{x - \tilde{x}} & \text{if}\; \tilde x\in \intdom{h} \\
            \infty & \text{otherwise.}
        \end{cases}%
    \end{equation*}}%
\end{definition}
\vspace{-0.5\baselineskip}
We say that the Bregman distance $D_h$ is \emph{induced by} the \emph{kernel function} $h$.
Examples include
\begin{inlinelist}
    \item the \rev[quadTerm]{}{square of} Euclidean distance with $h(x) = \norm{x}^2$; and
    \item the Kullback-Leibler divergence with $h(p) = \sum_i p_i\ln p_i$ being the negative entropy function.
\end{inlinelist}
The \emph{mirror-descent method} \cite{beck2003mirror} generalizes \rev{}{the} classical gradient descent method, which, when applied to \eqref{eq: reg_loss}, can be written as 
\begin{equation*}
    \theta^{k+1} = \argmin_{\theta} \regLoss(\theta^k) + \langle \nabla \regLoss(\theta^k), \theta - \theta^k \rangle + \tfrac{1}{2} \Vert \theta - \theta^k \Vert^2.
\end{equation*}
The mirror descent method replaces the quadratic term $\tfrac{1}{2}\Vert \theta - \theta^k \Vert^2$ by some Bregman distance $D_h(\theta, \theta^k)$.
We now show that \EMpp (cf. \cref{alg:em-conceptual}) can be interpreted as a 
mirror descent method for solving the problem
\eqref{eq: reg_loss}.
\begin{proposition}[Mirror descent interpretation] \label{prop:em_bregman}
    Consider the iterates $\left\{ \theta^{k} \right\}_{k \in \N}$ generated by \EMpp (cf. \cref{alg:em-conceptual}).
    Under \cref{asp: subsys_model,asp: regularizers,asp: smoothness} we have
    \begin{equation} \label{eq:bpg-update}
        \theta^{k+1} = \argmin_{\theta} \regLoss(\theta^k) + \inner{\nabla \regLoss(\theta^k)}{\theta - \theta^k} + D_{\kernel[k]}\left(\theta, \theta^k \right)
    \end{equation}
    for all $k \in \N$, where $\kernel[k]$ is a strictly convex function
    \begin{equation} \label{eq:bregman-function}
        \begin{aligned}
            \kernel[k](\theta) & \dfn
            \Qi[1][k]{\dTheta} + \Qi[2][k]{\dBeta} + \reg(\theta).
        \end{aligned}
    \end{equation}
\end{proposition}
\vspace{-\baselineskip}
\begin{pf}
    \Cref{alg:em-conceptual} computes $\theta^{k+1} = \argmin_{\theta} \regSurrogate[k](\theta)$.
    From \eqref{eq: Q_equal_func}, we know that \(
        \regSurrogate[k](\theta^k) = \regLoss (\theta^k)
    \).
    It therefore remains to show that
    \begin{equation} \label{eq:main-eq-to-show}
        \regSurrogate[k](\theta) = \regSurrogate[k](\theta^k) 
        + \inner{\nabla \regLoss({\theta^k})}{\theta - \theta^k} 
        + D_{\kernel[k]} \big(\theta, {\theta^k}\big).
    \end{equation}
    Under \cref{asp: smoothness},
    the functions $\regSurrogate[k]$, $\regLoss$ and $\reg$ are differentiable.
    Consequently, \cref{lem: directional_derivative} implies that 
    $\nabla \regLoss(\thk) = \nabla \regSurrogate[k](\thk)$, and 
    by definition of $\regSurrogate[k]$ we have
    \begin{align*}
        \inner{\nabla \regLoss({\theta^k})}{&\theta - \thk} 
        = \inner{\nabla \Qi[1][k]{\dTheta^k}}{\dTheta - \dTheta^k} \\
          & + \inner{\nabla \Qi[2][k]{\dBeta^k}}{\dBeta - \dBeta^k} 
            + \inner{\nabla \reg(\thk)}{\theta - \thk}. 
    \end{align*}
    Thus, we obtain that
    \begin{align*}
        \regSurrogate[k](\theta) - &\regSurrogate[k](\theta^k) - 
            \inner{\nabla \regLoss(\theta^k)}{\theta - \theta^k} \\
        = & D_{\Qi[1][k]{}}(\dTheta, \dTheta^k) + D_{\Qi[2][k]{}}(\dBeta, \dBeta^k) + D_{\reg}(\theta, \thk)\\
        = & D_{\kernel[k]}\left(\theta, \theta^k \right),
    \end{align*}
    where the last step follows from \cite[Prop. 3.5]{bauschke1997legendre}.
    This proves \eqref{eq:main-eq-to-show} and concludes the proof.
\qed\end{pf}
\vspace{-\baselineskip}
\revnew{2}{}{%
We show that the kernel functions $\kernel[k](\theta)$ in \eqref{eq:bregman-function}
are locally strongly convex on any compact set for all $k \in \N$.
\begin{lemma}[Local strong convexity]\label{lem: local_strong_convexity_of_kernels}
    Under \cref{asp: subsys_model,asp: regularizers,asp: smoothness},
    on any compact set $\Omega_{\theta} \subseteq \mathcal{X}_{\theta}$,
    there exists an $\varepsilon > 0$, such that
    \begin{equation}\label{eq: local_strong_convexity_of_kernels}
        \textstyle \inner{\nabla \kernel[k](\theta) - \nabla \kernel[k](\tilde{\theta})}{\theta - \tilde{\theta}} \geq \varepsilon \Vert \theta - \tilde{\theta} \Vert^2
    \end{equation}
    for all $\theta, \tilde{\theta} \in \Omega_{\theta}$ and for all $k \in \N$.
\end{lemma}
\begin{pf}
    See \cref{prf: local_strong_convexity_of_kernels}.
\qed\end{pf}
}

The convergence analysis of mirror descent-like methods typically relies on a descent lemma, see e.g., \cite[Lemma 1]{bauschke_descent_2017}, \cite[Lemma 2.1]{bolte2018first}. 
The following lemma describes a similar relation for \EMpp with respect to the variable kernel Bregman distance from \cref{prop:em_bregman}.

\begin{lemma}[Descent lemma] \label{lem:descent}
    Consider the iterates $\left\{ \theta^{k} \right\}_{k \in \N}$ generated by \EMpp (cf. \cref{alg:em-conceptual}).
    Under \cref{asp: subsys_model,asp: regularizers,asp: smoothness}, we have for all $k \in \N$ that
    \begin{equation*}
        \regLoss(\theta^{k+1}) \leq \regLoss(\theta^k) + \inner{\nabla \regLoss(\theta^k)}{\theta^{k+1} - \theta^k} + D_{\kernel[k]}\left(\theta^{k+1}, \theta^k \right).
    \end{equation*}
\end{lemma}
\begin{pf}
    Immediate from \eqref{eq:main-eq-to-show}, 
    \eqref{eq: Q_pos}, and \eqref{eq: Q_equal_func}.
\qed\end{pf}
\vspace{-\baselineskip}
We emphasize that the kernel function $\kernel[k]$ in \cref{eq:bpg-update,lem:descent} depends on the iterate $\thk$.
In contrast to the analyses of mirror-descent like methods in \cite{kunstner2021homeomorphic,bolte2018first},
where the kernel function is iterate-invariant,
this dependency introduces a significant challenge in showing asymptotic convergence.
Akin to \cite{bolte2018first}, which proves global \emph{sequential} convergence of so-called \textit{gradient-like descent sequences} under a KL property,
we show in the next lemma that the iterates of \EMpp constitute such a sequence, 
despite the variable kernel functions.
\begin{lemma}\label{lem: gradient_like_descent_sequence}
    Under \cref{asp: subsys_model,asp: regularizers,asp: smoothness},
    the sequence $\{\thk\}$ generated by \EMpp (cf. \cref{alg:em-conceptual}) is a \emph{gradient-like descent sequence} for $\regLoss$, i.e.,
    \begin{lemitem}
        \item There exists a positive scalar $\rho_1$ such that 
            $\rho_1 \norm{\thkp - \thk}^2 \leq \regLoss(\thk) - \regLoss(\thkp)$ for all $k \in \N$;
            \label{item: Sufficient_descrease}
        \item There exists a positive scalar $\rho_2$ such that 
            $ \norm{\nabla \regLoss(\thkp)} \leq \rho_2 \norm{\thkp - \thk}$ for all $k \in \N$;
            \label{item: gradient_lower_bound}
        \item Let $\bar{\theta}$ be a limit point of a subsequence $\{\thk\}_{k \in \mathcal{K}}$, 
            then $\lim\sup_{k\in \mathcal{K}\subset \N} \regLoss(\thk) \leq \regLoss(\bar{\theta})$.
            \label{item: limsup_regLoss}
    \end{lemitem}%
\end{lemma}
\vspace{-\baselineskip}
\begin{pf}
    \rev[gradient-like-descent-sequence]{
    The optimality conditions defining \eqref{eq:bpg-update} yield
\begin{equation} \label{eq:bpg-update-optimality-relation}
    \nabla \regLoss(\theta^k) = \nabla \kernel[k](\theta^k) - \nabla \kernel[k](\theta^{k+1}).
\end{equation}
By \cref{lem:descent} and \cref{def:bregman-distance} this yields
\begin{align}
    \regLoss(\theta^{k+1}) - \regLoss(\theta^k)
    \leq &\;\nonumber \langle \nabla \kernel[k](\theta^k) - \nabla \kernel[k](\theta^{k+1}), \theta^{k+1} - \theta^k \rangle \\
        &\;\;\nonumber + D_{\kernel[k]}(\theta^{k+1}, \theta^k)\\
    = &\;\;- D_{\kernel[k]}(\theta^k, \theta^{k+1}).\label{prf:bregman-distance-convergence:ineq-1}
\end{align}
\revnew{2}{%
By \eqref{eq:bregman-function}, we have for all $\theta, \tilde{\theta} \in \mathcal{X}_\theta$,
\vspace{-\baselineskip}%
\vspace{-\baselineskip}%
{\small%
\begin{align*}
    \inner{\nabla \kernel[k](\theta) - \nabla \kernel[k](\tilde{\theta})}{\theta - \tilde{\theta}}
    = &\; \inner{\nabla \Qi[1][k]{\dTheta} - \nabla \Qi[1][k]{\tilde{\dTheta}}}{\dTheta - \tilde{\dTheta}} \\
      &\; +\inner{\nabla \Qi[2][k]{\dBeta} - \nabla \Qi[2][k]{\tilde{\dBeta}}}{\dBeta - \tilde{\dBeta}} \\
      &\; + \inner{\nabla \reg(\theta) - \nabla \reg(\tilde{\theta})}{\theta - \tilde{\theta}} \\
    \geq  &\; \inner{\nabla \reg(\theta) - \nabla \reg(\tilde{\theta})}{\theta - \tilde{\theta}},
\end{align*}}%
where the inequality follows from the convexity of $\Qi[1][k]{}$ and $\Qi[2][k]{}$.%
}{%
Let $\hbar(\theta) = \frac{\varepsilon}{2} \norm{\theta}^2$ with $\varepsilon > 0$ satisfing \cref{lem: local_strong_convexity_of_kernels}.}
By \revnew{2}{}{\cref{lem: local_strong_convexity_of_kernels} and} \cite[Proposition 1.1]{doi:10.1137/16M1099546},
\revnew{2}{
\(
    D_{\kernel[k]}(\theta, \theta^k) \geq D_{\reg}(\theta, \theta^k).    
\)
}{%
we have that
\(
    D_{\kernel[k]}(\theta, \theta^k) \geq D_{\hbar}(\theta, \theta^k).    
\)}
Combined with 
\eqref{prf:bregman-distance-convergence:ineq-1}, 
this yields 
\revnew{2}{\(
    \regLoss(\thk) - \regLoss(\thkp) \geq D_{\reg}(\thk, \thkp).
\)}{%
\(
    \regLoss(\thk) - \regLoss(\thkp) \geq D_{\hbar}(\thk, \thkp).
\)}
\revnew{2}{
By \cref{asp: regularizers} and \cref{lem: general_compactness},
the iterates $\thk$ remain in a compact set $\Omega \subset \mathcal{X}_\theta$.
Since $\nabla^2 \reg(\theta)$ is continuous and $\nabla^2 \reg(\theta) \succ 0$ for all $\theta$, as assumed in \cref{asp: smoothness},
the function $\reg$ is locally strongly convex on $\Omega$.
Hence,
there exists a constant $\rho_1 > 0$, such that 
\(
    \regLoss(\thk) - \regLoss(\thkp) \geq \rho_1\norm{\thk - \thkp}^2,
\)
proving \crefpart{lem: gradient_like_descent_sequence}{item: Sufficient_descrease}.}{%
Since $D_{\hbar}(\thk, \thkp) = \frac{\varepsilon}{2} \norm{\thk - \thkp}^2$,
the claim \crefpart{lem: gradient_like_descent_sequence}{item: Sufficient_descrease} follows with $\rho_1 = \frac{\varepsilon}{2}$.
}
As for the second claim, we have by \eqref{eq:bpg-update-optimality-relation} that
\(
    \nabla \regLoss (\thkp) = \nabla \regLoss (\thkp) - 
        \nabla \regLoss (\thk) - \nabla \kernel[k](\thkp) + \nabla \kernel[k](\thk),
\)
and hence
\begin{align}\label{eq: upper_grad_L}
    \norm{\nabla \regLoss (\thkp)}
    \leq &\; \phantom{+} \norm{\nabla \regLoss (\thkp) - \nabla \regLoss (\thk)} \nonumber \\
     & + \norm{\nabla \kernel[k](\thkp) - \nabla \kernel[k](\thk)} .
\end{align}
By \cref{asp: smoothness}, 
the loss function is composed of Lipschitz smooth functions on the compact set $\Omega$ that contains all iterates $\left\{ \theta^k \right\}$.
Hence, there exists $L_\Omega > 0$ such that
\(
    \norm{\nabla \regLoss(\thkp) - \nabla \regLoss(\thk)} \leq L_{\Omega} \norm{\thkp - \thk}.    
\)
Moreover, by \eqref{eq:bregman-function}, we can bound the second term of \eqref{eq: upper_grad_L} by
\vspace{-\baselineskip}
\vspace{-\baselineskip}
{\small%
\begin{align*}
    &\norm{\nabla \kernel[k](\thkp) - \nabla \kernel[k](\thk)} 
    \leq \norm{\nabla \Qi[1][k]{\dTheta^{k+1}} - \nabla \Qi[1][k]{\dTheta^k}} \\
    &\;\; + \norm{\nabla \Qi[2][k]{\dBeta^{k+1}} - \nabla \Qi[2][k]{\dBeta^{k}}}
    + \norm{\nabla \reg(\thkp) - \nabla \reg(\thk)}.
\end{align*}
}%
Using \eqref{eq:switch-surrogate}
and the Cauchy-Schwarz inequality,
we further bound
\vspace{-\baselineskip}%
{%
    \small%
    \begin{equation}\label{eq: diff_grad_Q2}
        \begin{aligned}
            &\; \norm{\nabla \Qi[1][k]{\dTheta^{k+1}} - \nabla \Qi[1][k]{\dTheta^k}} \\
            &\; \quad \quad \leq 
            \big\lVert\Pi(\thk)\big\rVert
            \big\lVert\nabla \Psi(\dTheta^{k+1}) - \nabla \Psi(\dTheta^{k})\big\rVert.
        \end{aligned}
    \end{equation}
}%
Recall that $\Pi(\thk) \in [0, 1]^{d^{T+1}}$,
\rev{}{and that} the gradient of \rev{function}{} $\lse$ is Lipschitz continuous on the compact set $\Omega \in \mathcal{X}_{\theta}$.
Consequently, 
$\Psi$ has Lipschitz continuous gradients on $\Omega$,
and
there exists a constant $M_{1, \Omega} > 0$ such that 
$\norm{\nabla \Qi[1][k]{\dTheta^{k+1}} - \nabla \Qi[1][k]{\dTheta^k}} \leq M_{1, \Omega}\norm{\dTheta^{k+1} - \dTheta^{k}}$ for all $\thk\in \mathcal{X}_{\theta}$.
\rev[newCompact]{
In a similar way,
let $G, F_t, Z_t \in \re^{d^{T+1}}$
with elements
$G_{\modeSeq}(\dBeta)        = \sum_{t = 0}^{T-1} g_t(\beta_{\xi_{t+1}}),$
$Z_{t, \modeSeq}(\dBeta)     = \ell_t(\beta_{\xi_{t+1}}),                $ and
$F_{t, \modeSeq}(\dBeta^k)   = f'(\ell_t(\beta^k_{\xi_{t+1}}))           $
for all $\modeSeq \in \Xi^{T+1}$ and $t \in \N_{[0, T-1]}$.
Hence, $ \nabla \hat{\Phi}^k(\dBeta) = \sum_{t= 0}^{T-1} F_{t}(\dBeta^k) \nabla Z_{t}(\dBeta) + \nabla G(\dBeta)$.
}{
    Similarly,
    let 
    $Z_{t, \modeSeq}(\dBeta)        = G_t(\beta_{\xi_{t+1}}),$ and 
    $F_{t, \modeSeq}(\dBeta^k)      = \nabla f\big(G_t(\beta^k_{\xi_{t+1}})\big)$ 
    for all $\modeSeq \in \Xi^{T+1}$ and $t \in \N_{[0, T-1]}$.
    Recalling \eqref{eq:subsys-surrogate} and applying the chain rule to each entry $\hat{\Phi}_{\modeSeq}^k(\dBeta)$ given in \eqref{eq: def_elem_surrogate_subsys},
    we obtain
}
\rev{Using \eqref{eq:subsys-surrogate}, the triangle and Cauchy-Schwarz inequality,
we further bound}{}%
\vspace{-\baselineskip}%
\vspace{-\baselineskip}%
{%
    \small%
    \rev{
    \begin{equation}
        \begin{aligned}
            &\; \norm{\nabla \Qi[2][k]{\dBeta^{k+1}} - \nabla \Qi[2][k]{\dBeta^{k}}} \\
            \leq &\; \big\lVert
                \Pi(\thk) 
                \tlsum_{t=0}^{T-1} F_{t}(\dBeta^k) \left(\nabla Z_t(\dBeta^{k+1}) - \nabla Z_t(\dBeta^{k})\right)
                \big\rVert \\
            &\; + \big\lVert\Pi(\thk) \left(
                  \nabla G(\dBeta^{k+1})
                - \nabla G(\dBeta^{k})\right)
                \big\rVert \\
            \leq &\; \tlsum_{t=0}^{T-1} \norm{\Pi(\thk) F_t(\dBeta^k)} \norm{
                \nabla Z_t(\dBeta^{k+1})
                - \nabla Z_t^1(\dBeta^{k})
                } \\ 
            &\; + \big\lVert \Pi(\thk) \big\rVert \big\lVert
                    \nabla G_t(\dBeta^{k+1})
                - \nabla G_t(\dBeta^{k})
                \big\rVert.
        \end{aligned}
    \end{equation}
    }{
        \begin{align}
            &\; \norm{\nabla \Qi[2][k]{\dBeta^{k+1}} - \nabla \Qi[2][k]{\dBeta^{k}}} \label{eq: diff_grad_Q1}\\
            = &\; \big\lVert\tlsum_{\modeSeq \in \Xi^{T+1}}\big(\nabla \hat{\Phi}^k_{\modeSeq}(\dBeta^{k+1}) - \nabla \hat{\Phi}^k_{\modeSeq}(\dBeta^{k})\big)\Pi_{\modeSeq}(\thk)\big\rVert \nonumber\\
            = &\; \big\lVert\tlsum_{\modeSeq \in \Xi^{T+1}}\!\tlsum_{t= 0}^{T-1} \big(\nabla Z_{t, \modeSeq}(\dBeta^{k+1}) - \nabla Z_{t, \modeSeq}(\dBeta^{k}) \big)F_{t, \modeSeq}(\dBeta^k)\Pi_{\modeSeq}(\thk)\big\rVert \nonumber\\
            \leq &\; \tlsum_{\modeSeq \in \Xi^{T+1}}\!\tlsum_{t= 0}^{T-1} 
                \big\lVert \nabla Z_{t, \modeSeq}(\dBeta^{k+1}) - \nabla Z_{t, \modeSeq}(\dBeta^{k}) \big\rVert 
                \big\lVert F_{t, \modeSeq}(\dBeta^k)\Pi_{\modeSeq}(\thk)\big\rVert \nonumber
        \end{align}
        where the inequality follows from the triangle and Cauchy-Schwarz inequality.
    }%
}%
As assumed in \Cref{asp: smoothness}, 
\rev{$Z_t$, $G$}{$G_t$}
have Lipschitz continuous gradients on $\Omega$.
Since \rev{$\Pi(\thk) \!\in\! [0, 1]^{d^{T+1}}$}{$\Pi_{\modeSeq}(\thk) \!\in\! [0, 1]$}
and \rev{$F_t(\dBeta^k) \!\in\! [0, \bar{u}]^{d^{T+1}}$}{$F_{t, \modeSeq}(\dBeta^k) \!\in\! [0, \bar{u}]^M$},
following from \cref{asp: subsys_model},
\eqref{eq: diff_grad_Q1} implies that
there exists $M_{2, \Omega} \!>\! 0$ such that 
$\norm{\nabla \Qi[2][k]{\dBeta^{k+1}} \!-\! \nabla \Qi[2][k]{\dBeta^{k}}} 
    \leq M_{2, \Omega} \norm{\dBeta^{k+1} \!-\! \dBeta^{k}}$ 
for all $\thk \!\in\! \mathcal{X}_\theta$.
Since $\reg$ is of class $\mathcal{C}^2$ by \cref{asp: smoothness},
it has Lipschitz continuous gradients on $\Omega$.
Thus, there exists a constant $\rho_2 > 0$ such that
$\norm{\nabla \regLoss (\thkp)} \leq \rho_2 \norm{\thkp - \thk}$.
This proves \crefpart{lem: gradient_like_descent_sequence}{item: gradient_lower_bound}.
Finally, \crefpart{lem: gradient_like_descent_sequence}{item: limsup_regLoss} follows directly from the continuity of $\regLoss$ on the set $\Omega \subset \mathcal{X}_\theta$.

    }{See \cref{prf: gradient_like_descent_sequence}}
\qed\end{pf}

\begin{definition}[KL property \cite{bolte2007lojasiewicz,lojasiewicz1963propriete}] \label{def:kl}
    A proper and lower semicontinuous function $\eta : \R^n \to \re\cup \{\infty\}$ satisfies the Kurdyka-\L ojasiewicz (KL) property at $\bar \theta \in \mathcal{X}_\theta$ if there exists a concave KL function $\psi : [0, b] \to [0, +\infty)$ with $b > 0$ and a neighborhood $U_{\bar \theta}$ such that
    \begin{lemitem}
        \item $\psi(0) = 0$;
        \item $\psi \in \mathcal{C}^1$, with $\psi' > 0$ on $(0, b)$;
        \item $\forall \theta \in U_{\bar \theta}$: if $\eta(\bar \theta) < \eta(\theta) < \eta(\bar \theta) + b$, then \\$\psi'(\eta(\theta) - \eta(\bar \theta)) \mathbf{dist}(0, \partial \eta(\theta)) \geq 1$.
    \end{lemitem}
\end{definition}
The KL property is a mild requirement
that holds for real-analytic and semi-algebraic functions, and subanalytic functions that are continuous on their domain \cite{bolte2007lojasiewicz}.
In fact, for these first two classes of functions, the KL function of $\regLoss$ can be taken of the form $\psi(s) = c s^{1 - \alpha}$,
with $c > 0, \alpha \in [0, 1)$ \cite{attouch_proximal_2010}.
We highlight that often \rev[changeNotationKL]{$f, \ell_{t}, g_{t}$}{$f, G_{t}$} are real-analytic for $t \in \N_{[0, T-1]}$.
Since sums, products and compositions of real-analytic functions are real-analytic, $\regLoss$ trivially satisfies the KL property in such cases.  
The next theorem shows that if $\regLoss$ satisfies the KL property, the whole sequence of iterates converges with a rate depending on the specific form of the KL function. 

\begin{theorem}[Global convergence] \label{prop:convergence-kl}
    Consider the iterates $\{\theta^{k} \}_{k \in \N}$ generated by \EMpp (cf. \cref{alg:em-conceptual}).
    If $\regLoss$ satisfies the KL property, and if \cref{asp: subsys_model,asp: regularizers,asp: smoothness} hold, 
    then $\{\theta^{k} \}_{k \in \N}$ converges to a stationary point, i.e.,
    \begin{equation*}
        \lim_{k \to \infty} \theta^k = \theta^\star, \quad \text{with } \nabla \regLoss(\theta^\star) = 0.
    \end{equation*}
    If, additionally, the KL function $\psi$ of $\regLoss$ is of the form $\psi(s) = c s^{1 - \alpha}$ with $c > 0, \alpha \in [0, 1)$, then
    \begin{lemitem}
        \item if $\alpha = 0$ then $\left\{ \theta^k \right\}$ converges to $\theta^\star$ in a \emph{finite} number of steps;
        \item if $\alpha \in (0, 1/2]$ then there exist $w > 0$ and $\tau \in [0, 1)$ such that
        \(
            \Vert \theta^k - \theta^\star \Vert \leq \omega \tau^k;
        \)
        \item if $\alpha \in (1/2, 1)$ then there exist $\omega > 0$ such that
        \(
            \Vert \theta^k - \theta^\star \Vert \leq \omega k^{-\frac{1 - \alpha}{2 \alpha - 1}}.
        \)
    \end{lemitem}
\end{theorem}

\begin{pf}
    By \cref{lem: gradient_like_descent_sequence}, $\left\{ \theta^k \right\}$ is a gradient-like descent sequence.
    The claims follow by \cite[Theorem 6.2 \& 6.3]{bolte2018first}.
\qed\end{pf}

\section{Efficient evaluation of the surrogate} \label{sec: e_step}
An effective implementation of \EMpp (cf. \cref{alg:em-conceptual}) requires efficient evaluation of the surrogate function $\surrogate$.
Despite the exponential complexity inherent in the definition of $\surrogate$ in \eqref{eq: surrogate},
its evaluation can be accomplished efficiently using the following representation.

\begin{proposition}\label{prop:separable-explicit}
	Given a trajectory $\ve{y}$
	of system \eqref{eq: discrete_latent} that is initialized at $z_0$,
    the terms of $\surrogate(\theta) = \Qi[0][k]{} + \Qi[1][k]{\dTheta} + \Qi[2][k]{\dBeta}$ as defined in \eqref{eq: surrogate} can be represented as%
	\begin{subequations}\label{eq: surrogate-explicit}
		\vspace{-\baselineskip}
		\vspace{-\baselineskip}
		{%
			\small%
			\begin{align}
				\label{eq: initial-surrogate-explicit}
				\Qi[0][k]{} \dfn & \;
            \!-\!\tlsum_{i \in \Xi}\!\pi_{0}^i(\thk)\!
				\big[
					\ln p(\xi_0\!=\!i, z_0)\!+\!T \ln C
					\big]
				\!+\! c_{\thk},                \\
				\Qi[1][k]{\dTheta}
				\dfn             & \;
				\tlsum_{t=0}^{T-1} \tlsum_{i,j\in \Xi} \pi_{t}^{i,j}(\thk)
				\Big(
                    \lse\big(\Theta_{i}^\top z_t\big) - \Theta_{i, j}^\top z_t
				\Big), \label{eq:switch-surrogate-explicit}\\
                \ifshowboth
                    \rev{\Qi[2][k]{\dBeta}
                    \dfn}{}             
                    & \; \rev{\hat{c}_{\thk} {+}
                    \tlsum_{t=0}^{T-1}\!\tlsum_{i \in \Xi}\!\pi_{t+1}^i(\thk)
                    \Big(f'\big(
                    \ell_{t}(\beta_{i}^k)
                    \big)\, \ell_{t}(\beta_{i})
                    + g_{t}(\beta_{i})
                    \Big)}{}, \\
                \fi
                \rev{}{\Qi[2][k]{\dBeta}\dfn}   
                & \; \rev{}{
                    \hat{c}_{\thk} \!{+}\!
                    \tlsum_{t=0}^{T-1}\!\tlsum_{i \in \Xi}\!\pi_{t+1}^i\!(\thk)
                    \inner{\nabla\! f\big(
                    G_{t}(\beta_{i}^k)
                    \big)}{G_{t}(\beta_{i})
                    }
                },  
                  \label{eq:subsys-surrogate-explicit}                  
			\end{align}
		}%
	\end{subequations}
    where we have used the shorthands
	\begin{align*}
        \vspace{-\baselineskip}
        \pi_{t}^{i}(\thk)   & \dfn p(\xi_t=i \mid \ve{y}, z_0\midsc \thk);                \\
        \pi_{t}^{i,j}(\thk) & \dfn p(\xi_{t}=i, \xi_{t+1}=j \mid \ve{y}, z_0\midsc \thk);
    \end{align*} 
    \rev{and $
        \hat{c}_{\thk} = \sum_{t=0}^{T-1} \sum_{i \in \Xi} 
         - \pi_{t+1}^i(\thk) f '(\ell_{t}(\beta_i^k))\, \ell_{t}(\beta_{i}^k) + \\ \pi_{t+1}^i(\thk) f(\ell_{t}(\beta^k_i)) 
    $.\\}{and
        \(
            \hat{c}_{\thk} = \sum_{t=0}^{T-1} \sum_{i \in \Xi} 
         - \pi_{t+1}^i(\thk)
         \inner{\nabla f(G_{t}(\beta_i^k))}{G_{t}(\beta_{i}^k)} \\+ \pi_{t+1}^i(\thk) f(G_{t}(\beta^k_i)).
        \)
    }
\end{proposition}
\vspace{-\baselineskip}
\begin{pf}
We recall that each entry $\Pi_{\modeSeq}(\thk) = p(\modeSeq \mid \stateSeq, z_0 \midsc \thk) = p(\xi_0, \dots, \xi_{T} \mid \stateSeq, z_0 \midsc \thk)$ for some $\modeSeq \in \Xi^{T+1}$.
A distribution $p(\xi_0 \mid \stateSeq, z_0 \midsc \thk)$ can be obtained by \emph{marginalizing out} the random variables $\xi_1, \dots \xi_T$, i.e.,
\vspace{-\baselineskip}
\vspace{-\baselineskip}
{%
    \small%
    \begin{equation*}
    p(\xi_0 \mid \stateSeq, z_0 \midsc \thk) = \tlsum_{(\xi_1,\dots,\xi_T) \in \Xi^{T}} p(\xi_0, \dots, \xi_{T} \mid \stateSeq, z_0 \midsc \thk).
    \end{equation*}
}%
Considering $\Qi[0][k]{}$ from \eqref{eq: initial-surrogate},
we can apply this marginalization since $\ve{c}$ only depends on the latent variable $\xi_0$:
\vspace{-\baselineskip}
\vspace{-\baselineskip}
{%
    \small%
    \begin{align*}
        \Qi[0][k]{} 
        =&\; c_{\thk} - \tlsum_{\modeSeq \in \Xi^{T+1}}\Pi_{\modeSeq}(\thk) \ve{c}_{\modeSeq} \\
        =&\; c_{\thk} - \quad\tlsum_{\mathclap{(\xi_0, \dots, \xi_{T}) \in \Xi^{T+1}}}p(\xi_0, \dots, \xi_T\mid \stateSeq, z_0 \midsc \thk) [\ln p(\xi_0, z_0) + T\ln C]\\
        =&\; c_{\thk} - \;\;\tlsum_{\xi_0 \in \Xi}p(\xi_0 \mid \stateSeq, z_0 \midsc \thk) [\ln p(\xi_0, z_0) + T\ln C]
    \end{align*} 
}%
This yields \eqref{eq: initial-surrogate-explicit}.
We apply the same approach to $\Qi[1][k]{}, \Qi[2][k]{}$ from \cref{eq:subsys-surrogate,eq:switch-surrogate}. 
In particular, the terms of $\Psi_{\modeSeq}(\dTheta)$ in \eqref{eq:Lt2-def} only depends on $\xi_{t}, \xi_{t+1}$,
and the terms of $\hat{\Phi}_{\modeSeq}(\dBeta)$ in \eqref{eq: def_elem_surrogate_subsys} only depends on $\xi_{t+1}$.
By marginalizing out the other latent variables, we obtain \cref{eq:subsys-surrogate-explicit,eq:switch-surrogate-explicit}. 
\qed\end{pf}

Since the switching probability of $\xi_{t+1}$ depends on $\xi_t$,
as evident from \eqref{eq: switch_model},
we utilize the forward-backward algorithm \cite{rabiner1989tutorial} to compute the distributions 
$\pi_{t+1}^i(\thk)$,
$\pi_{t}^{i, j}(\thk)$.
The standard forward-backward algorithm computes the posterior distribution of latent variables in a hidden Markov model.
Unlike hidden Markov models where the 
observations $y_t$ are conditionally independent given the latent variable $\xi_t$,
the dynamics of $y_t$ described in \eqref{eq: subsys_model} necessitates some modifications to the original forward-backward scheme, 
as we now formalize.
\begin{proposition} \label{prop: forward_backward}
    Let $
        \rho_t(\xi_{t}, \xi_{t+1} \midsc \thk) 
        \dfn p(\xi_{t}, \xi_{t+1}, \ve{y}_{t+1:T} \mid \ve{y}_{0:t}, z_0\midsc \thk)$.
    The posterior distributions expressed as
    \begin{subequations}\label{eq: forward_backward_summary}
        \begin{align}
            \pi_{t+1}^i(\thk)  = &\; \textstyle \sum_{j \in \Xi} \pi^{j, i}_{t}(\thk) \label{eq: forward_backward} \\
            \pi_{t}^{i, j}(\thk) = &\;
                \frac{
                    \rho_t(\xi_{t} = i, \xi_{t+1} = j \midsc \thk)
                }{
                    \sum_{m, n \in \Xi^2} \rho_t(\xi_{t} = m, \xi_{t+1} = n\midsc \thk)
                } \label{eq: forward_backward_joint} 
        \end{align}
    can be computed with
    \vspace{-\baselineskip}
    \vspace{-\baselineskip}
    {%
        \small%
        \begin{equation}\label{eq: rho_joint}
            \rho_t(\xi_{t}, \xi_{t+1}\midsc \thk)
            = \zeta_t(\xi_{t+1}\midsc \thk)
            \, p(\xi_{t+1} \mid z_t, \xi_t\midsc \dTheta^k)
            \, \alpha_t(\xi_t\midsc \thk),
        \end{equation}
    }%
    where 
    $\alpha_t(\xi_t) \dfn p(\xi_t\mid \ve{y}_{0:t}, z_0)$, 
    $\zeta_t(\xi_{t+1}) \dfn p(\ve{y}_{t+1:T} \mid \xi_{t+1}, z_t)$.
    Let $q_{t-1}(\xi_t) \dfn p(\xi_t\mid \ve{y}_{0: t-1}, z_0)$.
    The probability $\alpha_t(\xi_t)$ and the likelihood $\zeta_t(\xi_{t+1})$
    are updated recursively for all $t\in \N_{[1, T]}$ through
    \vspace{-\baselineskip}
    \vspace{-\baselineskip}
    {%
        \small%
        \begin{alignat}{2}
            q_{t-1}(\xi_t \midsc \thk) 
            = \!&\;
            \tlsum_{\mathclap{\xi_{t-1} \in \Xi}} p(\xi_t \!\mid z_{t-1}, \xi_{t-1}\!\midsc \dTheta^k\!)\alpha_{t-1}(\xi_{t-1}\!\midsc \thk\!), 
            \label{eq: p_predict}
            \\
            \alpha_t(\xi_t\midsc \thk)
            = &\; \frac{
                p(y_t \mid z_{t-1}, \xi_t\midsc \dBeta^k)\, q_{t-1}(\xi_t \midsc \thk)
            }{
                \sum_{\xi_{t} \in \Xi} p(y_t \mid z_{t-1}, \xi_t\midsc \dBeta^k) q_{t-1}(\xi_t \!\midsc\! \thk)
            },  \label{eq: forward_filtering} \\
            \zeta_{t-1}(\xi_{t}\midsc \thk) 
            =&\;  \textstyle \label{eq: backward_filtering}
            \sum_{\xi_{t+1}\in \Xi} \zeta_t(\xi_{t+1}\midsc \thk) \times \\
               &\; \quad p(\xi_{t+1}\mid z_t, \xi_t\midsc \dTheta^k) \, p(y_{t} \mid z_{t-1}, \xi_t\midsc \dBeta^k) \nonumber
        \end{alignat}
    }%
    \end{subequations}
    with initialization $\alpha_0(\xi_0)$, $\zeta_T(\xi_{T+1}) = 1$.
\end{proposition}
\begin{pf}
    \rev[forward-backward]{
        For brevity, we omit the dependence on $\thk$ as it is fixed for all computations involved in \eqref{eq: forward_backward_summary}.
Recall the definition of $\pi_{t+1}^i$, $\pi_{t}^{i, j}$ in \cref{prop:separable}, 
equation \eqref{eq: forward_backward_joint} and \eqref{eq: forward_backward}
directly follow from the definition of the conditional distribution and the marginalization.
The distribution $\rho_t(\xi_t, \xi_{t+1})$ can be factorized as
\vspace{-\baselineskip}
\vspace{-0.5\baselineskip}
{%
    \small%
    \begin{equation}\label{eq: factorization_rho}
        \begin{aligned}
            \rho_t(\xi_t, \xi_{t+1})
        = &\;  p(\xi_t, \xi_{t+1}, \ve{y}_{t+1:T} \mid \ve{y}_{0:t}, z_0) \\
        = &\; p(\ve{y}_{t+1:T} \mid \xi_{t+1}, \xi_t, \ve{y}_{0:t}, z_0) \times \\
        &\; p(\xi_{t+1}\mid \xi_t, \ve{y}_{0:t}, z_0)
        \, p(\xi_t\mid \ve{y}_{0:t}, z_0).
        \end{aligned}
    \end{equation} 
}%
Recall from \eqref{eq: def_history} that $z_t$ is the trajectory history. 
Given $\xi_t, z_t$, the distribution of $\xi_{t+1}$ is independent of $\ve{y}_{0:t-t_y}$ (cf.\,\eqref{eq: switch_model}).
Hence, $p(\xi_{t+1}\mid \xi_t, \ve{y}_{0:t}, z_0) = p(\xi_{t+1} \mid z_t, \xi_t)$.
Likewise, the future trajectory $\ve{y}_{t+1:T}$ is independent of the past trajectory $\ve{y}_{0:t-t_y}$.
Thus, $p(\ve{y}_{t+1:T} \mid \xi_{t+1}, \xi_t, \ve{y}_{0:t}, z_0)$ simplifies to $p(\ve{y}_{t+1:T} \mid \xi_{t+1}, z_t)$.
By definition of $\alpha_t(\xi_t)$ and $\zeta_t(\xi_{t+1})$, \eqref{eq: factorization_rho} is equivalent to \eqref{eq: rho_joint}.

By definition of $\alpha_t(\xi_t)$ and $q_{t-1}(\xi_t)$,
\eqref{eq: p_predict} directly follows from the law of total probability.
Recall that 
\begin{equation}\label{eq: alpha_intermediate}
    \alpha_t(\xi_t) \dfn p(\xi_t\mid \ve{y}_{0:t}, z_0)
    =\frac{p(y_t, \xi_t \mid \ve{y}_{0: t-1}, z_0)}{p(y_t \mid \ve{y}_{0: t-1}, z_0)},
\end{equation}
where 
the denominator $p(y_t \mid \ve{y}_{0: t-1}, z_0) \!=\!\sum_{\xi_t\in \Xi} p(y_t, \xi_t \mid \ve{y}_{0: t-1}, z_0)$,
and the numerator
$p(y_t, \xi_t \mid \ve{y}_{0: t-1}) \!=\! p(y_t \mid \xi_t, \ve{y}_{0: t-1}, z_0)p(\xi_t \mid \ve{y}_{0: t-1}, z_0)$,
Given $\xi_{t}$ and $z_{t-1}$,
$y_{t}$ is independent of the trajectory $\ve{y}_{0: t-t_y-1}$.
Thus,
$p(y_t \mid \xi_t, \ve{y}_{0: t-1}, z_0)\! =\! p(y_t \mid z_{t-1}, \xi_t)$.
Plugging this equation and the definition of $q_{t-1}(\xi_t)$ in \eqref{eq: alpha_intermediate} yields \eqref{eq: forward_filtering}.

The likelihood $\zeta_{t-1}(\xi_{t})$ is obtained via marginalization:
\vspace{-\baselineskip}
\vspace{-0.5\baselineskip}
{%
    \small%
    \begin{equation}\label{eq: likelihood_to_be_proven}
        \zeta_{t-1}(\xi_{t})\!\dfn p(\ve{y}_{t:T}\! \mid\! \xi_{t}, z_{t-1})
        =\tlsum_{\mathclap{\xi_{t{+}1}\!\in \Xi}} p(\ve{y}_{t:T}, \xi_{t+1} \mid \xi_{t}, z_{t-1})
    \end{equation}  
}%
where $p(\ve{y}_{t:T}, \xi_{t+1} \mid \xi_{t}, z_{t-1})$
can be factorized as
\begin{align*}
    p(\ve{y}_{t:T}, \xi_{t+1} \mid \xi_{t}, z_{t-1})
    &= p(\ve{y}_{t+1:T} \mid \xi_{t+1}, \xi_{t}, y_t, z_{t-1}) \times \\
      & p(\xi_{t+1}\mid \xi_{t}, y_t, z_{t-1}) 
      p(y_{t} \mid \xi_{t}, z_{t-1}).
\end{align*}
As $z_t$ is a function of the trajectory history,
$p(\xi_{t+1}\mid \xi_{t}, y_t, z_{t-1}) = p(\xi_{t+1}\mid z_{t}, \xi_{t})$.
Applying the definition of $\zeta_t$ and \eqref{eq: discrete_latent} in 
\eqref{eq: likelihood_to_be_proven} yields \eqref{eq: backward_filtering}.

    }{See \cref{prf: forward_backward}.}
\qed\end{pf}
\begin{remark}[State-dependent switching]
    A special case arises when the switching probability only depends on the trajectory history,
    as in \eqref{eq: state_dependent_switching}, where the distribution
    $p(\xi_{t+1} \mid \xi_t, z_t\midsc \thk) = p(\xi_{t+1} \mid z_t\midsc \thk)$ for $t \in \N_{[0, T-1]}$.
    In this case, the posterior distribution
    \vspace{-\baselineskip}
    \vspace{-\baselineskip}
    {%
        \small%
        \begin{equation}\label{eq: posterior_state_dependent}
            \pi_{t+1}^i\!(\thk\!) \!=\! p(\xi_{t+1} {=} i\!\mid\! \ve{y}, z_0 \!\midsc\!\thk\!) 
                \!=\! \frac{p(y_{t+1}, \xi_{t+1} {=} i\mid z_t \!\midsc\! \thk\!)}{p(y_{t+1} \mid z_t\!\midsc \thk)},
        \end{equation} 
    }%
    with $p(y_{t+1}, \xi_{t+1}\mid z_t \midsc \thk) = p(y_{t+1}\mid \xi_{t+1}, z_t \midsc \thk)p(\xi_{t+1}\mid z_t \midsc \thk)$.
    The joint distribution $\pi_{t}^{i, j} (\thk) = p(\xi_{t} = i, \xi_{t+1} = j \mid \ve{y}, z_0 \midsc \thk) = p(\xi_{t} = i\mid \ve{y}, z_0 \midsc \thk)\,p(\xi_{t+1} = j \mid \ve{y}, z_0 \midsc \thk)$
    as a result of conditional independence. 
    The computational complexity of \eqref{eq: posterior_state_dependent}
    is $\mathcal{O}(dT)$.
    In contrast, the computation of \eqref{eq: forward_backward_summary} has a complexity of $\mathcal{O}(d^2 T)$.
    Leveraging the independence in \eqref{eq: posterior_state_dependent} is thus more efficient,
    especially when the number of modes $d$ is large.
    In addition,
    we highlight that \eqref{eq: posterior_state_dependent} enables parallel computation across all $\xi_t, t \in \N_{[1, T]}$.
    This parallelization capability significantly improves the scalability of the method compared to the recursive method \eqref{eq: forward_backward},
    in particular when the number of data points $T$ is large.
\end{remark}
\begin{remark}[Initialization]
    As discussed in \cref{prop: forward_backward},
    the recursion \eqref{eq: forward_backward_summary} requires initialization 
    $\alpha_0(\xi_0) = p(\xi_0 \mid y_0, z_0) = p(\xi_0 \mid z_0)$,
    which is necessary to start \EMpp (cf. \cref{alg:em-conceptual}).
    From the definition of $\alpha_0(\xi_0)$ and \eqref{eq: initial-surrogate-explicit}, 
    we observe that
    $\Qi[0][k]{} = -\sum_{i \in \Xi}\pi_{0}^i(\thk) [\ln \alpha_0(\xi_0 = i) + \ln p(z_0)] + c_{\thk}$
    is indeed a function of $\alpha_0(\xi_0)$.
    The weighted sum $H_0 \dfn -\sum_{i \in \Xi}\pi_{0}^i(\thk) [\ln \alpha_0(\xi_0 = i)]$
    is the cross entropy between $\pi_{0}^{\xi_0}(\thk)$ and $\alpha_0(\xi_0)$.
    This cross entropy $H_0$ is minimized when $\alpha_0(\xi_0) = \pi_{0}^{\xi_0}(\thk) $.
    Consequently, we initialize $\alpha_0(\xi_0)$ with $\pi_{0}^{\xi_0}(\theta^{k}) $ at subsequent iteration $k + 1$ for all $k \geq 0$,
    thus consistently applying the minimizer from the previous iteration.
\end{remark}
\begin{algorithm}
	\caption{Construction of $\surrogate$}
	\label{alg:construct-Q}
	\begin{algorithmic}[1]
		\Require $\thk$, $\alpha_0(\xi_0)$, $\zeta_T = 1$
		\For{$t = 1,\dots,T$}
            \State Update $\alpha_t(\xi_t\!)$ and $\zeta_{T{-}t}(\xi_{T{-}t{+}1}\!)$ using \eqref{eq: forward_filtering}, \eqref{eq: backward_filtering}
		\EndFor
        \State Compute $\pi_{t+1}^i$ and $\pi_{t}^{i, j}$ using \eqref{eq: forward_backward}, \eqref{eq: forward_backward_joint}
        \State Compute $\Qi[1][k]{}$ and $\Qi[2][k]{}$ using \eqref{eq:switch-surrogate-explicit}, \eqref{eq:subsys-surrogate-explicit}
	\end{algorithmic}
\end{algorithm}

\section{Numerical Experiments}\label{sec: numerical_experiment}
We evaluate the proposed model \eqref{eq: discrete_latent} and algorithm \EMpp (cf. \cref{alg:em-conceptual}) in various aspects.
First, \cref{sec: solver_comparison}
demonstrates on a synthetic example that \EMpp is more effective in identifying the switching system parameter than the popular alternatives \texttt{BFGS} \cite[\S 6.1]{nocedal1999numerical} and \texttt{Adam} \cite{kingma2014adam}.
Then, \cref{sec: comparison_tailored_algorithms} assesses the prediction accuracy of \EMpp against tailored algorithms for different switching systems.
Moreover,
we explore a robust identification scenario involving outliers to illustrate the importance of flexible modelling.
Finally, \cref{sec: case_study} compares the prediction accuracy of the proposed model \eqref{eq: discrete_latent} against different black box models on a nonlinear benchmark.
For all experiments, 
we use Gaussian and Student's t-distributions (cf. \cref{tab: example_distr}) with regularization terms
\begin{align}
    \reg_1(\dTheta) = &\; \textstyle \sum_{i = 1}^d \frac{\gamma_1}{2} \norm{\Theta_i}^2_{\mathrm{F}}, \label{eq: reg_exp}\\
    \reg_2(\dBeta) = &\;\textstyle \frac{1}{2}\sum_{i = 1}^d \gamma_2 \left[
            \tr(\Lambda_i) - \ln\det(\Lambda_i) 
        \right] + \gamma_3 \norm{B_i}^2_{\Lambda_i^{-1}}. \nonumber
\end{align}
\rev[reg-closed-form]{}{The choice of $\reg_2$ enables a closed-form solution for \eqref{eq: prob_subsys_iden}.}
\rev[old_reviewer2point1]{Softmax functions are translation invariant, 
i.e., $\sigma(x) = \sigma(x - a\mathbf{1}_d)$ with $a \neq 0$.
Through a change of variable $\tilde{x}_i = x_i - x_d$ for all $i \in \N_{[1, d]}$,
we maintain the same output, reduce the dimensionality of the unknown variables,
and ensure that $\sigma(\tilde{x}) = \sigma(\bar{x})$ if and only if $\tilde{x} = \bar{x}$ with $\tilde{x}_d = \bar{x}_d = 0$.
Consequently,
the parameter $\Theta_i \in \re^{n_z \times (d-1)}$ 
and \eqref{eq: switch_model} takes the form $\dTheta \mapsto \sigma\big(\bsmat{z_t^\top \Theta_i  & 0}^\top\big)$.}{}%
At iteration $k=0$, $\alpha_0(\xi_0)$ in \cref{alg:construct-Q} is randomly initialized. 
\rev[reviewer2point3]{}{%
    If no prior knowledge about the mode number $d$ is available,
    we grid search $d$ over a separate validation set.%
}
All experiments run on an Intel Core i7-11700 @ 2.50 GHz machine in \texttt{Python} 3.9, and \texttt{MATLAB 2018b}.

\subsection{Synthetic example}\label{sec: synthetic_example}
\begin{figure}[!tb]
    \centering
    \revnew{2}{%
    \includegraphics[width=0.4\textwidth]{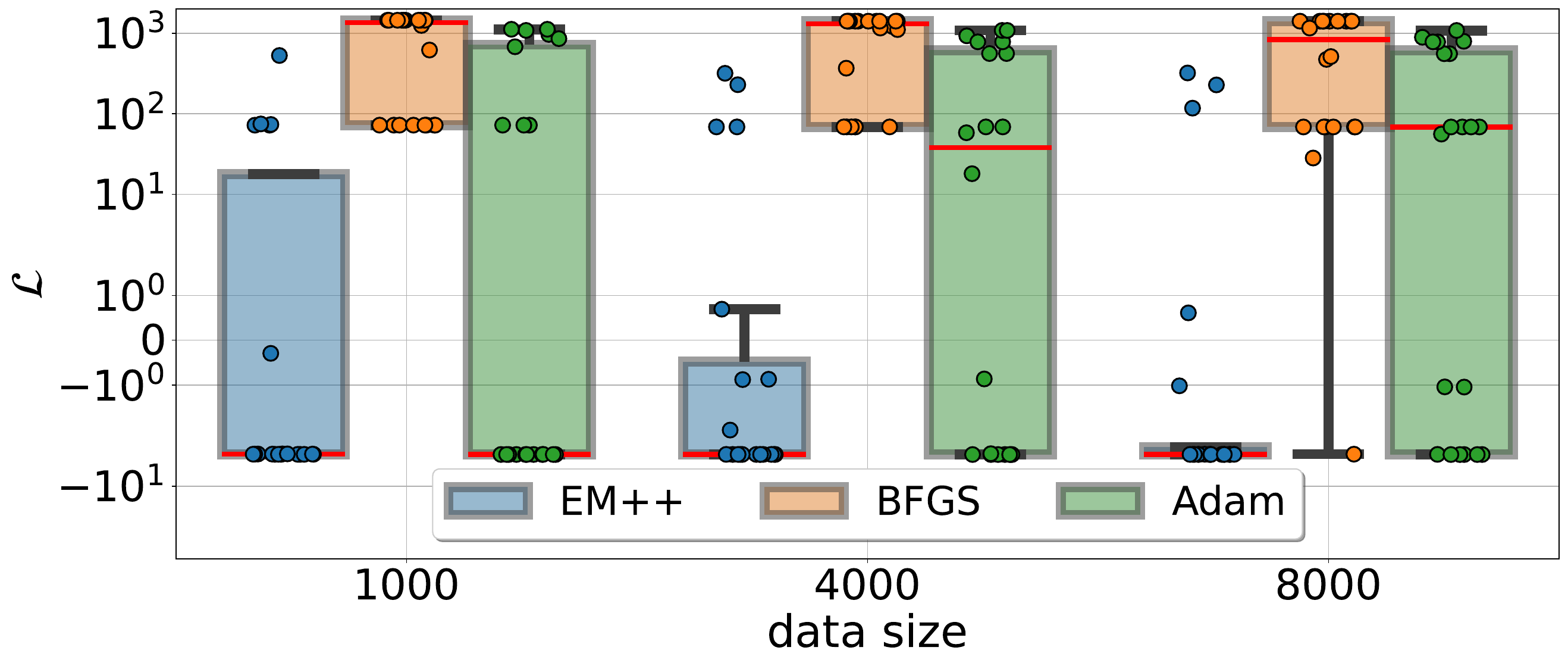}
    }{%
    \includegraphics[width=0.4\textwidth]{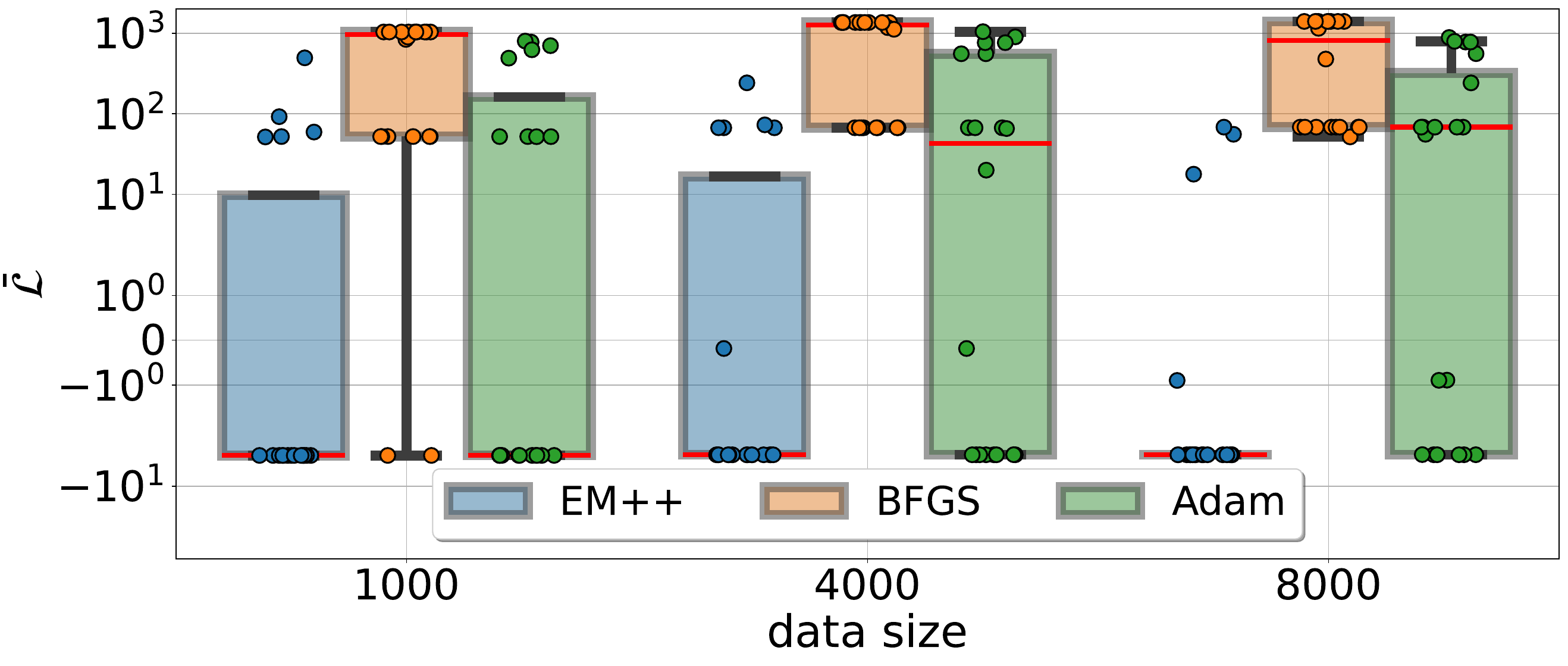}
    }
    \caption{Box plot and raw data of \revnew{2}{NLL $\loss$ \eqref{eq: loss}}{$\regLoss$ \eqref{eq: reg_loss}} using estimated parameters.
    The box draws from the first quartile $Q_1$ to the third quartile $Q_3$.  
    Whiskers extend up to $1.5$ times the interquartile range $\mathrm{IQR} = Q_3 - Q_1$ from $Q_1$ and $Q_3$.
    The median of each category is indicated by a red line.
    The individual value of \revnew{2}{NLL $\loss$}{$\regLoss$} for each method are shown as dots.}\label{fig: nll_boxplot}
\end{figure}

We compare \EMpp with \texttt{BFGS} \cite[\S 6.1]{nocedal1999numerical} and \texttt{Adam} \cite{kingma2014adam}
for identifying a switching system with 
\begin{equation*}
    \begin{aligned}
        A_1 = &\; \bsmat{
            0.9912& 0.1307 & 0.2 \\ 
            -0.1305&  0.9914 & 0.06},
        A_2 = \bsmat{
            0.94& 0.15 & -0.01\\ 
            -0.15&  0.94 & -0.13},\\
        A_3 = &\; \bsmat{
            0.97 & 0.4 & 0.1 \\ 
            -0.4&  0.97 & 0.1}, \\
        \Theta_1 = &\; \bsmat{ 30 &    10   \\  1 &   -16.07 \\-10 &    10},
        \Theta_2 = \bsmat{30  & 30 \\ 20 & -10 \\  0 &  0},
        \Theta_3 = \bsmat{24.8  &  0   \\ 11.38 & -28.62 \\-57.73 &   7.07}.
    \end{aligned}
\end{equation*}
The three subsystems are all of the form $y_{t+1} = A_i \bsmat{y_t \\ 1} + w_t$,
where the noise $w_t \sim \gauss(0, \Sigma_i)$ with
$\Sigma_i = 10^{-3}I$ for $i = 1, 2, 3$.
The covariance matrix is known a-priori as
\texttt{BFGS} and \texttt{Adam} cannot enforce a positive definite matrix constraint.
The implementations from \texttt{SciPy} \cite{2020SciPy-NMeth} and \texttt{PyTorch} \cite{NEURIPS2019_9015} are used with default hyperparameters.
For this experiment, we set
\revnew{2}{%
\begin{resetcolor}
\rev[typoIndex]{$\gamma_0 = \gamma_1 = \gamma_2 = 0$}{$\gamma_1 = \gamma_2 = \gamma_3 = 0$}
\end{resetcolor}
}{%
$\gamma_1=\gamma_3=10^{-10}$, $\gamma_2=0$
}
in \eqref{eq: reg_exp}\revnew{2}{.}{, as $\Lambda_i = \Sigma_i^{-1}$ are known.} 
As the problem \eqref{eq: reg_loss} is nonconvex,
we repeatedly solve it for 20 random initial guesses.
The same set of initial guesses is used for each solver to ensure a fair comparison.
All algorithms terminate when $\norm{\nabla \revnew{2}{\loss}{\regLoss}} \leq 10^{-\revnew{2}{4}{3}}$.
If a solver exceeds 30000 iterations,
the solution at the last iterate is used for comparison.
The runtime performance is detailed in \cref{tab: runtime_synthetic}.
The solution quality is evaluated
by the \revnew{2}{}{regularized} negative log-likelihood \revnew{2}{$\loss$ \eqref{eq: loss}}{$\regLoss$ \eqref{eq: reg_loss}} on a separate validation set with $T = 10000$.
The box plot of \revnew{2}{$\loss$}{$\regLoss$} is shown in \cref{fig: nll_boxplot}.
\begin{table}[tb]
        \caption{Mean value of runtime and number of iterations across 20 random initial guesses with different training data size}\label{tab: runtime_synthetic}
        \centering
    
\revnew{2}{
    \begin{resetcolor}
\rev[update-table]{
    \begin{table}[tb]
    \parbox{\linewidth}{Mean value of runtime and number of iterations across 20 random initial guesses with different training data size}\label{tab: runtime_synthetic_old}
    \centering
    \begin{tabular}{c c c r}
        \toprule 
        data size & method & \rev{runtime / iter. (\si{\second})}{} & num. of iter. \\
        \midrule 
        \multirow{3}{*}{$1000$} & \EMpp              &    \rev{0.1128}{} & 34 \\
                                & \texttt{BFGS}      &    \rev{0.0159}{} & 160\\
                                & \texttt{Adam}      &    \rev{0.0093}{} & 27697\\[-2pt]
        \midrule
        \multirow{3}{*}{$4000$} & \EMpp              &    \rev{0.2235}{} & 22 \\
                                & \texttt{BFGS}      &    \rev{0.0574}{} & 178 \\
                                & \texttt{Adam}      &    \rev{0.0347}{} & 29491 \\[-2pt]
        \midrule
        \multirow{3}{*}{$8000$} & \EMpp              &    \rev{0.4102}{} & 17 \\
                                & \texttt{BFGS}      &    \rev{0.1306}{} & 229 \\
                                & \texttt{Adam}      &    \rev{0.0774}{} & 30000\\[-2pt]
        \bottomrule
    \end{tabular}
    \end{table}
}{
        \begin{tabular}{c c r r}
            \toprule 
            data size               & method             & \rev{}{runtime (\si{\second})} & num. of iter. \\
            \midrule 
            \multirow{3}{*}{$1000$} & \EMpp              &    \textcolor{newcolor}{3.8352}    & 34 \\
                                    & \texttt{BFGS}      &    \textcolor{newcolor}{2.5440}    & 160\\
                                    & \texttt{Adam}      &    \textcolor{newcolor}{257.5821}  & 27697\\[-2pt]
            \midrule
            \multirow{3}{*}{$4000$} & \EMpp              &    \textcolor{newcolor}{4.9170}    & 22   \\
                                    & \texttt{BFGS}      &    \textcolor{newcolor}{10.2172}   & 178   \\
                                    & \texttt{Adam}      &    \textcolor{newcolor}{1023.3377} & 29491 \\[-2pt]
            \midrule
            \multirow{3}{*}{$8000$} & \EMpp              &    \textcolor{newcolor}{6.9734}    & 17   \\
                                    & \texttt{BFGS}      &    \textcolor{newcolor}{29.9074}   & 229   \\
                                    & \texttt{Adam}      &    \textcolor{newcolor}{2322.0000} & 30000\\[-2pt]
            \bottomrule
        \end{tabular}
}
    \end{resetcolor}
}{}
    \begin{tabular}{c c r r}
        \toprule 
        data size               & method             & \rev{}{runtime (\si{\second})} & num. of iter. \\
        \midrule 
        \multirow{3}{*}{$1000$} & \EMpp              &    \revnew{2}{}{2.3621}    & \revnew{2}{}{13} \\
                                & \texttt{BFGS}      &    \revnew{2}{}{2.4215}    & \revnew{2}{}{145}\\
                                & \texttt{Adam}      &    \revnew{2}{}{240.1643}  & \revnew{2}{}{25280}\\[-2pt]
        \midrule
        \multirow{3}{*}{$4000$} & \EMpp              &    \revnew{2}{}{4.3552}    & \revnew{2}{}{16}   \\
                                & \texttt{BFGS}      &    \revnew{2}{}{9.2570}   & \revnew{2}{}{155}   \\
                                & \texttt{Adam}      &    \revnew{2}{}{962.2217} & \revnew{2}{}{27029} \\[-2pt]
        \midrule
        \multirow{3}{*}{$8000$} & \EMpp              &    \revnew{2}{}{7.0467}    & \revnew{2}{}{15}   \\
                                & \texttt{BFGS}      &    \revnew{2}{}{19.4957}   & \revnew{2}{}{149}   \\
                                & \texttt{Adam}      &    \revnew{2}{}{2274.1030} & \revnew{2}{}{28110}\\[-2pt]
        \bottomrule
    \end{tabular}
\end{table}
\rev{
    Although individual iterations of \EMpp are slower,
    due to solving an internal optimization problem in $\texttt{Python}$,
    it requires significantly fewer iterations, %
    especially for larger training datasets.
    In contrast,
    both \texttt{BFGS} and \texttt{Adam} need more iterations as the dataset size increases.
    Despite fast individual iterations,
    the overall process of both algorithms becomes less efficient when given larger dataset.
}{
    Compared to \texttt{BFGS} and \texttt{Adam}, which requires more iterations as the dataset size increases,
    \EMpp requires significantly fewer iterations to reach the termination criterion.
    Consequently, \EMpp spends less runtime when the dataset size increases.
}
The effectiveness of the \EMpp method is further demonstrated by its convergence to points with lower loss,
as evidenced in \cref{fig: nll_boxplot}.
The \revnew{2}{NLL $\loss$}{regularized negative log-likelihood $\regLoss$} improves with larger training dataset\revnew{2}{s}{},
whereas both \texttt{BFGS} and \texttt{Adam} do not significantly benefit from the additional data.
\label{sec: solver_comparison}
\subsection{Comparison with tailored algorithms}\label{sec: comparison_tailored_algorithms}
We compare \EMpp against tailored algorithms
\cite{bemporad2018fitting,bemporad2022piecewise}
for different switching systems,
using 3 formulations for \eqref{eq: switch_model}: 
\begin{inlinelist}
    \item \texttt{only-mode}: $\dTheta \mapsto \sigma\big(\bsmat{\Theta_i & 0}^\top\big)$;
    \item \texttt{only-state}: $\dTheta \mapsto \sigma\big(\bsmat{z^\top\Theta & 0}^\top\big)$ where $\Theta_i = \Theta$; and
    \item \texttt{full-dependence}: $\dTheta \mapsto \sigma\big(\bsmat{z^\top\Theta_i & 0}^\top\big)$,
\end{inlinelist}
for all $i \in \Xi$.
To assess the flexibility of \eqref{eq: subsys_model} in a robust identification scenario,
we consider both the Gaussian distribution and the Student's t-distribution,
a common choice in robust system identification.
Unless otherwise specified,
each algorithm is trained with 5 initial guesses,
with hyperparameters tuned by hold-out validation.
The best solution on a separate validation data set is selected to compare prediction accuracy.
For trajectory prediction initialization,
both framework \cite{bemporad2018fitting} and our options \texttt{only-mode}, \texttt{full-dependence} use training and validation set.
The latter two options use \eqref{eq: forward_backward_summary} to compute initial mode distribution.
Denoting the ground truth \rev{with}{by} $y_n$ and the prediction \rev{with}{by} $\hat{y}_n$,
the prediction accuracy is measured by the $R^2$ score:
\vspace{-\baselineskip}
\vspace{-0.5\baselineskip}
{%
    \small%
\[
    R^2 = 1 - \frac{
        \sum_{n=1}^N \left(y_n - \hat{y}_n\right)^2
    }{
        \sum_{n=1}^N \left(y_n - \tfrac{1}{N} \sum_{n=1}^N y_n\right)^2
    }.
\]%
}%

\subsubsection{Switched Markov ARX system}\label{sec: switch_marx}
\begin{figure}[!bt]
    \centering
    \includegraphics[width=0.68\linewidth]{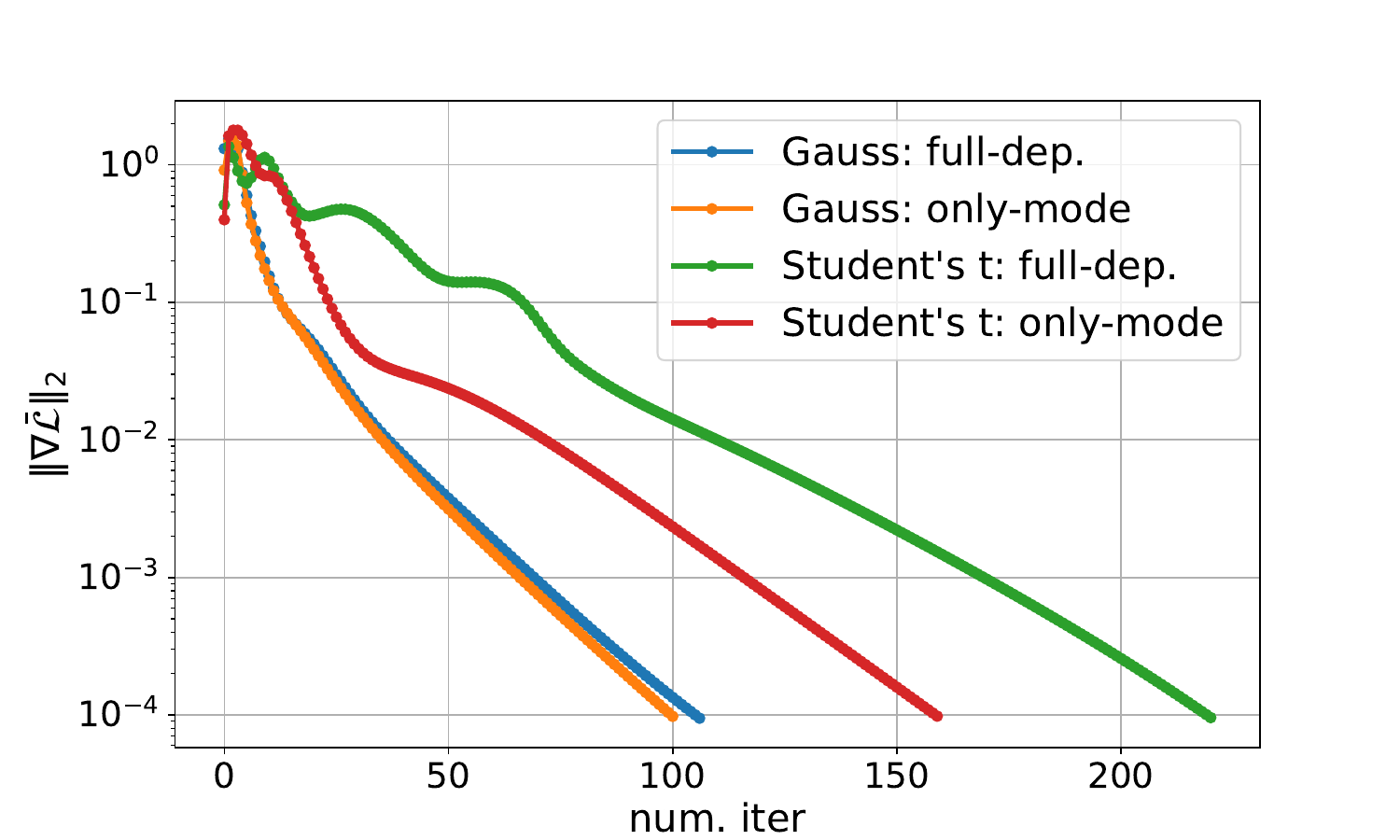}
    \caption{Gradient norm evaluation for Markov jump ARX system \eqref{eq: sys_MJARX} with different noise and switching options.}
    \label{fig: grad_norm_mjarx}
\end{figure}

\begin{table*}[tb]
    \centering
    \begin{threeparttable}
    \caption{$R^2$ score of recursive one-step-ahead prediction for switched Markov ARX system \eqref{eq: sys_MJARX} ($\uparrow$)}
    \label{tab: r2_markov_arx}
    \begin{tabular}{l  c  cc cc}
        \toprule \\
        [-10pt]& & \multicolumn{2}{c}{\EMpp (Gaussian distribution)} & \multicolumn{2}{c}{\EMpp (Student's t-distribution)} \\ [-2pt]
         & framework \cite{bemporad2018fitting} & \texttt{only-mode} & \texttt{full-dependence} & \texttt{only-mode} & \texttt{full-dependence} \\ [-2pt]
        \midrule \\ [-10pt]
        $p = 0$     & 0.9286  & 0.9526 & 0.9540 & 0.9536 & $\boldsymbol{0.9548}$ \\[-2pt]
        $p = 0.01$  & 0.9452  & 0.9341 & 0.9357 & 0.9528 & $\boldsymbol{0.9560}$ \\ [-2pt]
        $p = 0.05$  & 0.9367  & 0.8521 & 0.7623 & 0.8810 & $\boldsymbol{0.9477}$\\ [-2pt]
        \bottomrule
    \end{tabular}
\end{threeparttable}
\end{table*}
We collect a trajectory with $T = 10000$ data points from a switching Markov ARX system:
\vspace{-0.5\baselineskip}
\begin{equation}\label{eq: sys_MJARX}
    y_t = \beta_{\xi_t}^\top z_{t-1} + w_t,
\end{equation}
where $z_{t-1} = \begin{bsmallmatrix}y_{t-1} &y_{t-2} & u_{t-1} & u_{t-2} \end{bsmallmatrix}^\top$.
It consists of three subsystems with parameters 
$\beta_1^\top = \begin{bsmallmatrix} 1.143  & -0.4346 &  0.0572 &  0.2415\end{bsmallmatrix}$, 
$\beta_2^\top = \begin{bsmallmatrix} 0.9534 & -0.0475 &  0.0618 &  0.0336\end{bsmallmatrix}$, 
$\beta_3^\top = \begin{bsmallmatrix} 1.178  & -0.09   &  0.089  &  0.15  \end{bsmallmatrix}$,
and additive noise $w_t\sim\mathcal{N}(0, 0.025)$. 
Switching occurs according to a transition matrix
\(
    P = \begin{bsmallmatrix}
        0.25 & 0.1  & 0.65 \\
        0.55 & 0.35 &  0.1 \\
        0.15 & 0.15 &  0.7 \\
    \end{bsmallmatrix}.    
\)
The control input is uniformly randomly sampled in $[-1, 1]$.
The 10000 data points are split into a training (5000), validation (2500), and test (2500) set.
Each training data point has probability $p \in \{0\%, 1\%, 5\%\}$ of 
being perturbed by noise uniformly sampled from $[-\max_t \abs{y_t}, \max_t \abs{y_t}]$.
\rev[pmjarx]{}{The parameters $\gamma_1 = 10^{-4}$, $\gamma_2 = \gamma_3 = 10^{-8}$ in \eqref{eq: reg_exp} are grid-searched on the validation set.}
\EMpp is compared to the framework \cite{bemporad2018fitting}.
The solutions are evaluated using the recursive one-step ahead prediction:
at each time step $t$,
the algorithm predicts one-step ahead,
then is updated by the true observation $y_t$ \cite[algorithm 3]{bemporad2018fitting}.
For our stochastic model,
we compute the mean over 20 trajectory samples for evaluation.
The result is in \cref{tab: r2_markov_arx}.
As $f$ for the Student's t-distribution is a nonlinear function (cf. \cref{tab: example_distr}),
constructing the surrogate function \eqref{eq:subsys-surrogate} requires evaluating the linearization of $f$.
In contrast, the Gaussian distribution does not necessitate this step as its $f$ is linear (cf. \cref{tab: example_distr}).
To evaluate the impact of this additional majorization step via linearization,
we \rev{in addition}{additionally} compare the performance of the algorithm with those two noise options.
The progress plot for both options with the same initial guess is illustrated in \cref{fig: grad_norm_mjarx}.

\EMpp with the Student's t-distribution model demonstrates superior performance compared to the framework \cite{bemporad2018fitting},
despite requiring more iterations.
The increased iteration count owes to the necessity of computing a linearization \eqref{eq:subsys-surrogate} for the Student's t distribution,
which results in a looser approximation of the original loss.
However, this trade-off is justified by its robustness to outliers.
While the model with the Gaussian distribution achieves high $R^2$ score without outliers,
the scores drop significantly as the proportion of outliers increases,
since the Gaussian parameter estimation is sensitive to outliers.
This comparison highlights the necessity of a subsystem model tailored to the specific problem.
Comparing the \texttt{full-dependence} and \texttt{only-mode} options,
we observe that the former requirs more iterations due to its increased number of variables.
Nevertheless, the model with option \texttt{full-dependence}
reaches a similar score as the model utilizing the prior knowledge on the switching (\texttt{only-mode}) in most cases.
This suggests that the proposed method is able to learn the system's intrinsic structure without prior knowledge.

\subsubsection{Piecewise affine system}\label{sec: pwd}
\begin{figure}[!bt]
    \centering
    \includegraphics[width=0.68\linewidth]{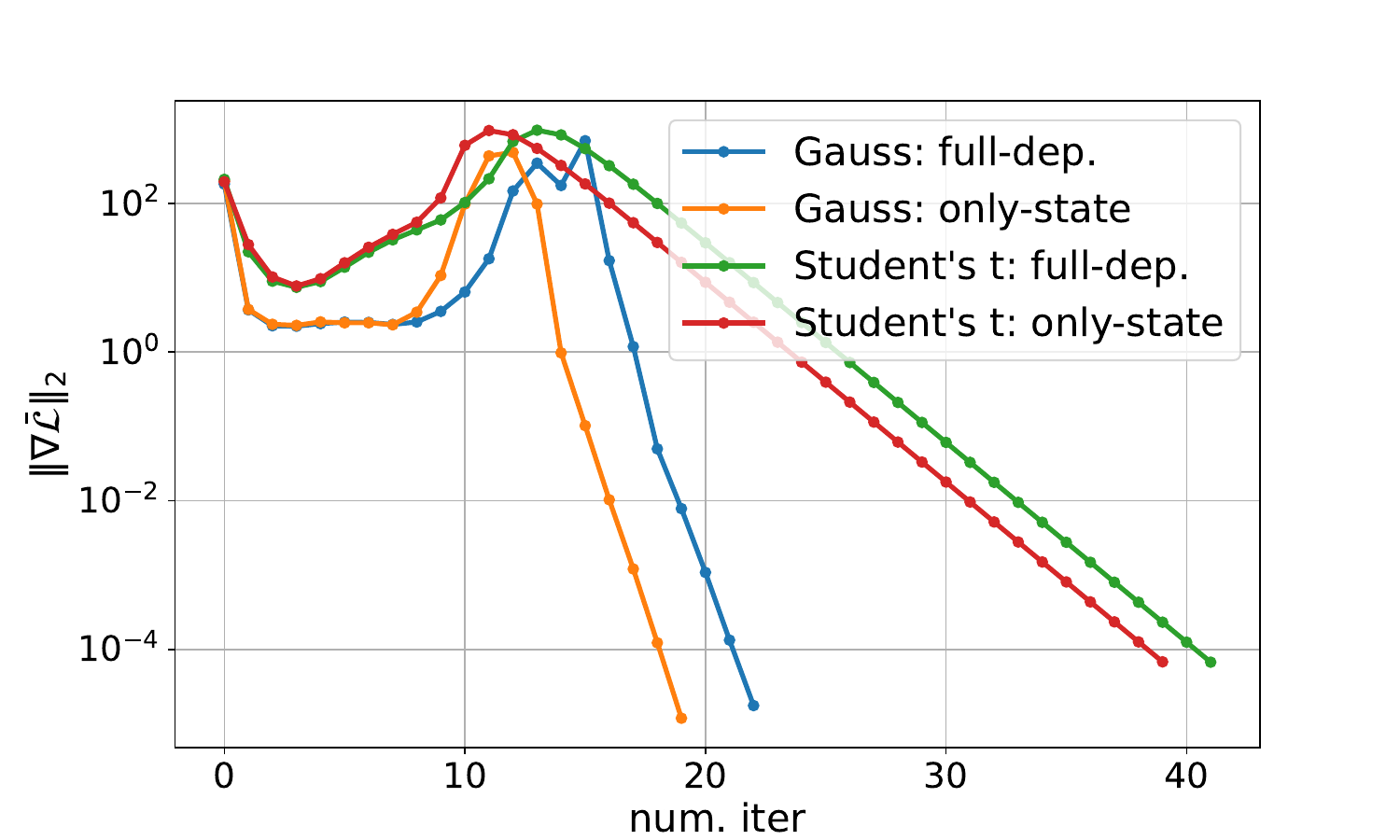}
    \caption{Gradient norm evaluation for piecewise affine system \eqref{eq: sys_PWA} with different noise and switching options.}
    \label{fig: grad_norm_pwa}
\end{figure}

We collect a trajectory with $T=10000$ data points from a piecewise affine system:
\vspace{-0.5\baselineskip}
\begin{equation}\label{eq: sys_PWA}
    \begin{cases}
        y_{t} = \beta_{1}^\top z_{t-1} + w_t & \text{if} \quad \beta_{0}^\top z_{t-1} \geq 0, \\
        y_{t} = \beta_{2}^\top z_{t-1} + w_t & \text{otherwise}, 
    \end{cases}%
\end{equation}%
where $z_{t-1} = \begin{bsmallmatrix}y_{t-1} & y_{t-2} & u_{t-1} & u_{t-2} & 1\end{bsmallmatrix}^\top$.
The system parameters are 
$\beta_0^\top = \begin{bsmallmatrix} 0.5 & 1   & 2    &-0.3 & 0.2 \end{bsmallmatrix}$, 
$\beta_1^\top = \begin{bsmallmatrix} 0.1 & 0.5 & -0.4 & 0.3 & 0   \end{bsmallmatrix}$,
$\beta_2^\top = \begin{bsmallmatrix} 0.2 & 0.4 &  0.1 & 0.4 & 0   \end{bsmallmatrix}$,
and the system is subject to an additive noise $w_t\sim\mathcal{N}(0, 10^{-4})$.
The control input is randomly sampled from $\mathcal{N}(0, 0.25)$.
The data are split as in \cref{sec: switch_marx}.
\rev[ppwa]{}{The parameter $\gamma_1 = 10^{-6}$, $\gamma_2 = \gamma_3 = 10^{-8}$ in \eqref{eq: reg_exp} are grid-searched on the validation set.}
\EMpp is compared with 
the framework \cite{bemporad2018fitting}, 
and a tailored algorithm for piecewise affine models, PARC \cite{bemporad2022piecewise}.
The identified models are evaluated using the open-loop prediction with the same control input that generates the test dataset.
For our stochastic model \eqref{eq: discrete_latent},
we employ $1\%$ trimmed mean of 500 sampled trajectories as the prediction,
since the trimmed mean is less sensitive to rare events.
The $R^2$ score is listed in \cref{tab: r2_pwa}.
Same as \rev{}{in} \cref{sec: switch_marx}, 
we compare the performance of the algorithm using both Gaussian and Student's t-distribution.
This comparison evaluates the impact of the additional majorization step described in \cref{prop:separable}.
The progress plot for both options with the same initial guess is in \cref{fig: grad_norm_pwa}.
\begin{table*}[tb]
    \centering
    \begin{threeparttable}
    \caption{$R^2$ score of open-loop prediction for piecewise affine system \eqref{eq: sys_PWA} ($\uparrow$)}
    \label{tab: r2_pwa}
    \begin{tabular}{l  cc  cc cc}
        \toprule \\
        [-10pt]& & & \multicolumn{2}{c}{\EMpp (Gaussian distribution)} & \multicolumn{2}{c}{\EMpp (Student's t-distribution)} \\
        & framework \cite{bemporad2018fitting} & PARC \cite{bemporad2022piecewise}& \texttt{only-state} & \texttt{full-dependence} & \texttt{only-state} & \texttt{full-dependence} \\
        \midrule \\ [-10pt]
        $p = 0$      & 0.8569 & 0.9883 & 0.9905 & $\boldsymbol{0.9909}$ & 0.9904                & 0.9905 \\
        $p = 0.01$   & 0.8601 & 0.9764 & 0.9537 & 0.9572                & 0.9901                & $\boldsymbol{0.9903}$ \\ 
        $p = 0.05$   & 0.8428 & 0.9050 & 0.8600 & 0.8689                & $\boldsymbol{0.9881}$ & $\boldsymbol{0.9881}$\\ 
    \bottomrule
    \end{tabular}
\end{threeparttable}
\end{table*}

Both PARC and \EMpp outperform the framework \cite{bemporad2018fitting}.
The enhanced performance can be attributed to 
the robustness against subsystem estimation error achieved by softmax modelling in \eqref{eq: switch_model},
compared to a Voronoi-distance-based method used by the framework \cite{bemporad2018fitting}.
As the number of outliers increases,
the performance of both PARC and of the proposed models with a Gaussian distribution decreases significantly.
In contrast, the models with the Student's t-distribution maintain high prediction accuracy,
despite a slight increase in the number of iterations due to the additional majorization step via linearization \eqref{eq:subsys-surrogate}.
We observe that option \texttt{full-dependence} requires more iterations than option \texttt{only-state} due to its increased number of variables.
Nevertheless, it achieves a similar $R^2$ score as the option \texttt{only-state},
akin to the previous example.
This suggests that our proposed method is able to learn the system's intrinsic structure without prior information.
This ability is particularly beneficial when (accurate) prior information is unavailable.

\subsubsection{Cart system}\label{sec: cart_system}
We further compare \EMpp with PARC \cite{bemporad2022piecewise} on a cart system \cite{bemporad2022piecewise} that is \rev{piecewise nonlinearity}{piecewise nonlinear}%
\rev[cart-description]{:
\vspace{-\baselineskip}
\vspace{-\baselineskip}
{%
\small%
\begin{align}
    &  
    \left\{
    \begin{aligned}
    \dot{p} &= \tfrac{1}{M+M_1} \left( F - \varphi_b(\dot{p}, b_1 + b) - \varphi_a(p-a-\alpha_1,k_1) \right) \\ \nonumber
    \dot{T} &= \tfrac{1}{\Theta} \left( \tfrac{T_0 - T}{R} + \tfrac{T_1-T}{R_1} \right), \quad \text{if } p-a \leq \alpha_1 \nonumber
    \end{aligned}
    \right. \\
    & \left\{
    \begin{aligned}
    \dot{p} &= \tfrac{1}{M+M_2} \left( F - \varphi_b(\dot{p}, b_2 + b) - \varphi_a(p+a-\alpha_2,k_2) \right) \\ \nonumber
    \dot{T} &= \tfrac{1}{\Theta} \left( \tfrac{T_0 - T}{R} + \tfrac{T_2-T}{R_2} \right), \quad \text{if } p+a \geq \alpha_2 \nonumber
    \end{aligned}
    \right. \\
    &
    \left\{
    \begin{aligned}
    \dot{p} &= \tfrac{1}{M} \left( F - \varphi_b(\dot{p}, b) \right) \\
    \dot{T} &= \tfrac{T_0 - T}{\Theta R} , \quad \text{otherwise} \label{eq: cart_system}
    \end{aligned}
    \right.
\end{align}    
}%
where
\(
        \varphi_a(\Delta p, k) = k \Delta p + \frac{k}{5} \Delta p^3,
        \varphi_b(\dot{p}, b) = b\dot{p} + \frac{b}{5} |\dot{p}| \dot{p}.
\)
The cart is controlled by a switching command $u \in \{-1, 0, 1\}$ that is generated every $T_s = \SI{0.5}{\second}$ according to a transition matrix:
\(
    P = \begin{bsmallmatrix}
        \frac{29}{30} & \frac{1}{60} & \frac{1}{60} \\
        \frac{1}{60}  & \frac{29}{30} & \frac{1}{60} \\
        \frac{1}{60}  & \frac{1}{60} &  \frac{29}{30} \\
    \end{bsmallmatrix}.
\)
The control input determines $F = u$. 
The other parameters can be found in \cite[TABLE \uppercase\expandafter{\romannumeral 8\relax}]{bemporad2022piecewise}.
}{%
.
    The cart moves longitudinally between two bumpers with states $p$ (position), $\dot{p}$ (velocity), and $T$ (temperature).
    A detailed system description and hyperparameter values can be found in \cite[\S VI]{bemporad2022piecewise}.
}
We collect $4000$ training data points and $550$ validation data points initialized at 
$x_0 = \begin{bsmallmatrix}
    p & \dot{p} & T     
    \end{bsmallmatrix} = \begin{bsmallmatrix}
        2 & 0 & 25 + 273.15
    \end{bsmallmatrix}$,
and $530$ testing data points initialized at
$\begin{bsmallmatrix}
    p & \dot{p} & T     
\end{bsmallmatrix} = \begin{bsmallmatrix}
    1.5 & 0 & 40 + 273.15
\end{bsmallmatrix}$.
We use the state-input pair $z_t = \bsmat{x_t^\top & u_t & 1}^\top$ as input and produce $x_{t+1}$ as output.
\rev{To}{Same as in \cite{bemporad2022piecewise}, to} ensure that all state variables are of the same order of magnitude, 
we convert the temperature $T$ into Celsius and divide it by $10$.
PARC is initialized with \texttt{K-means++} and with \rev{}{the} same hyperparameters as in \cite{bemporad2022piecewise},
while \EMpp uses random initialization.
\rev[pcart]{}{The parameters $\gamma_1 = \gamma_2 = \gamma_3 = 10^{-5}$ in \eqref{eq: reg_exp} are grid-searched on the validation set.}
The identified models are evaluated using the open-loop prediction on the last $500$ data points of the testing dataset,
using the same control input that generates the test data.
For \rev{}{the} option \texttt{full-dependence}, 
we estimate the initial mode distribution using the first $30$ points of testing dataset via \eqref{eq: forward_backward_summary}.
For our stochastic model \eqref{eq: discrete_latent},
we employ $1\%$ trimmed mean of 500 sampled trajectories as the prediction,
since the trimmed mean is less sensitive to rare events.
The $R^2$ score for each method is summarized in \cref{tab: r2_cart}.

\begin{table}[!htbp]
    \centering
    \caption{$R^2$ of open-loop prediction for cart system \rev{\eqref{eq: cart_system}}{} ($\uparrow$)}
    \label{tab: r2_cart}
    \begin{tabular}{c l ccc}
         \toprule \\[-10pt]
        $d$ & method name                   & $R^2(p)$ & $R^2(\dot{p})$ & $R^2(T)$ \\
        \midrule \\ [-10pt]
        \multirow{5}{*}{3} 
        & PARC \cite{bemporad2022piecewise} & ${0.927}$                 & ${0.900}$                & $\boldsymbol{0.920}$\\ 
        & \EMpp (\texttt{only-state})               & $\boldsymbol{0.931}$      & $\boldsymbol{0.909}$     & 0.912\\
        & \EMpp (\texttt{full-dependence})                     & 0.917                     & 0.901                    & 0.877\\
        \midrule
        \multirow{5}{*}{5} 
        & PARC \cite{bemporad2022piecewise} & $0.948$               & ${0.930}$             & $0.935$\\ 
        & \EMpp (\texttt{only-state})               & 0.948                 & $\boldsymbol{0.935}$  &  0.928\\
        & \EMpp (\texttt{full-dependence})                     & $\boldsymbol{0.972}$  & 0.931                 &  $\boldsymbol{0.976}$ \\
        \midrule
        \multirow{5}{*}{7} 
        & PARC \cite{bemporad2022piecewise} & $\boldsymbol{0.982}$ & $0.967$                & $\boldsymbol{0.978}$ \\ 
        & \EMpp (\texttt{only-state})               & 0.959                & 0.947                  & 0.946\\
        & \EMpp (\texttt{full-dependence})                     & 0.978                & $\boldsymbol{0.971}$   & 0.966 \\
        \bottomrule
    \end{tabular}
\end{table}
From \cref{tab: r2_cart}, we observe that \EMpp achieves comparable $R^2$ to the tailored algorithm \cite{bemporad2022piecewise}.
When the number of modes $d$ matches the underlying true model,
the case \texttt{only-state} exhibits higher prediction accuracy,
as the model aligns with the underlying switching mechanism.
However, by increasing the number of modes,
the case \texttt{full-dependence} has more flexibility to capture the underlying nonlinearity in the cart system \rev{\eqref{eq: cart_system}}{},
resulting in higher $R^2$ score.

\subsection{Case study: Coupled electric drives}\label{sec: case_study}
\begin{figure}[!tb]
    \centering
    \revnew{2}{%
        \includegraphics*[width=0.8\linewidth]{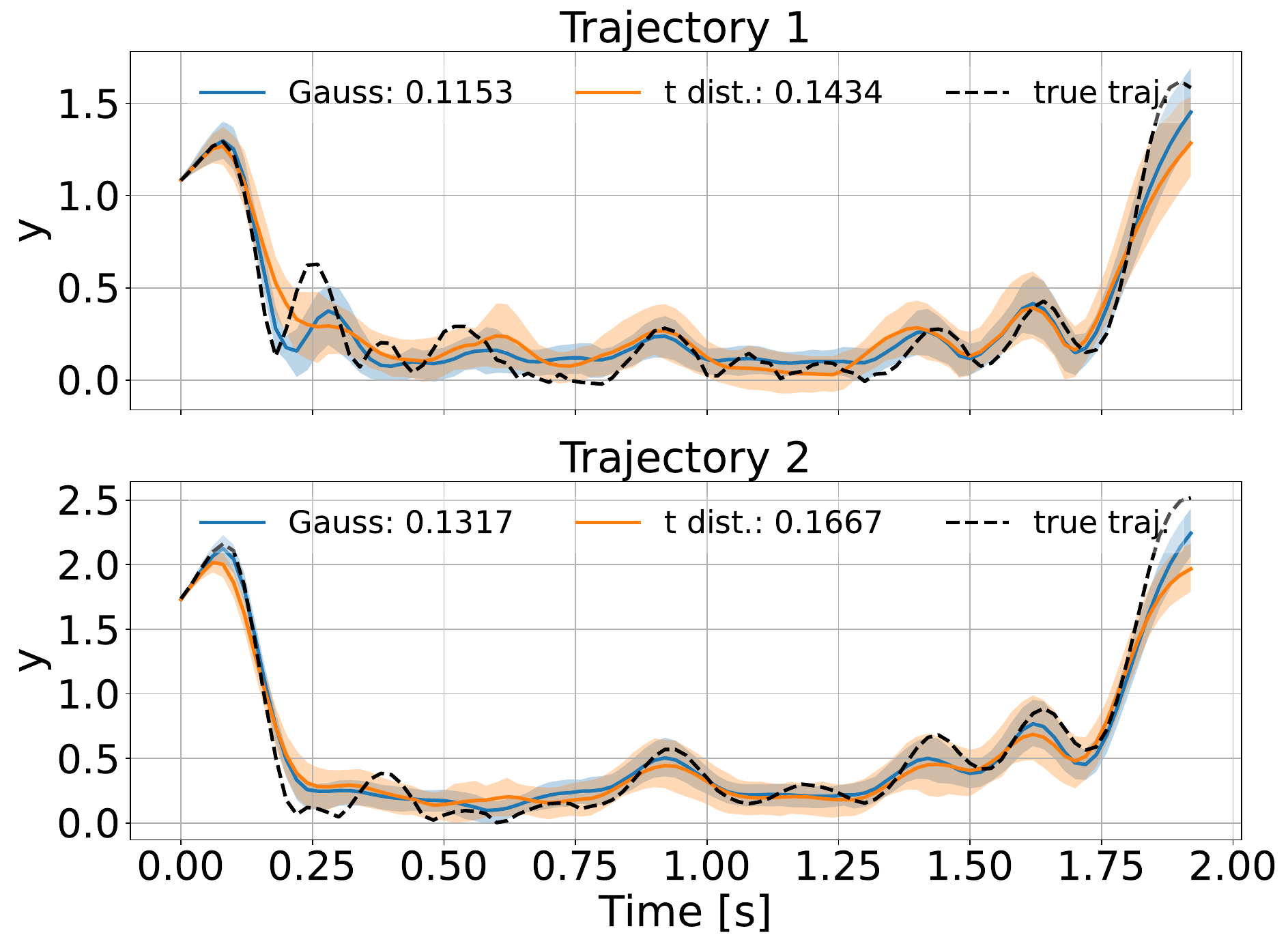}
    }{%
        \includegraphics*[width=0.8\linewidth]{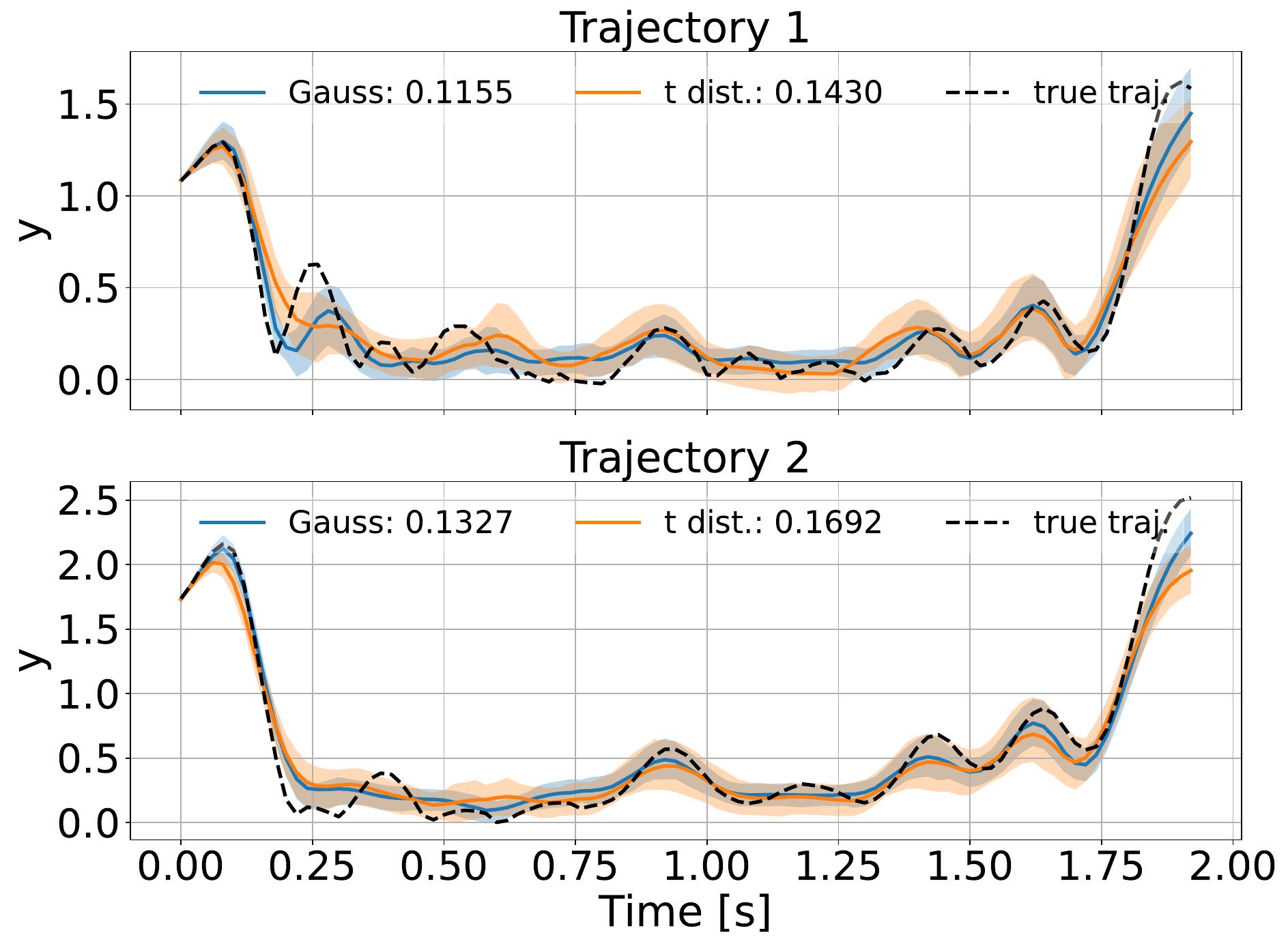}
    }
    \caption{Open-loop prediction for coupled electric drives. Colored area shows quantiles between $0.25$ and $0.75$.}
    \label{fig: open-loop_prediction}
\end{figure}

This section examines the expressiveness of our proposed method through a coupled electric drives system dataset \cite{wigren2017coupled}.
The system employs two electric motors to drive a pulley via a flexible belt.
The pulley is supported by a spring, resulting in lightly damped dynamics.
The system measures the absolute angular velocity of the pulley using a direction-insensitive sensor.
The dataset consists of two trajectories, each containing 500 data points sampled \rev{in}{at} 50 Hz.
The data points in each trajectory are allocated into training (300), validation (100), and testing (100) sets.
We choose the mode number $d = 8$ and use ARX form with 
$z_{t-1} = \begin{bsmallmatrix}y_{t-1}& \dots & y_{t-3}& u_{t-1}& \dots & u_{t-3} & 1\end{bsmallmatrix}^\top$.
This switching ARX model is compared against different types of ARX models listed in \cref{tab: coupled_drives}.
In particular, the neural network ARX (NN-ARX) consists of two hidden layers,
each with $32$ neurons. 
We use one layer of LSTM with hidden size equals to $32$.
The number of modes in PARC is $8$, same as for our model.
All models use the same lag value with hyperparameters grid-searched on the validation set.
Each model is trained for 20 times with different random initial guesses.
The best parameters on the validation set are chosen for comparison,
which evaluates the open-loop simulation error of each identified model.
\rev[pdrive]{}{The parameters $\gamma_1 = \gamma_2 = 10^{-6}$, \revnew{2}{$\gamma_3 = 0$}{$\gamma_3 = 10^{-10}$} in \eqref{eq: reg_exp} are grid-searched on the validation set}.
For simulation initialization, 
both LSTM and model \eqref{eq: discrete_latent} use the training and validation set.
The latter uses \eqref{eq: forward_backward_summary} to compute the initial mode distribution.
Same as \rev{}{in} \cite{beintema2023continuoustime},
we use root-mean-square error (RMSE) as comparison metric:
\(
    \mathrm{RMSE} = \sqrt{\frac{1}{N}\sum_{n=1}^N \left(y_n - \hat{y}_n\right)^2}
\)
where $y_n$ is the ground truth trajectory and $\hat{y}_n$ denotes the open-loop prediction.
As both our model \eqref{eq: discrete_latent} and Gaussian Process (GP)-ARX are stochastic models that output\rev{s}{} a trajectory distribution,
we sample 500 trajectories from each identified distribution. 
We compute the $1\%$ trimmed mean of the sampled trajectories as the prediction for the $\mathrm{RMSE}$,
which is less sensitive to rare events.
The open-loop prediction of identified proposed model in shown in \cref{fig: open-loop_prediction}.
The RMSE of all methods is summarized in \cref{tab: coupled_drives}.
\begin{table}[!tb]
    \centering
    \caption{RMSE of open-loop prediction ($\downarrow$)
    (The best three results are in bold, bold and italic, and italic font, respectively.)
    }\label{tab: coupled_drives}
    \begin{tabular}{l ll}
        \toprule\\ [-10pt]
        method & \multicolumn{2}{c}{RMSE} \\[-4pt]
        & traj. 1 & traj. 2 \\
        \midrule \\ [-10pt]
            ARX                              & 0.196 & 0.344 \\
            GP-ARX (with \texttt{rbf})       & 0.143 & 0.207 \\
            GP-ARX (with \texttt{Matern32})  & 0.149 & 0.204  \\
            NN-ARX (with \texttt{ReLU})      & 0.212 & 0.214 \\
            NN-ARX (with \texttt{SiLU})      & 0.161 & 0.179 \\
            LSTM                             & \textit{\textbf{0.104}} & \textit{\textbf{0.121}} \\
        \midrule \\ [-10pt]
            PARC \cite{bemporad2022piecewise}& 0.178 & 0.216 \\
            DT subspace encoder\cite{beintema2023continuoustime}                               & 0.130 & 0.145 \\ 
            CT subspace encoder ($\Delta t = 0.03$)\cite{beintema2023continuoustime}            & \textbf{0.085} & \textbf{0.072} \\ 
        \midrule \\ [-10pt]
        \EMpp (Gaussian distribution)                                                    &  \revnew{2}{\textit{0.115}}{\textit{0.116}} & \revnew{2}{\textit{0.132}}{\textit{0.133}}\\
        \EMpp (Student's t-distribution)                                                         &  \revnew{2}{0.143}{0.143} & \revnew{2}{0.167}{0.169} \\
        \bottomrule 
    \end{tabular}
\end{table}

As shown in \cref{tab: coupled_drives}, 
our method achieves better prediction accuracy compared to GP-ARX, NN-ARX, and PARC.
This superior performance demonstrates the expressiveness of the proposed model.
The enhanced result of CT subspace encoder may attribute to its continuous-time modelling approach.
However, our method is comparable to its discrete-time variant (DT subspace encoder) and LSTM,
with the added benefits of a simpler structure and improved interpretability.
This maintains the model's expressive power while ensuring an easier understanding and analysis.

\section{Conclusion}\label{sec: conclusion}
This work presented a general switching system model that encompasses various popular models.
Additionally, we proposed an algorithm \EMpp (cf. \cref{alg:em-conceptual}) to identify the parameters of this generic model by solving a regularized maximum likelihood estimation problem.
We proved that \EMpp converges to stationary points under suitable assumptions.
Finally, the effectiveness of both the presented model and the algorithm was demonstrated in a series of numerical experiments.
Future work includes: 
generalizing the model with continuous latent variables;
adapting \EMpp for online use;
and integrating the model into a model-based control framework such as model predictive control.

\appendix
\section{\revnew{2}{}{Auxiliary lemmas}}\label{app: auxiliary_lemmas}

\revnew{2}{}{
\begin{lemma}\label{lem: local_strong_convexity_on_a_compact_set}
    Consider $\mathscr{f}: \mathcal{X}_{\mathscr{f}} \to \re$ of the form
    \[
        \mathscr{f}(a) = \textstyle\sum_{t=0}^{T-1}\sum_{i=1}^{n_g} g_i\big(T_i(a, y_t)\big),
    \]
    where $y_{t} \in \re^{n}$ for all $t\in\N_{[0, T-1]}$,
    $T_i: \mathcal{X}_{\mathscr{f}} \times \re^n \to \mathcal{X}_g \subseteq \re^{m}$ is an affine transformation for $a$ for all $i\in\N_{[1, n_g]}$,
    and each function $g_i: \mathcal{X}_g \to \re$ is convex and differentiable.
    If there exists a nonempty index set $\mathcal{S} \subseteq \N_{[1, n_g]}$,
    such that 
    \begin{itemize}[topsep=-8pt]
        \item $g_i$ is of class $\mathcal{C}^2$ with $\nabla^2 g_i(x) \succ 0$ for all $i \in \mathcal{S}$ and all $x \in \mathcal{X}_g$; and
        \item $\sum_{t=0}^{T-1}\sum_{i\in\mathcal{S}} \jac{\!}^\top_{it}\jac{\!}_{it} \succ 0$,
        where $\jac{\!}_{it}$ is the Jacobian of $T_i(a, y_t)$ w.r.t. $a$,
    \end{itemize}
    then for any compact set $\Omega \subseteq \mathcal{X}_{\mathscr{f}}$, there exists an $\varepsilon > 0$, such that
    for all $a, \tilde{a}\in\Omega$,
    \[
        \inner{\nabla \mathscr{f}(a) - \nabla \mathscr{f}(\tilde{a})}{a - \tilde{a}} \geq \varepsilon \norm{a - \tilde{a}}^2.
    \]
\end{lemma}
\begin{pf}
    Let
    \[
        \tilde{\mathscr{f}}(a) \dfn \textstyle\sum_{t=0}^{T-1}\sum_{i\in\mathcal{S}} g_i\big(T_i(a, y_t)\big).
    \]
    By the convexity of $g_i$ for all $i \in \N_{[1, n_g]}$,
    we have that
    \begin{equation}\label{eq: gradient_diff_inner_lower_bound}
        \inner{\nabla \mathscr{f}(a) {-} \nabla \mathscr{f}(\tilde{a})}{a {-} \tilde{a}}
        {\geq} \inner{\nabla \tilde{\mathscr{f}}(a) {-} \nabla \tilde{\mathscr{f}}(\tilde{a})}{a {-} \tilde{a}}.
    \end{equation}
    By the chain rule, the Hessian of $\tilde{\mathscr{f}}$ is given by
    \begin{equation}\label{eq: hessian_tilde_f}
        \nabla^2 \tilde{\mathscr{f}}(a) 
        = \textstyle\sum_{t=0}^{T-1}\sum_{i\in\mathcal{S}} \jac{\!}_{it}^\top \nabla^2 g_i(T_i(a, y_t)) \jac{\!}_{it},
    \end{equation}
    As $\nabla^2 g_i$ is continuous and positive definite for all $i \in \mathcal{S}$ and all $x \in \mathcal{X}_g$,
    by the extreme value theorem,
    given a compact set $\Omega \subseteq \mathcal{X}_{\mathscr{f}}$,
    there exists a constant $\mu > 0$ such that
    \(\nabla^2 g_i(T_i(a, y_t)) \succeq \mu I_m\) for all $i \in \mathcal{S}$ and all $t \in \N_{[0, T-1]}$.
    Substituting this into \eqref{eq: hessian_tilde_f}, we have that
    \begin{equation}\label{eq: hessian_tilde_f_lower_bound_data}
        \nabla^2 \tilde{\mathscr{f}}(a)
        \succeq \mu \textstyle\sum_{t=0}^{T-1}\sum_{i\in\mathcal{S}} \jac{\!}_{it}^\top \jac{\!}_{it}
    \end{equation}
    Since $T_i$ is an affine transformation for $a$,
    the Jacobian $\jac{\!}_{it}$ only depends on the data $y_t$.
    Let $\lambda_{\min} > 0$ be the smallest eigenvalue of matrix $\textstyle\sum_{t=0}^{T-1}\sum_{i\in\mathcal{S}} \jac{\!}_{it}^\top \jac{\!}_{it}$.
    \eqref{eq: hessian_tilde_f_lower_bound_data} yields that
    \(
         \nabla^2 \tilde{\mathscr{f}}(a) \succeq \mu \lambda_{\min} I.
    \)
    Thus, $\tilde{\mathscr{f}}$ is locally strongly convex on $\Omega$.
    By \cite[proposition 12.60]{rockafellar2009variational},
    there exists an $\varepsilon = \mu \lambda_{\min} > 0$, such that
    \[
        \inner{\nabla \tilde{\mathscr{f}}(a) - \nabla \tilde{\mathscr{f}}(\tilde{a})}{a - \tilde{a}} \geq \varepsilon \norm{a - \tilde{a}}^2.
    \]
    Combining this with \eqref{eq: gradient_diff_inner_lower_bound} completes the proof.%
\qed\end{pf}
}

\revnew{2}{}{%
\begin{lemma}\label{lem: bounded_weights}
    Under \cref{asp: subsys_model},
    given a compact set $\Omega \subseteq \mathcal{X}_{\theta}$,
    there exist $\varepsilon_1, \varepsilon_2 > 0$, such that
    for all $\theta \in \Omega$,
    \begin{enumerate}[label=(\alph*), topsep=-8pt]
        \item $\Pi(\theta) \geq \varepsilon_1$ holds component-wise; \label{item: bounded_pi}
        \item $\nabla f[G_t(\beta_j)] \geq \varepsilon_2$ holds elementwise for all $j \in \Xi$ and $t\in \N_{[0, T-1]}$. \label{item: bounded_gradient}
    \end{enumerate}
\end{lemma}
}

\revnew{2}{}{%
\begin{pf}
    Denote $L_{\modeSeq}(\theta) = \ve{c}_{\modeSeq} - \Psi_{\modeSeq}(\dTheta) - \Phi_{\modeSeq}(\dBeta)$ for all $\modeSeq \in \Xi^{T+1}$.
    Since by \eqref{eq: loss} and \eqref{eq: expansion_nll} we have
    \(
        L_{\modeSeq}(\theta) = \ln \pjoint ,
    \)
    it follows that
    \(
        \exp(L_{\modeSeq}(\theta)) = \pjoint.
    \)
    Therefore, we have
    \begin{align*}
        \Pi_{\modeSeq}(\theta) 
        = &\; p(\modeSeq \mid \ve{y}, z_0\midsc \theta)
        = \frac{
                \pjoint 
            }{
                p(\ve{y}, z_0 \midsc \theta)
            } \\
        =&\; \frac{
                \pjoint 
            }{
                \sum_{\modeSeq' \in \Xi^{T+1}} p(\ve{y}, \modeSeq', z_0 \midsc \theta)
            }
        = 
        \sigma_{\modeSeq}(L(\theta))
    \end{align*}
Following from \cref{asp: subsys_model},
function $L_{\modeSeq}$ is continuous for all $\modeSeq \in \Xi^{T+1}$.
Thus, the function composition $\Pi_{\modeSeq}$ is continuous.
By the extreme value theorem and the definition of softmax function, we have
\[
    \Pi_{\modeSeq}(\theta) \geq \min_{\theta' \in \Omega, \modeSeq' \in \Xi^{T+1}} \Pi_{\modeSeq'}(\theta') \dfn \varepsilon_1 > 0
\]
for all $\theta \in \Omega$ and $\modeSeq \in \Xi^{T+1}$, which proves \crefpart{lem: bounded_weights}{item: bounded_pi}.

Similarly,
by \cref{asp: subsys_model},
$\beta \mapsto \nabla f(G_t(\beta))$ is continuous for all $t \in \N_{[0, T-1]}$.
By the extreme value theorem,
\[
    \partial_i f(G_t(\beta)) \geq \min_{\beta' \in \Omega_{\beta}, j\in\N_{[1, n_g]}} \partial_j f(G_t(\beta')) \dfn \varepsilon_{t} > 0
\]
for all $t \in \N_{[0, T-1]}$, $i\in\N_{[1, n_g]}$, and $\beta \in \Omega_{\beta}$,
where $\Omega_{\beta}$ is the projection of $\Omega$ onto the $\beta$-subspace,
and $\partial_i f$ denotes the partial derivative of $f$ w.r.t. the $i$-th component.
The positive lower bound follows from $\nabla f(x) > 0$ holding component-wise.
Let $\varepsilon_2 = \min_{t \in \N_{[0, T-1]}} \varepsilon_{t} > 0$,
which proves \crefpart{lem: bounded_weights}{item: bounded_gradient}.
\qed\end{pf}
}

\section{\rev{}{Deferred proofs}}\label{app: proofs}
\subsection{Proof of \cref{lem: general_compactness}}\label{prf: general_compactness}
By \cref{asp: regularizers}, \cref{eq:Lt1-def,eq:Lt2-def} we have
\begin{alignat*}{6}
    &\min_{\modeSeq \in \Xi^{T+1}}  \Psi_{\modeSeq} (\dTheta) \,&+&\, \reg_1(\dTheta) 
    &\,\geq&\, \kappa_1 \norm{\dTheta}  &+\, c_1, \\
    &\min_{\modeSeq \in \Xi^{T+1}} \Phi_{\modeSeq}(\dBeta)    \,&+&\, \reg_2(\dBeta)
    &\,\geq&\,  \kappa_2 \norm{\dBeta}  &+\, c_2.
\end{alignat*}%
Following from $\lse(x) \leq \max_i x_i + \ln n_x$,
\vspace{-\baselineskip}
\vspace{-\baselineskip}
{%
\small%
\begin{align*}
    \regLoss(\theta) = &\;-\lse\left( \ve{c} - \Psi(\dTheta) - \Phi(\dBeta) - \reg(\theta) \right)\\
                     \geq&\;  \min_{\modeSeq \in \Xi^{T+1}} -\ve{c}_{\modeSeq} + \Psi_{\modeSeq}(\dTheta) + \Phi_{\modeSeq}(\dBeta) + \reg(\theta)  - \ln d^{T+1},
\end{align*}
}%
where $\reg(\theta) = \reg_1(\dTheta) + \reg_2(\dBeta)$.
Thus,
there exist $\kappa_3 > 0$, $c_3 \in \R$ for which
\begin{align*}
    \regLoss(\theta) 
    &\geq \kappa_1 \Vert \dTheta \Vert + \kappa_2 \Vert \dBeta \Vert + c_3 \geq \kappa_3 \Vert \theta \Vert + c_3,
\end{align*}
and the last step follows by equivalence of norms in Euclidean spaces,
demonstrating that $\regLoss$ is bounded below and that
\(
    \liminf_{\Vert \theta \Vert \to \infty} \frac{\regLoss(\theta)}{\Vert \theta \Vert} = \kappa_3 > 0.
\)
Thus, $\regLoss$ is \textit{level-coercive} \cite[Definition 3.25]{rockafellar2009variational},
implying it has bounded level-sets \cite[Corollary 3.27]{rockafellar2009variational}.
As shown in \cref{cor: nonincreasing_sequence}, the sequence $\left\{ \regLoss(\theta^k) \right\}_{k \in \N}$ is nonincreasing.
Hence, the iterates $\left\{ \theta^k \right\}_{k \in \N}$ 
remain in the bounded level-set $\Omega \dfn \left\{ \theta \mid \regLoss(\theta) \leq \regLoss(\theta^0) \right\} \in \mathcal{X}_{\theta}$.
By continuity of $\regLoss$, $\Omega$ is closed \cite[Theorem 1.6(c)]{rockafellar2009variational}, and thus compact.

\subsection{Proof of \cref{lem: directional_derivative}}\label{prf: directional_derivative}

\subsection{Proof of \cref{lem: local_strong_convexity_of_kernels}}\label{prf: local_strong_convexity_of_kernels}
\revnew{2}{}{%
Following the definition \eqref{eq:bregman-function},
\vspace{-\baselineskip}%
\vspace{-\baselineskip}%
{\small%
\begin{align}
    \inner{\nabla \kernel[k](\theta) - \nabla \kernel[k](\tilde{\theta})}{\theta - \tilde{\theta}} \nonumber
    {=} & \inner{\nabla \Qi[1][k]{\dTheta} {-} \nabla \Qi[1][k]{\tilde{\dTheta}}}{\dTheta {-} \tilde{\dTheta}} \nonumber\\
      & {+}\inner{\nabla \Qi[2][k]{\dBeta} {-} \nabla \Qi[2][k]{\tilde{\dBeta}}}{\dBeta {-} \tilde{\dBeta}} \nonumber\\
      & {+} \inner{\nabla \reg(\theta) {-} \nabla \reg(\tilde{\theta})}{\theta {-} \tilde{\theta}}. \label{eq: kernel_diff_inner_product}
\end{align}}%
When \crefpart{asp: smoothness}{asp: reg_c2} applies,
there exists an $\varepsilon_\reg > 0$ such that
\[
    \inner{\nabla \reg(\theta) {-} \nabla \reg(\tilde{\theta})}{\theta {-} \tilde{\theta}} \geq \varepsilon_\reg \Vert\theta - \tilde{\theta}\Vert^2,
\]
as $\reg$ is locally strongly convex on any compact set $\Omega_\theta \subseteq \mathcal{X}_\theta$.
By applying this inequality to \eqref{eq: kernel_diff_inner_product},
\eqref{eq: local_strong_convexity_of_kernels} holds by following the convexity of $\Qi[1][k]{}$ and $\Qi[2][k]{}$,
and let $\varepsilon = \varepsilon_\reg$.
}

\revnew{2}{}{When \crefpart{asp: smoothness}{asp: sufficient_data} applies,
it suffices to show the first two terms on the right-hand side of \eqref{eq: kernel_diff_inner_product}
are each lower bounded by a quadratic term,
as $\reg$ is convex following \cref{asp: regularizers}.

We first examine the term involving $\dBeta$.
By \cref{lem: general_compactness},
all iterates $\{\theta^k\}_{k\in\N}$ remain in a compact set $\Omega_\theta \subseteq \mathcal{X}_\theta$.
Using the definition of weighted sum $\Qi[2][k]{}$ in \eqref{eq:subsys-surrogate},
and the lower bound on the probabilities $\Pi_{\modeSeq}(\thk) \geq \varepsilon_1$ from \cref{lem: bounded_weights},
we have
\begin{align*}
    & \inner{\nabla \Qi[2][k]{\dBeta} {-} \nabla \Qi[2][k]{\tilde{\dBeta}}}{\dBeta {-} \tilde{\dBeta}}\\
    \geq &\; \varepsilon_1 \textstyle\sum_{\modeSeq \in \Xi^{T+1}}
    \inner{\nabla \hat{\Phi}_{\modeSeq}^k(\dBeta) - \nabla \hat{\Phi}_{\modeSeq}^k(\tilde{\dBeta})}{\dBeta - \tilde{\dBeta}}.
\end{align*}
As $\hat{\Phi}_{\modeSeq}^k$ from \eqref{eq: def_elem_surrogate_subsys} is separable with respect to components $\beta_j$ for $j \in \Xi$,
we consider the partial derivative w.r.t. a specific $j\in\Xi$.
Following the chain rule to $\hat{\Phi}_{\modeSeq}^k$  yields
\vspace{-0.5\baselineskip}%
\[
    \partial_{\beta_j} \hat{\Phi}_{\modeSeq}^k(\dBeta) 
    {=} \tlsum_{\substack{t=0,\\ \xi_{t+1}=j}}^{T-1} \smashoperator{\tlsum_{i=1}^{n_g}} 
    \partial_i f(G_t(\beta^k_{j})) \jac{\!}_{it}^\top \nabla g_i(T_{it}(\beta_{j})),
\]
From \cref{lem: bounded_weights},
$\partial_i f(G_t(\beta_j^k)) \geq \varepsilon_2$ for all $k\in\N$ and $j\in\Xi$.
Consequently, the inner product can be lower bounded as
\vspace{-\baselineskip}%
\vspace{-\baselineskip}%
{\small%
\begin{align*}
    & \inner{\partial_{\beta_j} \hat{\Phi}_{\modeSeq}^k(\dBeta) - \partial_{\beta_j} \hat{\Phi}_{\modeSeq}^k(\tilde{\dBeta})}{\beta_j - \tilde{\beta}_j}\\
    \geq &\; \varepsilon_2 \tlsum_{\substack{t=0,\\ \xi_{t+1}=j}}^{T-1} \smashoperator{\tlsum_{i=1}^{n_g}}
    \inner{\nabla g_i(T_{it}(\beta_{j})) - \nabla g_i(T_{it}(\tilde{\beta}_{j}))}{T_{it}(\beta_j) - T_{it}(\tilde{\beta}_j)}
\end{align*}}%
Since the set of paths $\Xi^{T+1}$ is exhaustive,
for any $t \in \N_{[0, T-1]}$,
there exists a subset of paths where mode $j$ is active.
Therefore, the summation over $\modeSeq$ covers the full time horizon for every mode $j$:
\vspace{-\baselineskip}%
\vspace{-\baselineskip}%
{\small%
\begin{align}
    & \tlsum_{\modeSeq \in \Xi^{T+1}} 
    \inner{\partial_{\beta_j} \hat{\Phi}_{\modeSeq}^k(\dBeta) - \partial_{\beta_j} \hat{\Phi}_{\modeSeq}^k(\tilde{\dBeta})}{\beta_j - \tilde{\beta}_j}
    \nonumber\\
    \geq &\; \varepsilon_2 \tlsum_{t=0}^{T-1} \tlsum_{i=1}^{n_g}
    \inner{\nabla g_i(T_{it}(\beta_{j})) - \nabla g_i(T_{it}(\tilde{\beta}_{j}))}{T_{it}(\beta_j) - T_{it}(\tilde{\beta}_j)} \nonumber \\
    \geq &\; \varepsilon_2 \mu_j \Vert{\beta_j - \tilde{\beta}_j}\Vert^2, \label{eq: diff_grad_Phi_beta_j}
\end{align}}%
where the second inequality follows from \cref{lem: local_strong_convexity_on_a_compact_set} with $\mu_j > 0$.
Multiplying by $\varepsilon_1$ and summing \cref{eq: diff_grad_Phi_beta_j} over all $j \in \Xi$,
we obtain
\[
    \inner{\nabla \Qi[2][k]{\dBeta} {-} \nabla \Qi[2][k]{\tilde{\dBeta}}}{\dBeta {-} \tilde{\dBeta}} \geq \varepsilon_\beta \Vert \dBeta - \tilde{\dBeta} \Vert^2
\]
with $\varepsilon_\beta = \varepsilon_1\varepsilon_2 \min_{j \in \Xi} \mu_j > 0$.
}

\revnew{2}{}{%
The analysis for $\Qi[1][k]{}$ defined in \eqref{eq:switch-surrogate} takes analogous steps.
Bounding the probabilities $\Pi_{\modeSeq}(\thk)$ from \cref{lem: bounded_weights}
yields
\begin{align*}
    &\inner{\nabla \Qi[1][k]{\dTheta} {-} \nabla \Qi[1][k]{\tilde{\dTheta}}}{\dTheta {-} \tilde{\dTheta}} \\
    \geq &\; \varepsilon_1 \textstyle\sum_{\modeSeq \in \Xi^{T+1}}
    \inner{\nabla \Psi_{\modeSeq}(\dTheta) - \nabla \Psi_{\modeSeq}(\tilde{\dTheta})}{\dTheta - \tilde{\dTheta}}.
\end{align*}
As defined in \eqref{eq:Lt2-def},
$\Psi_{\modeSeq}$ is a sum of $\lse$ functions composed with $\Theta_{\xi_{t}}^\top z_t$.
By using the vec-trick of Kronecker product $\Theta_{\xi_{t}}^\top z_t = (I \otimes z_t^\top) \vec(\Theta_{\xi_t})$,
where $\vec$ flattens a matrix into a vector by stacking its columns and $\otimes$ denotes the Kronecker product,
we identify the jacobian of the linear mapping as $\jac{\!}_{\xi_t, t} = (I \otimes z_t^\top)$.
As remarked in \cref{rem: translation_invariance},
$x \mapsto \lse(x)$ is strictly convex with the constraint $x_d=0$,
and its Hessian is positive definite on a compact set.
Since the summation over $\modeSeq \in \Xi^{T+1}$ is exhaustive and covers the full time horizon for every mode, 
applying \cref{lem: local_strong_convexity_on_a_compact_set} yields a constant $\mu_\theta > 0$ such that
\[
    \inner{\nabla \Qi[1][k]{\dTheta} - \nabla \Qi[1][k]{\tilde{\dTheta}}}{\dTheta - \tilde{\dTheta}}
     \geq \mu_\theta \Vert \dTheta - \tilde{\dTheta} \Vert^2.
\]
}

\subsection{Proof of \cref{lem: gradient_like_descent_sequence}}\label{prf: gradient_like_descent_sequence}

\subsection{Proof of \cref{prop: forward_backward}}\label{prf: forward_backward}

\input{contents/revisionReferenceStyle.tex}

\bibliographystyle{plain}        %
\bibliography{reference}           %

\end{document}